\pdfoutput=1
\documentclass[12pt]{amsproc}

\newtheorem{theorem}{Theorem}[section]
\newtheorem{lemma}[theorem]{Lemma}
\newtheorem{proposition}[theorem]{Proposition}
\newtheorem{corollary}[theorem]{Corollary}

\theoremstyle{definition}

\theoremstyle{remark}
\newtheorem{remark}[theorem]{Remark}

\numberwithin{equation}{section}
\numberwithin{table}{section}
\numberwithin{figure}{section} %

\usepackage[margin=20mm,includeheadfoot,headheight=14pt]{geometry}
\AddToHook{cmd/@settitle/before}{\global\topskip=24pt\relax}
\usepackage{xcolor}
\usepackage{graphicx}
\usepackage{placeins}
\usepackage{amssymb}
\usepackage{mathrsfs}
\usepackage{tikz}
\usetikzlibrary{quotes,angles,arrows.meta}
\usepackage{subcaption}
\usepackage[hidelinks,bookmarksnumbered,bookmarksopen,bookmarksdepth=2,
  pdfpagemode=UseOutlines]{hyperref}

\renewcommand{\epsilon}{\varepsilon}

\begin{document}

\title[\texorpdfstring{Staggered Discontinuous Galerkin Neural Networks}{Staggered Discontinuous Galerkin Neural Networks}]{\texorpdfstring{A staggered discontinuous Galerkin neural network method with a computable error decomposition}{A staggered discontinuous Galerkin neural network method with a computable error decomposition}}

\author[G.-M. Gie]{Gung-Min Gie}  %
\address[Gung-Min Gie]{Department of Mathematics, University of Louisville,
Louisville, Kentucky, USA}
\email{gungmin.gie@louisville.edu}

\author[C.-Y. Jung]{Chang-Yeol Jung}
\address[Chang-Yeol Jung]{Department of Mathematical Sciences, Ulsan National Institute of
Science and Technology (UNIST), Ulsan 44919, Republic of Korea}
\email{cjung@unist.ac.kr}

\author[B. Kwon]{Bongsuk Kwon}
\address[Bongsuk Kwon]{Department of Mathematical Sciences, Ulsan National Institute of
Science and Technology (UNIST), Ulsan 44919, Republic of Korea}
\email{bkwon@unist.ac.kr}

\author[K. Kwon]{Kiwoong Kwon}
\address[Kiwoong Kwon]{Department of Mathematical Sciences, Ulsan National Institute of
Science and Technology (UNIST), Ulsan 44919, Republic of Korea}
\email{kwkwon@unist.ac.kr}
\thanks{Corresponding author: Kiwoong Kwon (kwkwon@unist.ac.kr).}

\begin{abstract}
\par
We develop a staggered discontinuous Galerkin neural network (SDGNN) method for scalar elliptic
problems with a computable error decomposition. The staggered discontinuous
Galerkin (SDG) method requires only local continuity of its trial functions and gives a square,
inf-sup stable system without stabilization parameters; we carry this structure over to neural
networks. Networks on the patches of a staggered mesh satisfy the SDG continuity
conditions for all parameter values, and the loss, the squared dual norm of the SDG weak
residual, is computed without the SDG finite element approximation or any solution data, so that the
training is unsupervised. By the inf-sup stability, the loss is equivalent,
with computable constants, to the square of the training error,
the distance between
the SDG projection of the network and the SDG finite element approximation, so that the full
error is bounded by the discretization error, a multiple of the square root of the loss, and
the network projection error. For a fixed network on uniform square meshes, this error decays
exponentially in the polynomial degree. The weak residual requires no derivatives of the
networks and no penalty or stabilization weights, and the scalar trial functions
may jump across the midlines of the primal cells. Numerical experiments show that the square
root of the loss tracks the error during training, that prescribed jumps
across a square and a circular interface are captured consistently over the initializations, and that the
method is more accurate than physics-informed neural networks at the same parameter counts.
\par
\end{abstract}


\maketitle
\makeatletter\@setaddresses\global\let\addresses\@empty\makeatother\enlargethispage{2pt}

\section{Introduction}
\label{sec:introduction}

\par
We study neural network approximations for scalar elliptic problems and
decompose their error into computable terms, one of which is equivalent to the
square root of the training loss. The idea of representing
the solution of a differential equation by a neural network and training it on
a loss derived from the equation goes back to the 1990s~\cite{1998Lagaris};
physics-informed neural networks (PINNs)~\cite{2019Raissi} and the deep
Galerkin method~\cite{2018Sirignano} minimize the squared strong-form residual
at sampled points, and the deep Ritz method~\cite{2018E} minimizes the energy
functional. For the residual losses, however, the original works do not
relate the value of the loss to the size of the error: the PINN work gives no error
estimate~\cite{2019Raissi}, and the deep Galerkin method is shown to converge
only as the loss tends to zero, without a
rate~\cite[Section~7]{2018Sirignano}. For scalar elliptic problems, we
carry the structure of a stable finite element method over to the neural
network trial functions. The loss is then equivalent, with computable
constants, to the square of the training error, the distance between the
finite element projection of the network and the finite element approximation.
The full error is bounded by the sum of the training error and two terms
determined by the finite element structure: the discretization error and the
distance between the network and its finite element projection.
\par

\par
A loss that controls the error has been sought in two
directions~\cite{2025Ernsta,2026Bachmayr}. The first direction changes the loss: it tests the equation
against test functions or minimizes a least-squares functional.
hp-VPINNs~\cite{2021Kharazmi} sum the squared weak residuals against local
polynomial test functions and add a boundary penalty, and weak adversarial
networks~\cite{2019Zang} approximate the dual norm of the weak residual by a maximization over a
second network that plays the test function. Robust
VPINNs~\cite{2024Rojas} measure the weak residual in the dual norm of a discrete test
space, and the deep Fourier residual method~\cite{2022Taylor} in the $H^{-1}$ norm, computed
on rectangular domains from a truncated sine and cosine expansion. These losses bound the error
up to the component that the test space or the truncation does
not see. For robust VPINNs, this component is controlled only under a local Fortin-type
assumption, with a remainder that is not estimated, or estimated by an adaptive
enrichment of the test space under a saturation
assumption~\cite{2026Udomworarat}. With the optimal test norm of discontinuous
Petrov--Galerkin methods, the dual norm of the residual equals the error, but
the norm is defined through the operator itself~\cite{2011Zitelli}. Deep
least-squares methods for first-order systems~\cite{2020Cai} minimize a
functional equivalent to the squared error when the boundary data are measured
in the $H^{1/2}$ and $H^{-1/2}$ norms, but not for the discrete functional, in which
these norms are replaced by weighted $L^2$ norms; for parametric problems, losses proportional to the squared
error are obtained from first-order least-squares and ultra-weak discontinuous
Petrov--Galerkin formulations with piecewise polynomial trial functions whose
coefficients are given by networks~\cite{2026Bachmayr}. Galerkin neural
networks~\cite{2021Ainsworth} add network basis functions one at a time and
estimate the error a posteriori.
\par

\par
The second direction changes the trial functions: it uses the structure of a finite element space.
In~\cite{2022Berrone} and in finite element interpolated neural
networks~\cite{2024Badia,2025Badia}, the network is interpolated onto a finite
element space before the residual is computed, and compatible
versions~\cite{2025Badiaa} preserve the discrete de Rham complex; the analysis
of~\cite{2022Berrone} concerns the interpolated function, and that
of~\cite{2025Badia} bounds the error of the network by the finite element
error, a training error, and the part of the network that the interpolation
does not see, which the loss does not control. When networks are placed element
by element, the coupling between elements is imposed through interior penalty
terms or collocation equations for the jumps~\cite{2024Sun}, through
additional trace unknowns~\cite{2024Dang}, or by enriching the local basis of
a discontinuous Galerkin scheme~\cite{2024Franck}.
\par

\par
In both directions, when
the network itself is the trial function, the part of its error that the loss
does not see is left out of the analysis or bounded under assumptions on the
network~\cite{2024Rojas,2022Taylor,2022Berrone,2021Ainsworth}; a posteriori estimates for
a given network bound its error from above by computable quantities,
independently of the training~\cite{2025Fuhrer,2025Ernsta,2025Hillebrecht}. In the method of this article, the network itself is
the trial function, and this part is the distance between the network and its finite element projection, which we
call the network projection error.
\par

The finite element method we use is the staggered discontinuous Galerkin (SDG)
method, which was introduced for wave propagation~\cite{2006Chung,2009Chung}
and has been analyzed for elliptic problems on rectangular
grids~\cite{KJN21} and for convection-dominated diffusion~\cite{ZP19}. Three
properties of this method allow us to carry its structure over to neural network trial functions.
 First, the SDG spaces require only local continuity, and each
continuity condition concerns one vertex-centered patch or one half of a primal
cell; we construct networks on these patches that satisfy the conditions for
all values of the network parameters, without penalty terms, interpolation, or additional unknowns.
Second, the first-order weak form acts on a trial pair only through its values
and traces, so the weak residual requires no derivatives of the networks.
Third, the SDG trial and test spaces have the same dimension, and the resulting
square system is inf-sup stable without stabilization parameters, also with an
upwind treatment of convection. In addition, the scalar trial
functions may jump across the midlines of the primal cells and can therefore
represent a prescribed jump across an interface along these midlines. Dirichlet
data enter the residual only through integrals over $\partial\Omega$, and the
trace is imposed by moments as in the staggered hybridization
of~\cite{2017Chung}, so the loss has no separate boundary weight. We call the resulting method the staggered discontinuous
Galerkin neural network (SDGNN) method.

\par
The main results of this article are the following.
\begin{enumerate}
\item \emph{A loss equivalent to the training error.} The SDGNN loss $J_h$, the squared dual
  norm of the SDG weak residual on the SDG test space, is equivalent to the
  square of the training error, the distance between the SDG projection of the network and the SDG
  finite element approximation (Theorem~\ref{prop:scalar-loss-identity}). The equivalence follows
  from the discrete inf-sup stability of the SDG method, as for discrete trial
  functions~\cite[Theorem~2.22]{2004Ern}, \cite[Theorem~2.1]{2014Carstensen}. With the optimal test
  norm, the loss would equal the squared training error, but each evaluation would require a solve
  with the SDG system (Remark~\ref{rem:trial-inner-product}). We use the Gram matrix $M_Y$ of the
  test space instead, as robust VPINNs do~\cite{2024Rojas}, so that the training is unsupervised:
  each evaluation of $J_h$ requires one solve with this block diagonal matrix and uses neither the
  SDG finite element approximation nor any solution data. For robust VPINNs, the loss is
  equivalent to the error modulo the trial functions that the test space does not
  see~\cite[Theorem~2]{2024Rojas}. Here the SDG system is square and inf-sup stable, so the loss
  sees the difference of two functions in the SDG space, and the constants of the equivalence are
  the extreme singular values $\sigma_{\min}$ and $\sigma_{\max}$ of the SDG matrix scaled by the
  trial and test Gram matrices. We compute $\sigma_{\min}$ from a generalized eigenvalue problem
  and, on uniform square meshes, bound it below independently of the mesh size
  (Corollary~\ref{cor:singular-value-bounds}).
\item \emph{An error decomposition with computable terms.} As for finite element
  interpolated neural networks~\cite[Proposition~3.13]{2025Badia}, the full error is
  bounded by the discretization error $D_h$, the term
  $\sigma_{\min}^{-1}J_h^{1/2}$, and the network projection error $P_h$, the
  distance between the network and its SDG projection
  (Theorem~\ref{prop:scalar-error-estimate}). The bound follows from
  Theorem~\ref{prop:scalar-loss-identity} and the triangle inequality, as in the a priori error
  analysis of finite element methods (cf.\ the proof of~\cite[Theorem~2.24]{2004Ern}). Here the
  middle term is computed from the loss of the trained network, without an assumption on the
  training such as the quasi-minimization assumed in~\cite{2025Badia}, and $P_h$ is
  estimated in Section~\ref{sec:projection-error-estimate}. For a trained network, the last
  two terms are computable, and $D_h$ is computable when the exact solution is
  known.
\item \emph{Degree-explicit estimates.} On uniform square meshes, we make the discrete
  inf-sup constants of the SDG method explicit in the polynomial degree $k$
  (Lemma~\ref{lem:sdg-b-infsup} and Proposition~\ref{prop:weak-form-infsup}),
  where~\cite[Theorem~1]{KJN21} gives a constant independent of the mesh size for the form $B$.
  The stability factor of the SDG projection then grows at most like $k^6$
  (Corollary~\ref{cor:sdg-projection-approximation}). With the Chebyshev approximation of
  analytic functions~\cite[Theorems~8.1 and~8.2]{T13}, the analyticity of the tanh units gives
  $P_h\le Ck^6R_*^{-k}$ for a fixed network, with $R_*>1$ computed from the weights
  (Theorem~\ref{thm:fixed-network-degree-convergence}). For a fixed network, the part of the
  error that the loss does not see thus decays exponentially as the degree increases, whereas for
  finite element interpolated neural networks this part is reduced by adapting the
  mesh~\cite{2025Badia}. The same analyticity gives the quadrature orders
  (Appendix~\ref{subsec:weight-quadrature}), and earlier adaptive
  quadratures for neural network solvers adapt the integration cells
  instead~\cite{2020Cai,2022Rivera,2024Magueresse}.
\item \emph{Structure instead of penalties.} The method requires no
  derivatives of the networks and no penalty or stabilization weights, and it
  represents a prescribed jump across the midlines of the primal cells
  without penalty terms or additional unknowns.
\item \emph{Numerical evidence.} The square root of the loss tracks the error during the training
  (Section~\ref{subsec:numerical-loss-error}),
  prescribed jumps across a square and a circular interface are captured
  consistently over the initializations, and the scalar errors are smaller than
  those of PINN and Mixed-PINN at the same parameter counts
  (Section~\ref{subsec:numerical-pinn}); a comparison with global networks
  trained on the SDGNN loss shows the effect of the patch structure.
\end{enumerate}
\par

\par
The scope of this article is as follows. The SDG finite element approximation
on the same mesh is more accurate than the SDGNN approximation and, for a
single linear problem, cheaper; we do not propose the method as its
replacement, but as a neural network trial function whose error has a computable
decomposition.
On smooth problems, a single global network trained on the SDGNN loss attains
comparable errors at the same parameter count; the
advantage of the SDG structure appears on the internal layer, where the patch
networks are at least as accurate as the global structures we tested, and on
the interface, where a continuous global network cannot represent the jump. The analysis
is carried out on rectangular meshes, with the estimates of
Section~\ref{sec:projection-error-estimate} on uniform square meshes, and a prescribed jump requires a
mesh fitted to the interface, which is the standard assumption of
discontinuous Galerkin methods.
\par

Section~\ref{sec:sdg-method} reviews the SDG finite element method with
Dirichlet data imposed through boundary moments, and
Section~\ref{sec:neural-trial-functions} constructs the neural network trial
functions. Section~\ref{sec:weak-residual-loss} defines the SDGNN loss
and proves the error decomposition.
Section~\ref{sec:projection-error-estimate}
bounds $\sigma_{\min}$, $D_h$, and $P_h$ on uniform square meshes, Section~\ref{sec:numerical-experiments} presents
the numerical experiments, and Section~\ref{sec:conclusions} concludes the paper.
Appendix~\ref{app:degree-b-proof} proves Lemma~\ref{lem:sdg-b-infsup},
Appendix~\ref{app:boundary-trace-extension} constructs the discrete extension
of boundary traces, Appendix~\ref{app:chebyshev} proves
Theorem~\ref{thm:fixed-network-degree-convergence}, and Appendix~\ref{subsec:weight-quadrature} describes the quadrature.

\section{The SDG finite element method}
\label{sec:sdg-method}

We review the SDG finite element method on which the SDGNN method is
built: the model problem, the staggered mesh and function spaces, and the SDG finite
element approximation, in which Dirichlet data enter through boundary moments.
The staggered mesh, the spaces, and the bilinear forms $B$ and $B^*$ are those
of~\cite{KJN21} on rectangular grids, the transport term is treated by upwinding as
in~\cite{ZP19}, and the boundary moments follow the staggered hybridization
of~\cite{2017Chung}. These references do not cover the present combination, so we show in
Section~\ref{sec:sdg-formulation} that the approximation exists and is unique.

\subsection{Model problem}
\label{sec:scalar-elliptic-problem}

Let $\Omega\subset\mathbb R^2$ be a bounded rectangular domain, and let
$f\in L^2(\Omega)$ and $g\in H^{1/2}(\partial\Omega)$. We write $\mathbf x=(x,y)$ for a spatial point and
use $x_1=x$ and $x_2=y$ in formulas with coordinate indices.
We assume that $a,c\in L^\infty(\Omega)$,
$\mathbf b\in W^{1,\infty}(\Omega)^2$, and that there exists a constant
$a_0>0$ such that
\begin{equation}
\label{eq:coefficient-assumptions}
a(\mathbf x)\ge a_0,
\qquad
c(\mathbf x)-\frac12\nabla\cdot\mathbf b(\mathbf x)\ge0
\qquad
\text{for almost every }\mathbf x\in\Omega.
\end{equation}
We consider the scalar elliptic problem
\[
\left\{
\begin{aligned}
-\nabla\cdot(a\nabla u+\mathbf b u)+cu&=f
&&\text{in }\Omega,\\
u&=g
&&\text{on }\partial\Omega.
\end{aligned}
\right.
\]
The diffusive variable $\mathbf p$ and the reciprocal diffusion coefficient
$\alpha$ are defined by
\[
\mathbf p:=a\nabla u,
\qquad
\alpha:=a^{-1}.
\]
With these definitions, the model problem can be written equivalently as
the first-order system
\[
\left\{
\begin{aligned}
\alpha\mathbf p-\nabla u&=\mathbf0
&&\text{in }\Omega,\\
-\nabla\cdot\mathbf p-\nabla\cdot(\mathbf b u)+cu&=f
&&\text{in }\Omega,\\
u&=g
&&\text{on }\partial\Omega.
\end{aligned}
\right.
\]

The weak solution, which exists and is unique by the Lax--Milgram lemma~\cite[Lemma~1.4]{DE12}, is the function $u\in H^1(\Omega)$ with $u=g$ on
$\partial\Omega$ in the sense of traces that satisfies
\begin{equation}
\label{eq:weak-formulation}
\int_\Omega(a\nabla u+\mathbf b u)\cdot\nabla v\,d\mathbf x
+\int_\Omega cuv\,d\mathbf x
=\int_\Omega fv\,d\mathbf x
\qquad\text{for every }v\in H_0^1(\Omega).
\end{equation}
The exact solution is the pair
$(u,\mathbf p)$ with $\mathbf p=a\nabla u$.

\subsection{Staggered mesh and function spaces}
\label{sec:staggered-framework}

We partition $\Omega$ into finitely many open axis-aligned rectangles, called primal cells, such
that for any two distinct primal cells $K$ and $K'$, the intersection
$\overline K\cap\overline{K'}$ is either empty, a common vertex, or a common
edge. Each primal cell $K$ has center $\nu_K$ and is divided into four
open subrectangles by joining $\nu_K$ to the midpoint of each edge. These
subrectangles form the refined partition $\mathcal T$.

The vertices of the primal cells, including those on $\partial\Omega$, are
called primal vertices. The set $\mathcal F$ of all edges of $\mathcal T$ is divided
into $\mathcal F_u$, the edges contained in the boundary of a primal cell, and
$\mathcal F_p:=\mathcal F\setminus\mathcal F_u$, the edges inside primal cells. An
edge is called interior if it is not contained in $\partial\Omega$. Every edge in
$\mathcal F_p$ is interior; $\mathcal F_u^0$ denotes the interior edges in
$\mathcal F_u$, and $\mathcal F_\partial:=\mathcal F_u\setminus\mathcal F_u^0$
the boundary edges.

For a primal vertex $\mathbf x$, the vertex-centered patch of $\mathbf x$ is the
interior of the union of the closures of the subrectangles that have $\mathbf x$ as a vertex.
We use three families of patches. The family $\mathcal P_u$ consists of the
vertex-centered patches, and the families $\mathcal P_{p_x}$ and
$\mathcal P_{p_y}$ consist of the open halves of primal cells formed by the two
subrectangles sharing a vertical and a horizontal edge in $\mathcal F_p$,
respectively. Each subrectangle $\tau\in\mathcal T$ has exactly one primal
vertex among its vertices and lies in exactly one primal cell, so for
$s\in\{u,p_x,p_y\}$ it lies in exactly one patch of $\mathcal P_s$, which we
denote by $P_s(\tau)$. The primal cells and the patches of $\mathcal P_u$
correspond to the two overlapping families of rectangles of the SDG method
in~\cite[Section~2.1]{KJN21}.

\begin{figure}[t]
\centering
\definecolor{sdgCell}{RGB}{120,120,120}
\definecolor{sdgOverlap}{RGB}{137,104,166}
\definecolor{scalar}{RGB}{55,105,160}
\definecolor{px}{RGB}{193,112,38}
\definecolor{py}{RGB}{47,139,115}

\begin{tikzpicture}[
  every node/.style={font=\scriptsize},
  fu edge/.style={draw=black, line width=0.85pt},
  fp edge/.style={draw=black!65, densely dashed, line width=0.65pt},
  point/.style={circle, fill=black, inner sep=1.25pt},
  midpoint/.style={circle, draw=black, fill=white, line width=0.5pt, inner sep=1.0pt},
  label box/.style={fill=white, fill opacity=0.88, text opacity=1, inner sep=1.2pt},
]
  \begin{scope}[x=1.05cm, y=1.05cm]
    \fill[sdgCell!12] (0,0) rectangle (3,3);
    \fill[sdgOverlap!32] (0,0) rectangle (1.5,1.5);

    \draw[fp edge] (1.5,0) -- (1.5,3);
    \draw[fp edge] (0,1.5) -- (3,1.5);
    \draw[fu edge] (0,0) rectangle (3,3);

    \foreach \p in {(0,0),(3,0),(0,3),(3,3)} {
      \node[point] at \p {};
    }
    \foreach \p in {(1.5,0),(3,1.5),(1.5,3),(0,1.5)} {
      \node[midpoint] at \p {};
    }
    \node[point, fill=sdgCell] at (1.5,1.5) {};

    \node[label box, anchor=north west] at (0.08,2.92) {$K$};
    \node[label box, anchor=south west] at (1.57,1.55) {$\nu_K$};
    \node[label box] at (0.75,0.75) {$\tau\in\mathcal T$};
    \node[label box, anchor=north east] at (-0.05,-0.05) {$\mathbf x$};

    \node[font=\small, anchor=north] at (1.5,-1.08)
      {\textbf{(a)} Primal cell refinement};
  \end{scope}

  \foreach \k/\col/\lab in {0/scalar/{P_u(\tau)},1/px/{P_{p_x}(\tau)},2/py/{P_{p_y}(\tau)}} {
    \begin{scope}[shift={(4.55cm+\k*2.35cm,0.35cm)}, x=0.5cm, y=0.5cm]
      \fill[sdgCell!12] (2,2) rectangle (4,4);
      \ifnum\k=0 \fill[\col!18] (1,1) rectangle (3,3);\fi
      \ifnum\k=1 \fill[\col!18] (2,2) rectangle (4,3);\fi
      \ifnum\k=2 \fill[\col!18] (2,2) rectangle (3,4);\fi
      \fill[\col!40] (2,2) rectangle (3,3);
      \foreach \cx in {1,3} {
        \foreach \cy in {1,3} {
          \draw[fp edge] (\cx-1,\cy) -- (\cx+1,\cy);
          \draw[fp edge] (\cx,\cy-1) -- (\cx,\cy+1);
        }
      }
      \foreach \g in {0,2,4} {
        \draw[fu edge] (\g,0) -- (\g,4);
        \draw[fu edge] (0,\g) -- (4,\g);
      }
      \ifnum\k=0
        \draw[draw=\col, line width=2.0pt] (2,1) -- (2,3);
        \draw[draw=\col, line width=2.0pt] (1,2) -- (3,2);
        \node[circle, fill=\col, inner sep=1.1pt] at (2,2) {};
        \node[label box, text=\col, anchor=north east, inner sep=0.6pt] at (1.95,1.95) {$\mathbf x$};
      \fi
      \ifnum\k=1 \draw[draw=\col, line width=2.0pt] (3,2) -- (3,3);\fi
      \ifnum\k=2 \draw[draw=\col, line width=2.0pt] (2,3) -- (3,3);\fi
      \node[label box, inner sep=0.6pt] at (2.5,2.5) {$\tau$};
      \node[inner sep=0.6pt] at (3.55,3.55) {$K$};
      \node[text=\col, anchor=north] at (2,-0.12) {$\lab$};
    \end{scope}
  }
  \node[font=\small, anchor=north] at (8.0cm,-1.08cm)
    {\textbf{(b)} Patches containing $\tau$};

  \begin{scope}[shift={(5.75cm,-0.58cm)}]
    \draw[fu edge] (-1.47,0) -- (-0.99,0);
    \node[anchor=west] at (-0.89,0) {$e\in\mathcal F_u$};
    \draw[fp edge] (0.48,0) -- (0.96,0);
    \node[anchor=west] at (1.06,0) {$e\in\mathcal F_p$};
  \end{scope}
\end{tikzpicture}
\caption{Rectangular staggered refinement and the patches containing a
subrectangle. Panel~(a) depicts a primal cell $K$ with a primal vertex
$\mathbf x$, its four subrectangles in $\mathcal T$, and the edge families
$\mathcal F_u$ (solid) and $\mathcal F_p$ (dashed). Panel~(b) depicts, for a
subrectangle $\tau$ of a primal cell $K$, the patches $P_u(\tau)$,
$P_{p_x}(\tau)$, and $P_{p_y}(\tau)$ that contain it, with the edges inside each
patch drawn thick; across these edges, the scalar, the first component, and the
second component of the vector field are continuous in the sense of traces, respectively.}
\label{fig:staggered-mesh-patches}
\end{figure}
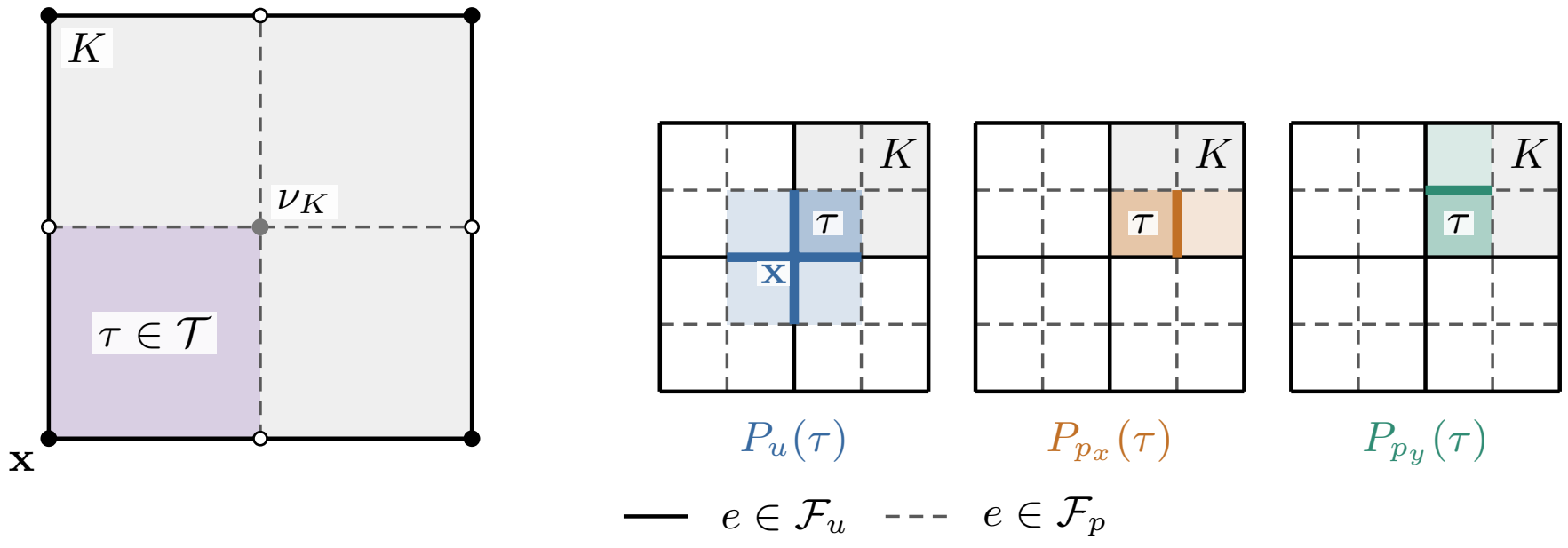

Figure~\ref{fig:staggered-mesh-patches} illustrates the refined partition
$\mathcal T$, the edge sets $\mathcal F_u$ and $\mathcal F_p$, and the three
patches that contain a subrectangle.

Each interior edge $e\in\mathcal F_u^0\cup\mathcal F_p$ is assigned a fixed unit normal
$\mathbf n_e$. The adjacent subrectangles $\tau_e^+$ and $\tau_e^-$ are
labeled so that $\mathbf n_e$ points from $\tau_e^+$ to $\tau_e^-$. For a
scalar function $v$ and a vector field $\mathbf z$ with well-defined traces
on $e$, the superscripts $+$ and $-$ indicate traces from these two
subrectangles. The scalar jump and the jump in the normal component are
defined by
\[
[v]_e:=v^+-v^-,
\qquad
[\mathbf z\cdot\mathbf n_e]_e
:=
(\mathbf z^+-\mathbf z^-)\cdot\mathbf n_e.
\]

The scalar space $\mathbb V$ and the vector space $\mathbb W$ are defined by
regularity on these patches:
\[
\begin{aligned}
\mathbb V
&:=
\bigl\{
v\in L^2(\Omega):
v|_P\in H^1(P)
\text{ for every }P\in\mathcal P_u
\bigr\},
\\
\mathbb W
&:=
\bigl\{
\mathbf z=(z_x,z_y)\in L^2(\Omega)^2:
z_x|_P\in H^1(P)\text{ for every }P\in\mathcal P_{p_x},
\\
&\hphantom{{}:=\bigl\{\mathbf z=(z_x,z_y)\in L^2(\Omega)^2:{}}
z_y|_P\in H^1(P)\text{ for every }P\in\mathcal P_{p_y}
\bigr\}.
\end{aligned}
\]
Each edge in $\mathcal F_u^0$ lies inside a patch of $\mathcal P_u$, and each
vertical or horizontal edge in $\mathcal F_p$ lies inside a patch of
$\mathcal P_{p_x}$ or $\mathcal P_{p_y}$, respectively. Hence functions in
$\mathbb V$ have continuous traces across $\mathcal F_u^0$, fields in
$\mathbb W$ have continuous normal components across $\mathcal F_p$, and both
have $L^2(e)$ traces on every edge. These are the continuity conditions of the
SDG spaces in~\cite{KJN21}.

For $m\in\mathbb N_0$ and a rectangle $D$, the space ${\mathbb Q}_m(D)$ consists
of polynomials of degree at most $m$ in each coordinate. For $k\in\mathbb N_0$, the
scalar space ${U_h^k}$, its subspace ${U_{h,0}^k}$ with zero boundary trace, and the vector
space ${W_h^k}$ are defined by
\[
\begin{aligned}
{U_h^k}
&:=
\left\{
v_h\in\mathbb V:
v_h|_\tau\in {\mathbb Q}_k(\tau)
\text{ for every }\tau\in\mathcal T
\right\},
\\
U_{h,0}^k
&:=
\left\{
v_h\in {U_h^k}:
v_h=0\text{ on }\partial\Omega
\right\},
\\
W_h^k
&:=
\left\{
\mathbf z_h\in\mathbb W:
\mathbf z_h|_\tau\in {\mathbb Q}_k(\tau)^2
\text{ for every }\tau\in\mathcal T
\right\}.
\end{aligned}
\]
By~\cite[Lemma~1.23]{DE12}, $U_h^k$ thus consists of the functions that lie in ${\mathbb Q}_k(\tau)$ on each $\tau\in\mathcal T$ and are continuous across $\mathcal F_u^0$, and $W_h^k$ of the fields that lie in ${\mathbb Q}_k(\tau)^2$ on each $\tau$ and have continuous normal components across $\mathcal F_p$.

For $(v,\mathbf z)\in\mathbb V\times\mathbb W$, $\nabla_hv$ and
$\nabla_h\cdot\mathbf z$ denote the gradient and divergence taken on each
subrectangle $\tau\in\mathcal T$.
For $v\in\mathbb V$ and $\mathbf z\in{\mathbb W}$,
the staggered bilinear forms are defined by
\[
\begin{aligned}
B(\mathbf z,v)
&:=
\sum_{\tau\in\mathcal T}
\int_\tau\mathbf z\cdot\nabla_hv\,d\mathbf x
-
\sum_{e\in\mathcal F_p}
\int_e(\mathbf z\cdot\mathbf n_e)[v]_e\,d\sigma,
\\
B^*(v,\mathbf z)
&:=
-
\sum_{\tau\in\mathcal T}
\int_\tau v\,\nabla_h\cdot\mathbf z\,d\mathbf x
+
\sum_{e\in\mathcal F_u^0}
\int_e v[\mathbf z\cdot\mathbf n_e]_e\,d\sigma.
\end{aligned}
\]
For these functions, integration by parts on each subrectangle gives
\[
B(\mathbf z,v)
=
B^*(v,\mathbf z)
+
\int_{\partial\Omega}v\,\mathbf z\cdot\mathbf n\,d\sigma,
\]
where $\mathbf n$ is the outward unit normal to $\partial\Omega$, and the
boundary integral is the sum over the boundary edges.
On an interior edge $e$, the contribution is $\int_e(v^+\mathbf z^+-v^-\mathbf z^-)\cdot\mathbf n_e\,d\sigma$, which equals $\int_e(\mathbf z\cdot\mathbf n_e)[v]_e\,d\sigma$ for $e\in\mathcal F_p$ by the normal continuity and $\int_e v[\mathbf z\cdot\mathbf n_e]_e\,d\sigma$ for $e\in\mathcal F_u^0$ by the scalar continuity. For
$(\mathbf z_h,v_h)\in {W_h^k}\times {U_{h,0}^k}$, the zero boundary trace of
$v_h$ yields
\begin{equation}
\label{eq:discrete-adjoint-identity}
B(\mathbf z_h,v_h)=B^*(v_h,\mathbf z_h).
\end{equation}

For each $e\in\mathcal F$, $h_e$ denotes the length of $e$. The staggered
broken $H^1$-seminorm on $\mathbb V$ is defined by
\[
\|v\|_{1,h}^2
:=
\sum_{\tau\in\mathcal T}
\|\nabla_hv\|_{L^2(\tau)}^2
+
\sum_{e\in\mathcal F_p}
h_e^{-1}\|[v]_e\|_{L^2(e)}^2.
\]
The mesh-dependent norm on $\mathbb W$ is defined by
\[
\|\mathbf z\|_{0,h}^2
:=\|\mathbf z\|_{L^2(\Omega)}^2
+\sum_{e\in\mathcal F_p}
h_e\|\mathbf z\cdot\mathbf n_e\|_{L^2(e)}^2.
\]

On $\mathbb V\times L^2(\Omega)^2$, we use the norm
\begin{equation}
\label{eq:full-error-norm}
\|(v,\mathbf p)\|_h^2
:=\|v\|_{1,h}^2
+\sum_{e\in\mathcal F_\partial}h_e^{-1}\|v\|_{L^2(e)}^2
+\|\mathbf p\|_{L^2(\Omega)}^2.
\end{equation}
On $\mathbb V\times{\mathbb W}$, we also use
\begin{equation}
\label{eq:projection-input-norm}
\begin{aligned}
\|(v,\mathbf p)\|_{*,h}^2
&:=\|v\|_{1,h}^2
+\sum_{e\in\mathcal F_\partial}h_e^{-1}\|v\|_{L^2(e)}^2
+\|\mathbf p\|_{0,h}^2\\
&=\|(v,\mathbf p)\|_h^2
+\sum_{e\in\mathcal F_p}h_e
\|\mathbf p\cdot\mathbf n_e\|_{L^2(e)}^2.
\end{aligned}
\end{equation}
The additional term controls the normal traces in $B$, and
$\|\cdot\|_h\le\|\cdot\|_{*,h}$ on $\mathbb V\times{\mathbb W}$.
For $v\in U_{h,0}^k$, the boundary term vanishes.

\subsection{SDG finite element approximation}
\label{sec:sdg-formulation}

The upwind trace of $v\in\mathbb V$ is defined pointwise on each
$e\in\mathcal F_p$ by the trace from the upstream side with
respect to $-\mathbf b$, the transport velocity of the term
$-\nabla\cdot(\mathbf b u)$:
\[
v^{\mathrm{up}}|_e
:=
\begin{cases}
v^+,&\mathbf b\cdot\mathbf n_e\le0,\\
v^-,&\mathbf b\cdot\mathbf n_e>0.
\end{cases}
\]
For $v\in\mathbb V$ and $w_h\in U_{h,0}^k$, the upwind transport bilinear
form is defined by
\[
C_h^{\mathrm{up}}(v,w_h)
:=
\sum_{\tau\in\mathcal T}
\int_\tau \mathbf b v\cdot\nabla_hw_h\,d\mathbf x
-
\sum_{e\in\mathcal F_p}
\int_e(\mathbf b\cdot\mathbf n_e)v^{\mathrm{up}}[w_h]_e\,d\sigma.
\]
This form discretizes the conservative transport term
$-\nabla\cdot(\mathbf b u)$; its edge sum involves only $\mathcal F_p$
because $w_h$ is continuous across $\mathcal F_u^0$ and vanishes on
$\partial\Omega$.

For $k\in\mathbb N_0$, the SDG trial and test spaces with zero boundary
trace are
\[
X_{h,0}^k:=U_{h,0}^k\times W_h^k,
\qquad
{Y_{h,0}^k}:=W_h^k\times U_{h,0}^k.
\]
We denote a trial pair by
$\Phi=(v,\mathbf p)\in\mathbb V\times{\mathbb W}$
and a test pair by $\Psi_h=(\mathbf z_h,w_h)\in {Y_{h,0}^k}$.
The weak form $\mathcal A$ is defined by
\[
\begin{aligned}
\mathcal A(\Phi,\Psi_h)
&:=
(\alpha\mathbf p,\mathbf z_h)_\Omega
-B^*(v,\mathbf z_h)
\\
&\quad+
B(\mathbf p,w_h)
+C_h^{\mathrm{up}}(v,w_h)
+(cv,w_h)_\Omega.
\end{aligned}
\]
The first two terms of $\mathcal A$ discretize the diffusive
constitutive relation, and the remaining terms discretize the equilibrium
equation with upwind transport. Here $(\cdot,\cdot)_S$ denotes the $L^2(S)$ inner product, where $S$ is $\Omega$, a subrectangle, or an edge.

\begin{proposition}[Discrete uniqueness]
\label{prop:scalar-discrete-uniqueness}
Suppose that $\Phi_h=(v_h,\mathbf p_h)\in {X_{h,0}^k}$ satisfies
\[
\mathcal A(\Phi_h,\Psi_h)=0
\qquad
\text{for every }\Psi_h\in {Y_{h,0}^k}.
\]
Then $v_h=0$ and $\mathbf p_h=\mathbf0$.
\end{proposition}

\begin{proof}
Taking $\Psi_h=(\mathbf p_h,v_h)$, using
\eqref{eq:discrete-adjoint-identity}, and integrating the transport term
by parts on each subrectangle gives, with the identity $\frac12(\mathbf b\cdot\mathbf n_e)[v_h^2]_e-(\mathbf b\cdot\mathbf n_e)v_h^{\mathrm{up}}[v_h]_e=\frac12|\mathbf b\cdot\mathbf n_e|\,[v_h]_e^2$ on each $e\in\mathcal F_p$, the continuity of $v_h$ across $\mathcal F_u^0$, and $v_h=0$ on $\partial\Omega$, the energy identity
\[
\begin{aligned}
0
&=
\mathcal A(\Phi_h,\Psi_h)
\\
&=
(\alpha\mathbf p_h,\mathbf p_h)_\Omega
+
\left(
\left(c-\frac12\nabla\cdot\mathbf b\right)v_h,
v_h
\right)_\Omega
\\
&\quad+
\frac12
\sum_{e\in\mathcal F_p}
\int_e
|\mathbf b\cdot\mathbf n_e|\,[v_h]_e^2\,d\sigma.
\end{aligned}
\]
By \eqref{eq:coefficient-assumptions}, each term is nonnegative, so
$(\alpha\mathbf p_h,\mathbf p_h)_\Omega=0$. Since $\alpha=a^{-1}>0$ almost everywhere in $\Omega$, it follows that
$\mathbf p_h=\mathbf0$.
With $\mathbf p_h=\mathbf0$, testing with $\Psi_h=(\mathbf z_h,0)$ and
using \eqref{eq:discrete-adjoint-identity} yields
\[
B(\mathbf z_h,v_h)=B^*(v_h,\mathbf z_h)=0
\qquad\text{for every }\mathbf z_h\in {W_h^k}.
\]
On each subrectangle $\tau$, the affine functions $\omega_{\tau,1}$ and
$\omega_{\tau,2}$ vanish on its vertical and horizontal edges in
$\mathcal F_p$, respectively, and equal one on the opposite edges
(cf.~\cite[proof of Theorem~1]{KJN21}).
The vector function $\mathbf z_h^{\nabla}$ is defined on each subrectangle by
\[
\mathbf z_h^{\nabla}|_\tau
:=
\bigl(\omega_{\tau,1}\partial_1v_h,\,
      \omega_{\tau,2}\partial_2v_h\bigr).
\]
Each component of $\mathbf z_h^{\nabla}|_\tau$ belongs to ${\mathbb Q}_k(\tau)$, and the
normal component of $\mathbf z_h^{\nabla}$ vanishes on every edge in $\mathcal F_p$. Thus,
$\mathbf z_h^{\nabla}\in {W_h^k}$, and testing with $\mathbf z_h=\mathbf z_h^{\nabla}$ gives
\[
0=B(\mathbf z_h^{\nabla},v_h)
=
\sum_{\tau\in\mathcal T}\int_\tau
\left(
\omega_{\tau,1}|\partial_1v_h|^2
+
\omega_{\tau,2}|\partial_2v_h|^2
\right)\,d\mathbf x.
\]
The positivity of both weights in the interior of each subrectangle
implies $\nabla_hv_h=0$, so $v_h$ is constant on each subrectangle.

For each $e\in\mathcal F_p$, the piecewise constant vector function
$\mathbf z_e$ equals $-[v_h]_e\mathbf n_e$ on the two adjacent
subrectangles and $\mathbf0$ elsewhere. Its normal component is continuous
across $e$ and vanishes on every other edge in $\mathcal F_p$, so
$\mathbf z_e\in {W_h^k}$. Since $\nabla_hv_h=0$, testing with
$\mathbf z_h=\mathbf z_e$ gives
\[
0=B(\mathbf z_e,v_h)
=\int_e[v_h]_e^2\,d\sigma.
\]
Thus, $[v_h]_e=0$ on every $e\in\mathcal F_p$. Since $v_h$ is already
continuous across $\mathcal F_u^0$, its constant values on adjacent
subrectangles agree. The connectedness of $\Omega$ implies that $v_h$ is
constant throughout $\Omega$, and the homogeneous boundary condition
gives $v_h=0$.
\end{proof}

For Dirichlet data $g$, we impose the boundary condition through
moments against the trace space
$U_h^k|_{\partial\Omega}:=\{v_h|_{\partial\Omega}:v_h\in U_h^k\}$.

For $\Psi_h=(\mathbf z_h,w_h)\in {Y_{h,0}^k}$, the load functional is given by
\begin{equation}
\label{eq:scalar-load-functional}
F_g(\Psi_h)
:=
(f,w_h)_\Omega
+
\int_{\partial\Omega}g\,\mathbf z_h\cdot\mathbf n\,d\sigma.
\end{equation}
We use the trial and test spaces $X_h^k:=U_h^k\times W_h^k$
and $Y_h^k:=W_h^k\times U_{h,0}^k\times U_h^k|_{\partial\Omega}$,
and we write a test triple as $(\Psi_h,\mu_h)$ with $\Psi_h\in Y_{h,0}^k$. For
$\Phi=(v,\mathbf p)\in\mathbb V\times{\mathbb W}$, the extended
weak form and load functional are
\[
\begin{aligned}
\mathcal A_\partial\bigl(\Phi,(\Psi_h,\mu_h)\bigr)
&:=\mathcal A(\Phi,\Psi_h)+\int_{\partial\Omega}v\,\mu_h\,d\sigma,\\
F_{g,\partial}(\Psi_h,\mu_h)
&:=F_g(\Psi_h)+\int_{\partial\Omega}g\,\mu_h\,d\sigma.
\end{aligned}
\]

The SDG finite element approximation with boundary data $g$ is the pair
$({u_{h,g}^k},{\mathbf p_{h,g}^k})\in X_h^k$ satisfying
\begin{equation}
\label{eq:sdg-boundary-problem}
\mathcal A_\partial\bigl(({u_{h,g}^k},{\mathbf p_{h,g}^k}),(\Psi_h,\mu_h)\bigr)
=
F_{g,\partial}(\Psi_h,\mu_h)
\qquad
\text{for every }(\Psi_h,\mu_h)\in Y_h^k.
\end{equation}
If $\Phi_h=(v_h,\mathbf p_h)\in X_h^k$ satisfies
\eqref{eq:sdg-boundary-problem} with $f=0$ and $g=0$, the moment equations give
$v_h|_{\partial\Omega}=0$, so $\Phi_h\in X_{h,0}^k$, and
Proposition~\ref{prop:scalar-discrete-uniqueness} gives $\Phi_h=0$. Since
$\dim X_h^k=\dim Y_h^k$, the approximation exists and is unique. The moment equations state that
$u_{h,g}^k|_{\partial\Omega}$ is the $L^2(\partial\Omega)$-orthogonal
projection of $g$ onto $U_h^k|_{\partial\Omega}$. Staggered hybridization
imposes the Dirichlet data through moments against polynomials of degree $k$ on
each boundary edge~\cite[eq.~(9) and Section~3.1]{2017Chung}; we use the trace
space of $U_h^k$, whose functions are continuous at the primal vertices on
$\partial\Omega$ and may jump at the midpoints of the primal boundary
edges, where edges of $\mathcal F_p$ meet $\partial\Omega$.

\section{Neural network trial functions with the SDG structure}
\label{sec:neural-trial-functions}

The scalar solution $u$ and the Cartesian components $p_x$ and $p_y$ of
$\mathbf p=a\nabla u$ are approximated by separate neural networks on the
patches of $\mathcal P_u$, $\mathcal P_{p_x}$, and $\mathcal P_{p_y}$ of
Section~\ref{sec:staggered-framework}, respectively.

\subsection{Neural networks on patches}
\label{subsec:subrectangle-local-interior-networks}

For a patch $P$ in $\mathcal P_u\cup\mathcal P_{p_x}\cup\mathcal P_{p_y}$, $x_{P,i}^{\min}$ and $x_{P,i}^{\max}$ denote
its lower and upper bounds in the $i$th coordinate. The affine map
$\boldsymbol\xi_P$ from $\overline P$ onto the reference square
$\widehat P=[-1,1]^2$ is given by
\[
\xi_{P,i}(\mathbf x)
=
2\frac{x_i-x_{P,i}^{\min}}
{x_{P,i}^{\max}-x_{P,i}^{\min}}
-1,
\qquad i=1,2,
\]
for $\mathbf x=(x_1,x_2)\in\overline P$. The same map provides the network input
$\boldsymbol\xi_P(\mathbf x)$ on every subrectangle in $P$.

The vector $\theta$ collects the parameters of all scalar networks,
and $\eta$ collects those of all networks for $p_x$ and $p_y$; $\Theta_u$
and $\Theta_{\mathbf p}$ are the Euclidean spaces of these parameter vectors.
For $\theta\in\Theta_u$ and $\eta\in\Theta_{\mathbf p}$, the
networks for $u$, $p_x$, and $p_y$ are denoted by
\[
N_{\theta,P}^u:\widehat P\to\mathbb R,
\qquad
N_{\eta,P}^{p_x}:\widehat P\to\mathbb R,
\qquad
N_{\eta,P}^{p_y}:\widehat P\to\mathbb R.
\]
In each map, $P$ denotes a patch for the component named in the superscript. Networks
on different patches or for different components have independent
parameters.

The networks have one hidden layer with local hidden width $r\in\mathbb N$
and hyperbolic tangent activation. Any network that is smooth on the closure of its patch gives
trial functions with the same continuity properties, and
Theorems~\ref{prop:scalar-loss-identity} and~\ref{prop:scalar-error-estimate} do not
depend on its form; the one-hidden-layer form is used in the quadrature of
Appendix~\ref{subsec:weight-quadrature} and in
Theorem~\ref{thm:fixed-network-degree-convergence}. For $\boldsymbol\xi\in\widehat P$,
the scalar network is given by
\[
N_{\theta,P}^{u}(\boldsymbol\xi)
=
c_P^{u}
+
\sum_{\ell=1}^{r}
a_{P,\ell}^{u}
\tanh\!\left(
(\mathbf w_{P,\ell}^{u})^{\mathsf T}\boldsymbol\xi
+
b_{P,\ell}^{u}
\right),
\]
and $N_{\eta,P}^{p_x}$ and $N_{\eta,P}^{p_y}$ have the same form with
the superscript $u$ replaced by $p_x$ and $p_y$, respectively.
For each component $s\in\{u,p_x,p_y\}$, the trainable parameters are the
scalars $c_P^s$, $a_{P,\ell}^s$, and $b_{P,\ell}^s$, and the vectors
$\mathbf w_{P,\ell}^s\in\mathbb R^2$. We call
$a_{P,\ell}^s\tanh\bigl((\mathbf w_{P,\ell}^s)^{\mathsf T}\boldsymbol\xi+b_{P,\ell}^s\bigr)$
the $\ell$th hidden unit of $N_{\theta,P}^s$ or $N_{\eta,P}^s$ and
$a_{P,\ell}^s$ its output weight.

\subsection{Neural network trial functions}
\label{subsec:neural-trial-functions}

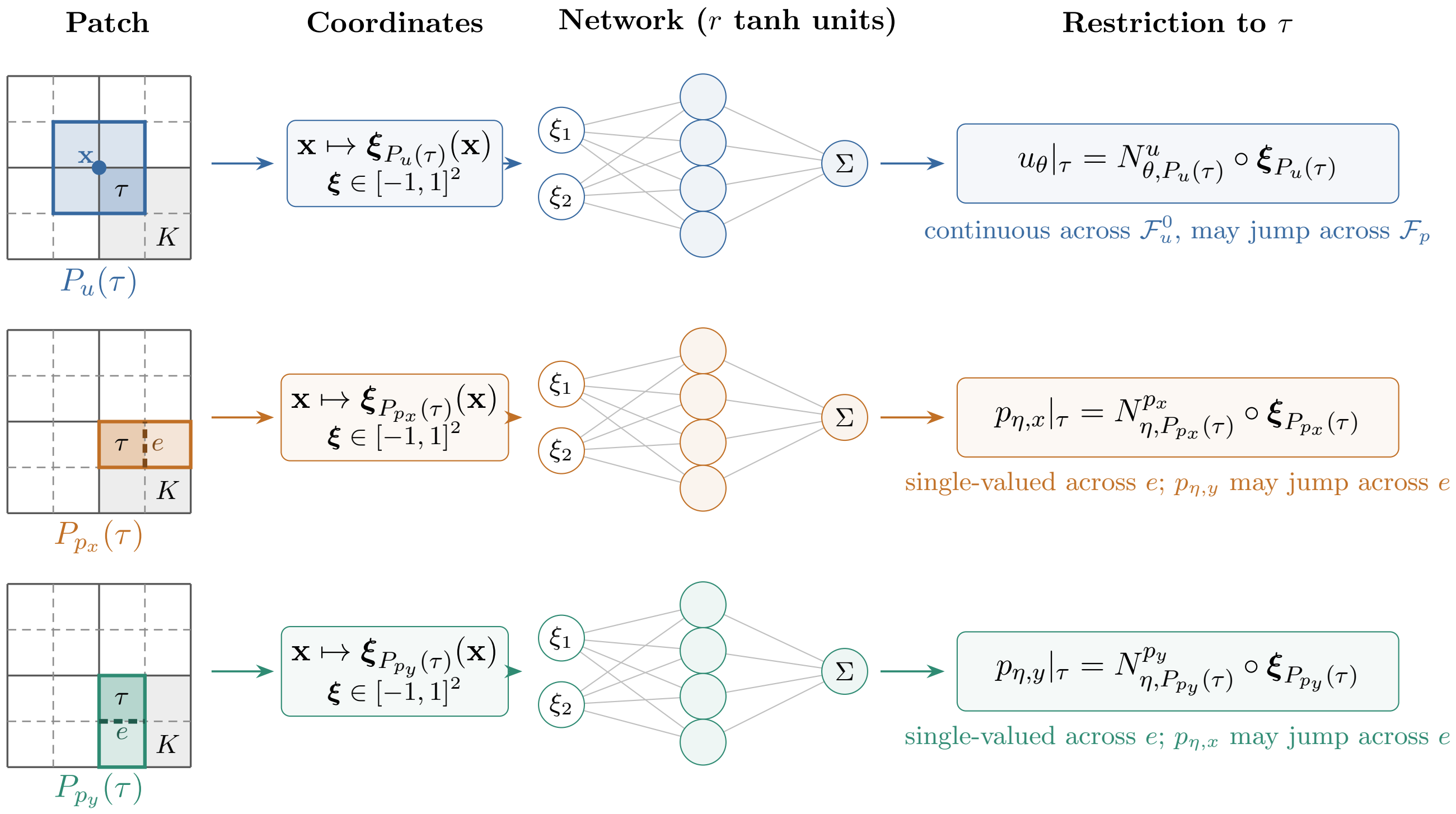
\begin{figure}[tbp]
\centering
\begingroup
\fontsize{10}{12}\selectfont
\def\SDGNNFigureOnly{}
\resizebox{\textwidth}{!}{\let\SDGNNFinish\relax
\definecolor{scalar}{RGB}{55,105,160}
\definecolor{px}{RGB}{193,112,38}
\definecolor{py}{RGB}{47,139,115}
\centering
\begin{tikzpicture}[x=1cm,y=1cm,>=Stealth,
  every node/.style={font=\small},
  neuron/.style={circle,draw,fill=white,minimum size=5.5mm,inner sep=0pt},
  wire/.style={draw=black!25,line width=.4pt},
  flow/.style={->,line width=.7pt},
  fu/.style={draw=black!65,line width=.65pt},
  fp/.style={draw=black!45,densely dashed,line width=.5pt}]
\path[use as bounding box] (-.3,-8.85) rectangle (17.1,1.0);
\node[font=\bfseries] at (1.2,.65) {Patch};
\node[font=\bfseries] at (4.65,.65) {Coordinates};
\node[font=\bfseries] at (8.65,.65) {Network ($r$ tanh units)};
\node[font=\bfseries] at (14.05,.65) {Restriction to $\tau$};

\foreach \row/\col/\comp in {0/scalar/u,1/px/p_x,2/py/p_y}{
  \pgfmathsetmacro{\yy}{-3.05*\row}
  \begin{scope}[shift={(0,\yy)}]
    \begin{scope}[x=.55cm,y=.55cm]
      \fill[black!7] (2,-4) rectangle (4,-2);
      \ifnum\row=0 \fill[\col!18] (1,-3) rectangle (3,-1);\fi
      \ifnum\row=1 \fill[\col!18] (2,-3) rectangle (4,-2);\fi
      \ifnum\row=2 \fill[\col!18] (2,-4) rectangle (3,-2);\fi
      \fill[\col!35] (2,-3) rectangle (3,-2);
      \foreach \g in {0,2,4}{
        \draw[fu] (\g,0)--(\g,-4);
        \draw[fu] (0,-\g)--(4,-\g);
      }
      \foreach \g in {1,3}{
        \draw[fp] (\g,0)--(\g,-4);
        \draw[fp] (0,-\g)--(4,-\g);
      }
      \node[font=\scriptsize,inner sep=1pt] at (3.5,-3.5) {$K$};
      \ifnum\row=0
        \draw[\col,very thick] (1,-3) rectangle (3,-1);
        \fill[\col] (2,-2) circle (2.5pt);
        \node[font=\scriptsize,text=\col,anchor=south east,inner sep=1.5pt]
          at (2,-2) {$\mathbf x$};
      \fi
      \ifnum\row=1
        \draw[\col,very thick] (2,-3) rectangle (4,-2);
        \draw[\col!65!black,line width=1.6pt,densely dashed] (3,-3)--(3,-2);
        \node[font=\scriptsize,text=\col!65!black,anchor=west,inner sep=2pt]
          at (3,-2.5) {$e$};
      \fi
      \ifnum\row=2
        \draw[\col,very thick] (2,-4) rectangle (3,-2);
        \draw[\col!65!black,line width=1.6pt,densely dashed] (2,-3)--(3,-3);
        \node[font=\scriptsize,text=\col!65!black,anchor=north,inner sep=2pt]
          at (2.5,-3) {$e$};
      \fi
      \node[font=\scriptsize] at (2.5,-2.5) {$\tau$};
    \end{scope}
    \draw[flow,\col] (2.45,-1.05)--(3.2,-1.05);
    \node[draw=\col,rounded corners=3pt,fill=\col!5,
      minimum width=2.55cm,minimum height=.95cm,align=center] (map\row) at (4.65,-1.05)
      {$\mathbf x\mapsto\boldsymbol\xi_{P_{\comp}(\tau)}(\mathbf x)$\\[-1pt]
       {\scriptsize $\boldsymbol\xi\in[-1,1]^2$}};
    \foreach \iy in {-.65,-1.45}{
      \foreach \hy in {-.25,-.8,-1.35,-1.9}{
        \draw[wire] (6.65,\iy)--(8.35,\hy);
      }
    }
    \foreach \hy in {-.25,-.8,-1.35,-1.9}{
      \draw[wire] (8.35,\hy)--(10.05,-1.05);
    }
    \draw[flow,\col] (map\row.east)--(6.18,-1.05);
    \node[neuron,draw=\col] at (6.65,-.65) {\scriptsize $\xi_1$};
    \node[neuron,draw=\col] at (6.65,-1.45) {\scriptsize $\xi_2$};
    \foreach \hy in {-.25,-.8,-1.35,-1.9}{
      \node[neuron,draw=\col,fill=\col!8] at (8.35,\hy) {};
    }
    \node[neuron,draw=\col,fill=\col!8] at (10.05,-1.05) {\scriptsize $\Sigma$};
    \draw[flow,\col] (10.48,-1.05)--(11.25,-1.05);
  \end{scope}
}
\node[text=scalar,font=\small] at (1.1,-2.48) {$P_u(\tau)$};
\node[text=px,font=\small] at (1.1,-5.53) {$P_{p_x}(\tau)$};
\node[text=py,font=\small] at (1.1,-8.58) {$P_{p_y}(\tau)$};
\node[draw=scalar,rounded corners=3pt,fill=scalar!5,
  minimum width=5.3cm,minimum height=.95cm,align=center] at (14.05,-1.05)
  {$u_\theta|_\tau=N^u_{\theta,P_u(\tau)}\circ\boldsymbol\xi_{P_u(\tau)}$};
\node[font=\scriptsize,align=center,text=scalar] at (14.05,-1.85)
  {continuous across $\mathcal F_u^0$, may jump across $\mathcal F_p$};
\node[draw=px,rounded corners=3pt,fill=px!5,
  minimum width=5.3cm,minimum height=.95cm,align=center] at (14.05,-4.1)
  {$p_{\eta,x}|_\tau=N^{p_x}_{\eta,P_{p_x}(\tau)}\circ\boldsymbol\xi_{P_{p_x}(\tau)}$};
\node[font=\scriptsize,align=center,text=px] at (14.05,-4.9)
  {single-valued across $e$; $p_{\eta,y}$ may jump across $e$};
\node[draw=py,rounded corners=3pt,fill=py!5,
  minimum width=5.3cm,minimum height=.95cm,align=center] at (14.05,-7.15)
  {$p_{\eta,y}|_\tau=N^{p_y}_{\eta,P_{p_y}(\tau)}\circ\boldsymbol\xi_{P_{p_y}(\tau)}$};
\node[font=\scriptsize,align=center,text=py] at (14.05,-7.95)
  {single-valued across $e$; $p_{\eta,x}$ may jump across $e$};
\end{tikzpicture}
\SDGNNFinish}
\endgroup
\caption{Construction of the neural network trial functions on a
subrectangle $\tau$. Each row shows, for $u$, $p_x$, and $p_y$, the patch
containing $\tau$, the affine map onto $[-1,1]^2$, the network on this patch
with four of its $r$ hidden units drawn, and the restriction of the trial
function to $\tau$.
Thin solid and dashed lines are edges in $\mathcal F_u$ and $\mathcal F_p$, respectively, and thick colored lines outline the patches. The continuity stated on the right holds for all
parameter values.}
\label{fig:neural-network-architecture}
\end{figure}

The restrictions of the scalar and diffusive
trial functions to each $\tau\in\mathcal T$ are defined by
\[
\begin{aligned}
u_\theta|_\tau(\mathbf x)
&:=
N_{\theta,P_u(\tau)}^u
\bigl(\boldsymbol\xi_{P_u(\tau)}(\mathbf x)\bigr),
\\
\mathbf p_\eta|_\tau(\mathbf x)
&:=
\begin{pmatrix}
N_{\eta,P_{p_x}(\tau)}^{p_x}
\bigl(\boldsymbol\xi_{P_{p_x}(\tau)}(\mathbf x)\bigr)
\\[2pt]
N_{\eta,P_{p_y}(\tau)}^{p_y}
\bigl(\boldsymbol\xi_{P_{p_y}(\tau)}(\mathbf x)\bigr)
\end{pmatrix},
\qquad \mathbf x\in\tau.
\end{aligned}
\]
Figure~\ref{fig:neural-network-architecture} shows this construction
for one subrectangle $\tau$.

On each patch of its family, each of $u_\theta$ and the components
$p_{\eta,x}$ and $p_{\eta,y}$ of $\mathbf p_\eta=(p_{\eta,x},p_{\eta,y})$ is one
network composed with one affine map and is therefore smooth on the closure of that patch.
Hence $u_\theta\in\mathbb V$ and $\mathbf p_\eta\in\mathbb W$ for all values of the network
parameters, and the SDG spaces require no other continuity between
subrectangles. The trial functions permit jumps in $u$ across $\mathcal F_p$,
in both components of $\mathbf p$ across $\mathcal F_u^0$, and in its
tangential component across $\mathcal F_p$, since the two subrectangles adjacent to such an edge lie in different patches of the family used for that component; the jumps in $u$
can approximate a
prescribed jump across an interface aligned with $\mathcal F_p$ (Examples~3 and~4 in
Section~\ref{sec:numerical-experiments}). The scalar networks are not
constrained on $\partial\Omega$; the boundary residual in
Section~\ref{subsec:weak-residual-functionals} imposes the data $g$ through
moments. Thus, the neural network trial class satisfies
\[
\left\{(u_\theta,\mathbf p_\eta):
(\theta,\eta)\in{\Theta_u\times\Theta_{\mathbf p}}\right\}
\subset \mathbb V\times{\mathbb W},
\]
and $\mathcal A_\partial$ is well-defined on this class.

\section{Weak residual loss and error decomposition}
\label{sec:weak-residual-loss}

\subsection{Weak residual functionals}
\label{subsec:weak-residual-functionals}

For each $(\theta,\eta)\in{\Theta_u\times\Theta_{\mathbf p}}$, the weak
residuals of the constitutive and equilibrium equations and the boundary
residual are the linear functionals $R_1(\theta,\eta;\cdot)\in({W_h^k})'$,
$R_2(\theta,\eta;\cdot)\in({U_{h,0}^k})'$, and
$R_\partial(\theta;\cdot)\in(U_h^k|_{\partial\Omega})'$, respectively. Their
actions on $\mathbf z_h\in {W_h^k}$, $w_h\in {U_{h,0}^k}$, and
$\mu_h\in U_h^k|_{\partial\Omega}$ are defined by
\begin{equation}
\label{eq:scalar-residual-components}
\begin{aligned}
R_1(\theta,\eta;\mathbf z_h)
&:=
(\alpha\mathbf p_\eta,\mathbf z_h)_\Omega
-B^*(u_\theta,\mathbf z_h)
-\int_{\partial\Omega}g\,\mathbf z_h\cdot\mathbf n\,d\sigma,
\\
R_2(\theta,\eta;w_h)
&:=
B(\mathbf p_\eta,w_h)
+C_h^{\mathrm{up}}(u_\theta,w_h)
+(cu_\theta,w_h)_\Omega
-(f,w_h)_\Omega,
\\
R_\partial(\theta;\mu_h)
&:=
\int_{\partial\Omega}(u_\theta-g)\,\mu_h\,d\sigma.
\end{aligned}
\end{equation}

The boundary residual is the moment equation of
Section~\ref{sec:sdg-formulation} with $u_\theta$ in place of $u_{h,g}^k$.
The forms $B^*$, $B$, and $C_h^{\mathrm{up}}$, the mass and reaction terms,
and the boundary moments act on the trial pair only through its values and
traces, and all derivatives fall on the test polynomials, so evaluating the
residuals requires no derivatives of the networks.

With $\mathcal A_\partial$ and $F_{g,\partial}$ of
Section~\ref{sec:sdg-formulation}, combining the three components in
\eqref{eq:scalar-residual-components} gives the full residual functional
$R_h(\theta,\eta)\in(Y_h^k)'$:
\begin{equation}
\label{eq:scalar-residual-functional}
\begin{aligned}
R_h\bigl(\theta,\eta;(\Psi_h,\mu_h)\bigr)
&:=
R_1(\theta,\eta;\mathbf z_h)
+R_2(\theta,\eta;w_h)
+R_\partial(\theta;\mu_h)
\\
&=
\mathcal A_\partial\bigl((u_\theta,\mathbf p_\eta),(\Psi_h,\mu_h)\bigr)
-F_{g,\partial}(\Psi_h,\mu_h).
\end{aligned}
\end{equation}

\subsection{Matrix representation and SDGNN loss}
\label{subsec:matrix-representation-norm}

Fixed bases of ${U_{h,0}^k}$ and ${W_h^k}$ are denoted by
$\{\varphi_i\}_{i=1}^{n_U}$ and
$\{\boldsymbol\zeta_j\}_{j=1}^{n_W}$, respectively, so that
$n_U=\dim U_{h,0}^k$ and $n_W=\dim W_h^k$. With
$n_\partial:=\dim U_h^k|_{\partial\Omega}$, we also fix
$\chi_m\in U_h^k$, $1\le m\le n_\partial$, whose traces
$\{\chi_m|_{\partial\Omega}\}_{m=1}^{n_\partial}$ form a basis of
$U_h^k|_{\partial\Omega}$, so that
$U_h^k=U_{h,0}^k\oplus\operatorname{span}\{\chi_m\}_{m=1}^{n_\partial}$.
Both $X_h^k$ and $Y_h^k$ then have dimension
$n_h:=n_U+n_W+n_\partial$. We order
the trial basis $\{\phi_\ell\}_{\ell=1}^{n_h}$ as $(\varphi_i,\mathbf0)$,
$(0,\boldsymbol\zeta_j)$, $(\chi_m,\mathbf0)$ and the test basis
$\{\psi_\ell\}_{\ell=1}^{n_h}$ as $(\boldsymbol\zeta_j,0,0)$,
$(\mathbf0,\varphi_i,0)$, $(\mathbf0,0,\chi_m|_{\partial\Omega})$, for
$1\le i\le n_U$, $1\le j\le n_W$, and $1\le m\le n_\partial$.

We write $\langle\cdot,\cdot\rangle_{X_h^k}$ for the inner
product on $X_h^k$ whose norm is the restriction of $\|\cdot\|_{*,h}$
in \eqref{eq:projection-input-norm}. The Gram matrix
$M_X\in\mathbb R^{n_h\times n_h}$ of the trial basis
$\{\phi_\ell\}_{\ell=1}^{n_h}$ with respect to this inner product has entries
\begin{equation}
\label{eq:trial-gram-matrix}
(M_X)_{\ell m}
:=
\langle\phi_m,\phi_\ell\rangle_{X_h^k},
\qquad
1\leq\ell,m\leq n_h.
\end{equation}
The matrix $M_X$ is symmetric positive definite.
We use it to measure errors in
Theorem~\ref{prop:scalar-loss-identity}, and it does not enter the loss.

The matrix ${A}_h\in\mathbb R^{n_h\times n_h}$ represents
$\mathcal A_\partial|_{X_h^k\times Y_h^k}$ in the trial
basis $\{\phi_j\}_{j=1}^{n_h}$ and the test basis $\{\psi_i\}_{i=1}^{n_h}$,
with entries
\begin{equation}
\label{eq:sdg-matrix-entries}
({A}_h)_{ij}
:=
\mathcal A_\partial(\phi_j,\psi_i),
\qquad
1\leq i,j\leq n_h.
\end{equation}
The residual vector is
\begin{equation}
\label{eq:full-residual-vector}
\mathbf r_h(\theta,\eta)
:=
\bigl(R_h(\theta,\eta;\psi_\ell)\bigr)_{\ell=1}^{n_h}\in\mathbb R^{n_h};
\end{equation}
by \eqref{eq:scalar-residual-functional}, its three blocks contain the values
of $R_1$, $R_2$, and $R_\partial$ on the bases of $W_h^k$, $U_{h,0}^k$, and
$U_h^k|_{\partial\Omega}$.

Since the SDG finite element approximation
\eqref{eq:sdg-boundary-problem} exists and is unique for all data
(Section~\ref{sec:sdg-formulation}), the matrix ${A}_h$ is nonsingular.

For $e\in\mathcal F_\partial$, let $\mathbf n_e:=\mathbf n$. On
$Y_h^k$, we use the inner product
\[
\begin{aligned}
\bigl\langle(\mathbf z_h,w_h,\mu_h),(\mathbf z_h',w_h',\mu_h')\bigr\rangle_{Y_h^k}
&:=(\alpha^2\mathbf z_h,\mathbf z_h')_\Omega
+\sum_{e\in\mathcal F_p\cup\mathcal F_\partial}h_e
(\alpha^2\,\mathbf z_h\cdot\mathbf n_e,\mathbf z_h'\cdot\mathbf n_e)_{e}\\
&\quad+\sum_{\tau\in\mathcal T}(\nabla_hw_h,\nabla_hw_h')_{\tau}
+\sum_{e\in\mathcal F_p}h_e^{-1}([w_h]_e,[w_h']_e)_{e}\\
&\quad+\sum_{e\in\mathcal F_\partial}h_e(\mu_h,\mu_h')_{e}.
\end{aligned}
\]
In the edge terms, $\alpha$ is assumed to have a single-valued trace on each
edge with values between $\alpha_{\min}:=\operatorname*{ess\,inf}_{\Omega}\alpha$ and $\alpha_{\max}:=\|\alpha\|_{L^\infty(\Omega)}$, which holds for the constant coefficients in Examples~1--3 of Section~\ref{sec:numerical-experiments}.
The weight $\alpha^2$ measures the constitutive residual on the scale of
the diffusive variable: for $u\in H^1(\Omega)$ and $\mathbf p\in L^2(\Omega)^2$, over all $\mathbf z\in L^2(\Omega)^2$, the dual norm of
$\mathbf z\mapsto(\alpha\mathbf p-\nabla u,\mathbf z)_\Omega$ with respect to
$\|\alpha\mathbf z\|_{L^2(\Omega)}$ is $\|\mathbf p-a\nabla u\|_{L^2(\Omega)}$.
The Gram matrix
$M_Y\in\mathbb R^{n_h\times n_h}$ of the test basis
$\{\psi_\ell\}_{\ell=1}^{n_h}$ has entries
\begin{equation}
\label{eq:residual-metric-matrix}
(M_Y)_{\ell m}
:=
\langle\psi_m,\psi_\ell\rangle_{Y_h^k},
\qquad
1\leq\ell,m\leq n_h.
\end{equation}
The matrix $M_Y$ is symmetric positive definite and block diagonal with
respect to the vector, scalar, and moment components. We write
$\|\cdot\|_{Y_h^k}$ for the norm induced by this inner product.

The SDGNN loss is the squared $(Y_h^k)'$-norm of
$R_h(\theta,\eta)$. Since $M_Y$ is the Gram matrix of the test basis and
$\mathbf r_h(\theta,\eta)$ contains the values of $R_h(\theta,\eta)$ on this
basis, the Riesz representation in $Y_h^k$ gives the second
expression in
\begin{equation}
\label{eq:sdgnn-loss}
\begin{aligned}
{J_h(k;\theta,\eta)}
&:=
\left\|
R_h(\theta,\eta)
\right\|_{(Y_h^k)'}^2
\\
&=
\mathbf r_h(\theta,\eta)^{\mathsf T}
M_Y^{-1}
\mathbf r_h(\theta,\eta).
\end{aligned}
\end{equation}

\subsection{SDG projection and error decomposition}
\label{subsec:sdgnn-loss-projection}

For $\Phi\in\mathbb V\times{\mathbb W}$, the SDG projection
${\Pi_h^k}\Phi\in X_h^k$ is defined by
\begin{equation}
\label{eq:sdg-projection}
\mathcal A_\partial\bigl({\Pi_h^k}\Phi,(\Psi_h,\mu_h)\bigr)
=\mathcal A_\partial\bigl(\Phi,(\Psi_h,\mu_h)\bigr)
\qquad\text{for every }(\Psi_h,\mu_h)\in Y_h^k.
\end{equation}
These equations also have the coefficient matrix ${A}_h$, so the
projection exists and is unique, and it satisfies ${\Pi_h^k}\Phi_h=\Phi_h$ for
every $\Phi_h\in X_h^k$. The uniqueness also shows that $\Pi_h^k$ is linear.

The following theorem states that the loss measures, up to the constants
$\sigma_{\min}$ and $\sigma_{\max}$, the distance between the SDG projection of
the network and the SDG finite element approximation. We call
$\|\Pi_h^k(u_\theta,\mathbf p_\eta)-({u_{h,g}^k},{\mathbf p_{h,g}^k})\|_{*,h}$ the training error.

\begin{theorem}[Loss and training error]
\label{prop:scalar-loss-identity}
For every $(\theta,\eta)\in{\Theta_u\times\Theta_{\mathbf p}}$, the SDGNN loss
satisfies
\begin{equation}
\label{eq:scalar-loss-representation}
\sigma_{\min}^2
\left\|
\Pi_h^k(u_\theta,\mathbf p_\eta)
-({u_{h,g}^k},{\mathbf p_{h,g}^k})
\right\|_{*,h}^2
\le
{J_h(k;\theta,\eta)}
\le
\sigma_{\max}^2
\left\|
\Pi_h^k(u_\theta,\mathbf p_\eta)
-({u_{h,g}^k},{\mathbf p_{h,g}^k})
\right\|_{*,h}^2,
\end{equation}
where $\sigma_{\min}$ and $\sigma_{\max}$ are the smallest and largest singular
values of $M_Y^{-1/2}{A}_hM_X^{-1/2}$; $\sigma_{\min}>0$ since ${A}_h$ is nonsingular. Equivalently,
\[
\sigma_{\min}
=
\min_{0\ne\Phi_h\in X_h^k}\,
\max_{0\ne(\Psi_h,\mu_h)\in Y_h^k}
\frac{\mathcal A_\partial\bigl(\Phi_h,(\Psi_h,\mu_h)\bigr)}
{\|\Phi_h\|_{*,h}\,\|(\Psi_h,\mu_h)\|_{Y_h^k}}
\]
is the discrete inf-sup constant of $\mathcal A_\partial$ in these norms, and
$\sigma_{\max}$ is the norm of $\mathcal A_\partial$ on
$X_h^k\times Y_h^k$.
\end{theorem}

\begin{proof}
The vector $\mathbf e_h(\theta,\eta)\in\mathbb R^{n_h}$ contains the
coefficients of
\[
\Pi_h^k(u_\theta,\mathbf p_\eta)
-({u_{h,g}^k},{\mathbf p_{h,g}^k})
\]
in the basis $\{\phi_\ell\}_{\ell=1}^{n_h}$ of $X_h^k$.
For every $(\Psi_h,\mu_h)\in Y_h^k$, subtracting
\eqref{eq:sdg-boundary-problem} from \eqref{eq:sdg-projection} with
$\Phi=(u_\theta,\mathbf p_\eta)$ and using
\eqref{eq:scalar-residual-functional} yields
\begin{equation}
\label{eq:projected-residual-functional}
\begin{aligned}
R_h\bigl(\theta,\eta;(\Psi_h,\mu_h)\bigr)
&=\mathcal A_\partial\bigl((u_\theta,\mathbf p_\eta),(\Psi_h,\mu_h)\bigr)
-F_{g,\partial}(\Psi_h,\mu_h)\\
&=\mathcal A_\partial\bigl(
\Pi_h^k(u_\theta,\mathbf p_\eta)
-({u_{h,g}^k},{\mathbf p_{h,g}^k}),(\Psi_h,\mu_h)
\bigr).
\end{aligned}
\end{equation}
Evaluating \eqref{eq:projected-residual-functional} at
$(\Psi_h,\mu_h)=\psi_\ell$ for $1\leq\ell\leq n_h$ and using
\eqref{eq:full-residual-vector} and \eqref{eq:sdg-matrix-entries} gives
\begin{equation}
\label{eq:residual-error-vector}
\mathbf r_h(\theta,\eta)
=
{A}_h\mathbf e_h(\theta,\eta).
\end{equation}
Let $T_h:=M_Y^{-1/2}{A}_hM_X^{-1/2}$ and
$\mathbf y_h:=M_X^{1/2}\mathbf e_h(\theta,\eta)$. Substituting
\eqref{eq:residual-error-vector} into \eqref{eq:sdgnn-loss} yields
\[
J_h(k;\theta,\eta)
=
\mathbf e_h^{\mathsf T}{A}_h^{\mathsf T}M_Y^{-1}{A}_h\mathbf e_h
=
|T_h\mathbf y_h|^2,
\]
and the Gram matrix definition \eqref{eq:trial-gram-matrix} gives
$|\mathbf y_h|^2=\mathbf e_h^{\mathsf T}M_X\mathbf e_h
=\|{\Pi_h^k}(u_\theta,\mathbf p_\eta)-({u_{h,g}^k},{\mathbf p_{h,g}^k})\|_{*,h}^2$.
The inequalities $\sigma_{\min}|\mathbf y_h|\le|T_h\mathbf y_h|\le\sigma_{\max}|\mathbf y_h|$
prove \eqref{eq:scalar-loss-representation}. For $\Phi_h$ with coefficient
vector $\mathbf c$ and $(\Psi_h,\mu_h)$ with coefficient vector $\mathbf d$,
\eqref{eq:sdg-matrix-entries}, \eqref{eq:trial-gram-matrix}, and
\eqref{eq:residual-metric-matrix} give
$\mathcal A_\partial(\Phi_h,(\Psi_h,\mu_h))=\mathbf d^{\mathsf T}{A}_h\mathbf c$,
$\|\Phi_h\|_{*,h}=|M_X^{1/2}\mathbf c|$, and
$\|(\Psi_h,\mu_h)\|_{Y_h^k}=|M_Y^{1/2}\mathbf d|$, which gives the
characterization of $\sigma_{\min}$ and $\sigma_{\max}$
(cf.~\cite[Remark~2.23(iii)]{2004Ern}).
\end{proof}

\par
The lower bound in \eqref{eq:scalar-loss-representation} is the a priori estimate of the discrete
Banach--Ne\v{c}as--Babu\v{s}ka theorem~\cite[Theorems~2.6 and~2.22]{2004Ern} (see
also~\cite[Theorem~2.1]{1971Babuska}) applied to \eqref{eq:projected-residual-functional}, and the
upper bound is the boundedness of $\mathcal A_\partial$ (cf.\ the two-sided residual bound for
discontinuous Petrov--Galerkin methods~\cite[Theorem~2.1]{2014Carstensen}).
\par

For the exact solution $(u,\mathbf p)$ with boundary data $g$, we
define the discretization error $D_h(k)$ and the network projection error
$P_h(k;\theta,\eta)$ in the norm \eqref{eq:full-error-norm} by
\begin{align}
\label{eq:sdg-discretization-error-degree}
{D_h(k)}
&:=\left\|
(u,\mathbf p)
-({u_{h,g}^k},{\mathbf p_{h,g}^k})
\right\|_h,\\
\label{eq:neural-sdg-projection-difference}
{P_h(k;\theta,\eta)}
&:=\left\|
(u_\theta,\mathbf p_\eta)
-{\Pi_h^k}(u_\theta,\mathbf p_\eta)
\right\|_h.
\end{align}

The following theorem decomposes the error of the network into three parts:
the discretization error, which does not depend on the network; the training
error, which the loss measures; and the network projection error, which the loss does not
see.

\begin{theorem}[Error decomposition]
\label{prop:scalar-error-estimate}
Let $(u,\mathbf p)$ be the exact solution with boundary data $g$. For every
$(\theta,\eta)\in{\Theta_u\times\Theta_{\mathbf p}}$,
\[
\left\|(u,\mathbf p)-(u_\theta,\mathbf p_\eta)\right\|_h
\leq
{D_h(k)}
+\sigma_{\min}^{-1}{J_h(k;\theta,\eta)}^{1/2}
+{P_h(k;\theta,\eta)},
\]
with $\sigma_{\min}$ as in Theorem~\ref{prop:scalar-loss-identity}.
\end{theorem}

\begin{proof}
The full error decomposes as
\begin{equation}
\label{eq:scalar-error-decomposition}
\begin{aligned}
(u,\mathbf p)-(u_\theta,\mathbf p_\eta)
&=\bigl[(u,\mathbf p)-({u_{h,g}^k},{\mathbf p_{h,g}^k})\bigr]\\
&\quad+\bigl[({u_{h,g}^k},{\mathbf p_{h,g}^k})
-{\Pi_h^k}(u_\theta,\mathbf p_\eta)\bigr]\\
&\quad+\bigl[{\Pi_h^k}(u_\theta,\mathbf p_\eta)
-(u_\theta,\mathbf p_\eta)\bigr].
\end{aligned}
\end{equation}
The second difference in \eqref{eq:scalar-error-decomposition} lies in
$X_h^k$, and $\|\cdot\|_h\le\|\cdot\|_{*,h}$. The lower bound in
\eqref{eq:scalar-loss-representation} therefore gives
\begin{equation}
\label{eq:projected-error-bound}
\begin{aligned}
\left\|({u_{h,g}^k},{\mathbf p_{h,g}^k})
-{\Pi_h^k}(u_\theta,\mathbf p_\eta)\right\|_h
&\leq\left\|({u_{h,g}^k},{\mathbf p_{h,g}^k})
-{\Pi_h^k}(u_\theta,\mathbf p_\eta)\right\|_{*,h}\\
&\le\sigma_{\min}^{-1}{J_h(k;\theta,\eta)}^{1/2}.
\end{aligned}
\end{equation}
Applying the triangle inequality to
\eqref{eq:scalar-error-decomposition} and using
\eqref{eq:projected-error-bound} proves the estimate.
\end{proof}

When ${\Pi_h^k}(u,\mathbf p)=({u_{h,g}^k},{\mathbf p_{h,g}^k})$, which
holds under the assumption of Corollary~\ref{cor:sdg-fe-quasi-optimality} on the meshes of Section~\ref{sec:projection-error-estimate}, the
second difference in \eqref{eq:scalar-error-decomposition} is
${\Pi_h^k}\bigl((u,\mathbf p)-(u_\theta,\mathbf p_\eta)\bigr)$ by the linearity of $\Pi_h^k$, so the loss
measures, up to the constants in \eqref{eq:scalar-loss-representation}, the
part of the error that the SDG projection retains.
For a given network, the three terms in
Theorem~\ref{prop:scalar-error-estimate} can be computed separately. The loss
$J_h$ is the training objective, $P_h$ requires one solve with ${A}_h$ for
the SDG projection and a quadrature of the difference, and $D_h$ requires the
SDG finite element approximation and the exact solution.
The constant $\sigma_{\min}$ depends on the mesh, the degree, and the
coefficients, but not on the network. Its square is the smallest
eigenvalue $\lambda$ of the generalized eigenvalue problem
${A}_h^{\mathsf T}M_Y^{-1}{A}_h\mathbf c=\lambda M_X\mathbf c$~\cite[Theorem~2]{1981Malkus},
\cite[Proposition~2.1]{1993Chapelle}.

\begin{remark}
\label{rem:trial-inner-product}
If $M_Y$ in \eqref{eq:sdgnn-loss} is replaced by
$G:={A}_hM_X^{-1}{A}_h^{\mathsf T}$, the argument in the proof of
Theorem~\ref{prop:scalar-loss-identity} gives for the resulting loss $J_h^G$ the identity
$J_h^G(k;\theta,\eta)=\|{\Pi_h^k}(u_\theta,\mathbf p_\eta)-({u_{h,g}^k},{\mathbf p_{h,g}^k})\|_{*,h}^2$.
The matrix $G$ is the Gram matrix of the norm
$\|(\Psi_h,\mu_h)\|_G:=\sup_{0\ne\Phi_h\in X_h^k}
\mathcal A_\partial(\Phi_h,(\Psi_h,\mu_h))/\|\Phi_h\|_{*,h}$ on
$Y_h^k$, which is the discrete counterpart of the optimal test norm of
discontinuous Petrov--Galerkin methods~\cite[Section~2.2 and
Proposition~2.1]{2011Zitelli} (see also~\cite[eqs.~(2.12)--(2.14)]{2012Dahmen}). The test norm $\|\cdot\|_{Y_h^k}$ induced
by $M_Y$ satisfies
$\sigma_{\min}\|\cdot\|_{Y_h^k}\le\|\cdot\|_G\le\sigma_{\max}\|\cdot\|_{Y_h^k}$.
It is thus an equivalent test norm in the sense of
\cite[Section~2.4]{2011Zitelli}, and Theorem~\ref{prop:scalar-loss-identity} is
the discrete form of the resulting norm equivalence
\cite[eq.~(2.18)]{2011Zitelli}.
With $G$, each evaluation of the loss requires a solve with ${A}_h$
and each gradient a further solve with ${A}_h^{\mathsf T}$. We use the test inner product in
\eqref{eq:residual-metric-matrix} instead, so that evaluating the loss does
not require the SDG finite element system.
\end{remark}

\begin{remark}
\label{rem:loss-kernel-example}
For a fixed degree $k$, the loss can vanish for a network whose error is large.
Let $a=1$, $\mathbf b=\mathbf 0$, $c=f=g=0$, so that the exact solution and the
SDG finite element approximation vanish. On a patch $P\in\mathcal P_{p_x}$, let $t\in[-1,1]$ be the horizontal reference coordinate, so that the
vertical edge of $\mathcal F_p$ in $P$ lies at $t=0$. For distinct
$W_1,\dots,W_{k+2}>0$, the $k+1$ linear conditions
$\sum_{\ell}a_\ell\int_0^1t^j\tanh(W_\ell t)\,dt=0$, $0\le j\le k$, on the
$k+2$ coefficients $a_\ell$ have a solution $\mathbf a\ne\mathbf 0$. Let
$\phi(t)=\sum_\ell a_\ell\tanh(W_\ell t)$. For the largest $W_\ell$ with
$a_\ell\ne0$, the pole of $\tanh(W_\ell t)$ at $i\pi/(2W_\ell)$ is closer to the
real axis than every pole of the other terms, so $\phi$ has a pole there and is
not identically zero. By construction, $\phi$ is orthogonal to the polynomials of degree at most
$k$ on $(0,1)$; since $\phi$ is odd, it is also orthogonal to them on
$(-1,0)$, and $\phi(0)=0$. The trial pair with $u_\theta=0$,
$\mathbf p_\eta=(\phi,0)$ on $P$, and $\mathbf p_\eta=\mathbf 0$ elsewhere is
given by networks of local hidden width $r\ge k+2$ (for instance, $r=6$ when
$k=4$, as in Section~\ref{sec:numerical-experiments}), and all residuals vanish by the orthogonality of $\phi$ and $\phi(0)=0$. Hence $J_h=0$ and
$\Pi_h^k(u_\theta,\mathbf p_\eta)=0$, whereas the full error equals
$P_h=\|\mathbf p_\eta\|_{L^2(\Omega)}$, which grows without bound as $\mathbf a$
is scaled. On the uniform meshes of Section~\ref{sec:projection-error-estimate},
Theorem~\ref{thm:fixed-network-degree-convergence} shows that, for a fixed network,
$P_h$ tends to zero as $k$ increases.
\end{remark}

\par
The integrals in \eqref{eq:scalar-residual-components} that contain the
networks are evaluated by Gauss--Legendre rules whose orders are predicted,
subrectangle by subrectangle, from the current network weights through the
analyticity of the tanh units, without evaluating residuals.
Appendix~\ref{subsec:weight-quadrature} describes this quadrature;
Section~\ref{subsec:numerical-settings} compares the computed losses with a
higher-order rule.
\par

\section{Estimates of the discretization and \texorpdfstring{network}{network} projection errors}
\label{sec:projection-error-estimate}

\subsection{Stability of the SDG projection}
\label{subsec:sdg-projection-approximation}

By \eqref{eq:sdg-projection}, the neural network trial pair
and its SDG projection have the same residual functional on
$Y_h^k$, so
the loss carries no information about their difference.
To bound $\sigma_{\min}$, $D_h(k)$, and $P_h(k;\theta,\eta)$ in
Theorem~\ref{prop:scalar-error-estimate}, we assume
in this section that $\Omega=(0,1)^2$, use a
uniform primal partition into squares of side length $h$, and take
$k\ge1$.

For test pairs $\Psi_h=(\mathbf z_h,w_h)\in {Y_{h,0}^k}$, we write
$\|\Psi_h\|_{*,h}:=\|(w_h,\mathbf z_h)\|_{*,h}$, and for test
triples $(\Psi_h,\mu_h)\in Y_h^k$,
\[
\|(\Psi_h,\mu_h)\|_{*,h}^2
:=\|\Psi_h\|_{*,h}^2
+\sum_{e\in\mathcal F_\partial}h_e\|\mu_h\|_{L^2(e)}^2.
\]
Unlike this norm, the norm $\|\cdot\|_{Y_h^k}$
induced by $M_Y$ in \eqref{eq:residual-metric-matrix}, whose dual defines the
SDGNN loss \eqref{eq:sdgnn-loss}{,} weights the vector component by
$\alpha$ and includes the normal traces on $\partial\Omega$.
For a function $\gamma$ on $\partial\Omega$, we write
\[
|\gamma|_{\partial,h}^2
:=\sum_{e\in\mathcal F_\partial}
h_e^{-1}\|\gamma\|_{L^2(e)}^2,
\]
and $Q_\partial$ denotes the $L^2(\partial\Omega)$-orthogonal projection onto
$U_h^k|_{\partial\Omega}$.

\begin{lemma}[Trace bounds for polynomials]
\label{lem:polynomial-trace}
For $d>0$ and an integer $n\ge0$, let $\mathbb P_{n}(0,d)$ denote the polynomials of degree at most $n$
on $(0,d)$. For $q\in\mathbb P_{n}(0,d)$,
\[
\|q\|_{L^\infty(0,d)}^2
\le\frac{(n+1)^2}{d}\|q\|_{L^2(0,d)}^2;
\]
in particular, this bound holds for $|q(0)|^2$ and $|q(d)|^2$.
For every $\mathbf z_h\in W_h^k$,
\begin{equation}
\label{eq:boundary-normal-trace}
\sum_{e\in\mathcal F_\partial}
h_e\|\mathbf z_h\cdot\mathbf n\|_{L^2(e)}^2
\le(k+1)^2\|\mathbf z_h\|_{L^2(\Omega)}^2
\end{equation}
and
\begin{equation}
\label{eq:degree-inverse-trace}
\|\mathbf z_h\|_{0,h}^2\le\Bigl(1+\frac{(k+1)^2}{2}\Bigr)\|\mathbf z_h\|_{L^2(\Omega)}^2.
\end{equation}
\end{lemma}
\begin{proof}
The Legendre polynomial $L_i$ of degree $i$ is normalized by
$L_i(1)=1$, so that $|L_i|\le1$ on $[-1,1]$~\cite{OLBC10}. The Cauchy--Schwarz inequality
applied to the expansion of $q$ in the orthonormal basis
$\{\sqrt{(2i+1)/d}\,L_i(2t/d-1)\}_{i=0}^{n}$ of $\mathbb P_n(0,d)$ and
$\sum_{i=0}^{n}(2i+1)=(n+1)^2$ give the first bound. On a subrectangle $\tau$, a square of side $d=h_e$, the normal
component of $\mathbf z_h$ is, for each value of the tangential coordinate,
a polynomial of degree at most $k$ in the normal coordinate, so integrating
the first bound along an edge $e$ of $\tau$ gives
$h_e\|\mathbf z_h\cdot\mathbf n_e\|_{L^2(e)}^2
\le(k+1)^2\|\mathbf z_h\cdot\mathbf n_e\|_{L^2(\tau)}^2$.
Each boundary edge belongs to one subrectangle, and the boundary edges
of a subrectangle have different normal directions, which gives
\eqref{eq:boundary-normal-trace}. Each subrectangle has one vertical and
one horizontal edge in $\mathcal F_p$, and summing over the subrectangles counts
every single-valued normal trace twice; the trace contribution to
$\|\mathbf z_h\|_{0,h}^2$ is therefore at most
$(k+1)^2\|\mathbf z_h\|_{L^2(\Omega)}^2/2$
(cf.~\cite[Lemma~1.46 and Remark~1.48]{DE12}).
\end{proof}

In these norms, the boundary integral in $\mathcal A$ is bounded through
\eqref{eq:boundary-normal-trace}, which brings in the factor $k+1$, so the
next proposition leaves it out.

\begin{proposition}[Continuity of the weak form]
\label{prop:weak-form-continuity}
For $\Phi=(v,\mathbf p)\in\mathbb V\times{\mathbb W}$
and $\Psi_h=(\mathbf z_h,w_h)\in Y_{h,0}^k$, let
\[
\mathcal A_0(\Phi,\Psi_h)
:=\mathcal A(\Phi,\Psi_h)
-\int_{\partial\Omega}v\,\mathbf z_h\cdot\mathbf n\,d\sigma.
\]
There exists a constant $C_A>0$, independent of $h$ and $k$, such that
\begin{equation}
\label{eq:weak-form-continuity}
|\mathcal A_0(\Phi,\Psi_h)|
\le C_A\|\Phi\|_{*,h}\|\Psi_h\|_{*,h}
\qquad
\text{for all }{\Phi\in\mathbb V\times{\mathbb W}},
\quad \Psi_h\in {Y_{h,0}^k}.
\end{equation}
The constant $C_A$ depends on the fixed domain and $\|\alpha\|_{L^\infty(\Omega)}$, $\|\mathbf b\|_{L^\infty(\Omega)}$, and $\|c\|_{L^\infty(\Omega)}$.
\end{proposition}
\begin{proof}
For $\Phi=(v,\mathbf p)$ and $\Psi_h=(\mathbf z_h,w_h)$,
the identity
$B(\mathbf z,v)=B^*(v,\mathbf z)+\int_{\partial\Omega}v\,\mathbf z\cdot\mathbf n\,d\sigma$
of Section~\ref{sec:staggered-framework} gives
\[
\mathcal A_0(\Phi,\Psi_h)
=(\alpha\mathbf p,\mathbf z_h)_\Omega
-B(\mathbf z_h,v)
+B(\mathbf p,w_h)
+C_h^{\mathrm{up}}(v,w_h)
+(cv,w_h)_\Omega.
\]
The Cauchy--Schwarz inequality gives
the standard SDG bounds (cf.~\cite[eq.~(14)]{KJN21})
\[
|B(\mathbf p,w_h)|
\le\|\mathbf p\|_{0,h}\|w_h\|_{1,h},
\qquad
|B(\mathbf z_h,v)|
\le\|\mathbf z_h\|_{0,h}\|v\|_{1,h};
\]
the second bound requires no boundary condition on $v$.
On the uniform square meshes, the broken Poincar\'e--Friedrichs inequality
\cite[eq.~(1.8) and Corollary~6.3]{B03} gives
\begin{equation}
\label{eq:projection-poincare}
\|v\|_{L^2(\Omega)}
\le C_P\bigl(\|v\|_{1,h}^2+\|v\|_{L^2(\partial\Omega)}^2\bigr)^{1/2}
\le C_P\|(v,\mathbf0)\|_{*,h},
\qquad v\in\mathbb V,
\end{equation}
where $C_P$ is independent of $h$ and $k$, and the last step uses
$h_e\le1$.
For the transport term, the continuous trace inequality
\cite[Lemma~1.49]{DE12},
$|v^{\mathrm{up}}|^2\le |v^+|^2+|v^-|^2$, and
\eqref{eq:projection-poincare} yield
\[
|C_h^{\mathrm{up}}(v,w_h)|
\le C\|\mathbf b\|_{L^\infty(\Omega)}
\|(v,\mathbf0)\|_{*,h}\|w_h\|_{1,h},
\]
where $C$ is independent of $h$ and $k$; see also the upwind estimates in
\cite[proof of Theorem~3.1]{ZP19}.
The mass term is bounded by
$\|\alpha\|_{L^\infty(\Omega)}\|\mathbf p\|_{0,h}\|\mathbf z_h\|_{0,h}$,
and \eqref{eq:projection-poincare}{, applied to $v$ and to $w_h$, which
vanishes on $\partial\Omega$,} bounds the reaction term by
$C_P^2\|c\|_{L^\infty(\Omega)}{\|(v,\mathbf0)\|_{*,h}}\|w_h\|_{1,h}$.
Summing these estimates proves \eqref{eq:weak-form-continuity}.
\end{proof}

\begin{lemma}[Discrete inf-sup estimate for $B$; {cf.~\cite[Theorem~1]{KJN21}}]
\label{lem:sdg-b-infsup}
For every $k\ge1$ and $v_h\in {U_{h,0}^k}$,
\begin{equation}
\label{eq:sdg-b-infsup}
{\beta_B(k)}\|v_h\|_{1,h}
\le
\sup_{\mathbf0\ne\mathbf z_h\in {W_h^k}}
\frac{B(\mathbf z_h,v_h)}{\|\mathbf z_h\|_{0,h}},
\qquad
{\beta_B(k)}=\frac1{8k^2}.
\end{equation}
More precisely, for every $v_h\in {U_{h,0}^k}$, there exists
$\mathbf z_{v_h}\in {W_h^k}$ such that
\begin{equation}
\label{eq:sdg-b-test-function}
\|\mathbf z_{v_h}\|_{0,h}\le\|v_h\|_{1,h},
\qquad
B(\mathbf z_{v_h},v_h)\ge{\beta_B(k)}\|v_h\|_{1,h}^2.
\end{equation}
\end{lemma}
The proof is given in Appendix~\ref{app:degree-b-proof}.

We combine Lemma~\ref{lem:sdg-b-infsup} with the upwind energy argument
of \cite[Lemma~3.1]{ZP19}, adapted to the present weak form.
We track the dependence on $k$ to estimate the network projection error.
In the rest of this section,
$\mathcal E_h:U_h^k|_{\partial\Omega}\to U_h^k$ denotes the discrete trace extension
of Lemma~\ref{lem:boundary-trace-extension}, which satisfies
$(\mathcal E_h\gamma)|_{\partial\Omega}=\gamma$ and
$\|\mathcal E_h\gamma\|_{1,h}\le C_{\mathcal E}(k)|\gamma|_{\partial,h}$ with
$C_{\mathcal E}(k)=\frac54(k+1)^2$.

\begin{proposition}[Discrete inf-sup stability]
\label{prop:weak-form-infsup}
(a) For each fixed $k\ge1$, there exists a constant ${\beta_A(k)}>0$,
independent of $h$, such that
\begin{equation}
\label{eq:weak-form-infsup}
{\beta_A(k)}\|\Phi_h\|_{*,h}
\le\sup_{0\ne\Psi_h\in {Y_{h,0}^k}}
\frac{|\mathcal A(\Phi_h,\Psi_h)|}{\|\Psi_h\|_{*,h}}
\qquad\text{for every }\Phi_h\in {X_{h,0}^k}.
\end{equation}
With $\alpha_{\min}=\operatorname*{ess\,inf}_{\Omega}\alpha>0$ and
$\alpha_{\max}=\|\alpha\|_{L^\infty(\Omega)}$, the constant
${\beta_A(k)}$ can be chosen such that
\begin{equation}
\label{eq:weak-form-infsup-degree}
{\beta_A(k)}^{-1}\le C_\beta k^4,
\qquad
C_\beta=\left(1+\frac{\alpha_{\min}}{8\alpha_{\max}^2}\right)
\max\left\{\frac6{\alpha_{\min}},
\frac{128\alpha_{\max}^2}{\alpha_{\min}}\right\}.
\end{equation}
(b) For each fixed $k\ge1$, there exists a constant
$\beta_\partial(k)>0$, independent of $h$, such that
\begin{equation}
\label{eq:extended-weak-form-infsup}
\beta_\partial(k)\|\Phi_h\|_{*,h}
\le\sup_{0\ne(\Psi_h,\mu_h)\in Y_h^k}
\frac{\bigl|\mathcal A_\partial\bigl(\Phi_h,(\Psi_h,\mu_h)\bigr)\bigr|}
{\|(\Psi_h,\mu_h)\|_{*,h}}
\qquad\text{for every }\Phi_h\in X_h^k,
\end{equation}
and $\beta_\partial(k)$ can be chosen such that
\begin{equation}
\label{eq:extended-weak-form-infsup-degree}
\beta_\partial(k)^{-1}
\le\beta_A(k)^{-1}\bigl(k+2+C_A(1+C_{\mathcal E}(k))\bigr)+1+C_{\mathcal E}(k).
\end{equation}
\end{proposition}
\begin{proof}
For $\Phi_h=(v_h,\mathbf p_h)$, we use the test function
\[
\Psi_h=(\mathbf p_h-{\delta(k)}\mathbf z_{v_h},v_h),
\qquad
{\delta(k)}=\frac{\alpha_{\min}{\beta_B(k)}}{\alpha_{\max}^2},
\]
where $\mathbf z_{v_h}$ satisfies \eqref{eq:sdg-b-test-function}.
By the definition of $\mathcal A$ and \eqref{eq:discrete-adjoint-identity},
$\mathcal A(\Phi_h,(\mathbf z_{v_h},0))=(\alpha\mathbf p_h,\mathbf z_{v_h})_\Omega-B(\mathbf z_{v_h},v_h)$.
This identity, the energy identity in the proof of
Proposition~\ref{prop:scalar-discrete-uniqueness},
\eqref{eq:coefficient-assumptions}, \eqref{eq:sdg-b-test-function}, and Young's inequality yield
\[
\begin{aligned}
\mathcal A(\Phi_h,\Psi_h)
&\ge\alpha_{\min}\|\mathbf p_h\|_{L^2(\Omega)}^2
+{\delta(k)}{\beta_B(k)}\|v_h\|_{1,h}^2
-{\delta(k)}\alpha_{\max}\|\mathbf p_h\|_{L^2(\Omega)}\|v_h\|_{1,h}\\
&\ge\frac{\alpha_{\min}}2\|\mathbf p_h\|_{L^2(\Omega)}^2
+\frac{{\delta(k)}{\beta_B(k)}}2\|v_h\|_{1,h}^2.
\end{aligned}
\]
Moreover, \eqref{eq:sdg-b-test-function} gives
$\|\Psi_h\|_{*,h}\le(1+{\delta(k)})\|\Phi_h\|_{*,h}$.
Applying \eqref{eq:degree-inverse-trace} to the first term of the lower
bound and dividing by $\|\Psi_h\|_{*,h}$ gives
\eqref{eq:weak-form-infsup} with
\[
\begin{aligned}
\beta_A(k)^{-1}
&=(1+\delta(k))\max\left\{\frac{2+(k+1)^2}{\alpha_{\min}},
\frac{2}{\delta(k)\beta_B(k)}\right\}\\
&=\left(1+\frac{\alpha_{\min}}{8\alpha_{\max}^2k^2}\right)
\max\left\{\frac{2+(k+1)^2}{\alpha_{\min}},
\frac{128\alpha_{\max}^2k^4}{\alpha_{\min}}\right\},
\end{aligned}
\]
and \eqref{eq:weak-form-infsup-degree} follows from
$2+(k+1)^2\le6k^2\le6k^4$ for $k\ge1$.

For (b), let $\Phi_h=(v_h,\mathbf p_h)\in X_h^k$ and
$\gamma:=v_h|_{\partial\Omega}$. Then
$\Phi_{h,0}:=\Phi_h-(\mathcal E_h\gamma,\mathbf0)\in X_{h,0}^k$ and
$\|(\mathcal E_h\gamma,\mathbf0)\|_{*,h}\le(1+C_{\mathcal E}(k))|\gamma|_{\partial,h}$. Since
$\Phi_{h,0}\in X_{h,0}^k$, the definition of $\mathcal A_0$ gives
$\mathcal A(\Phi_{h,0},\Psi_h)=\mathcal A_0(\Phi_h,\Psi_h)
-\mathcal A_0((\mathcal E_h\gamma,\mathbf0),\Psi_h)$ for $\Psi_h\in Y_{h,0}^k$, and
\eqref{eq:weak-form-continuity} bounds the last term by
$C_A(1+C_{\mathcal E}(k))|\gamma|_{\partial,h}\|\Psi_h\|_{*,h}$. Part (a) applied to
$\Phi_{h,0}$ and the triangle inequality
$\|\Phi_h\|_{*,h}\le\|\Phi_{h,0}\|_{*,h}+\|(\mathcal E_h\gamma,\mathbf0)\|_{*,h}$ give
\begin{equation}
\label{eq:trace-splitting-estimate}
\|\Phi_h\|_{*,h}
\le\beta_A(k)^{-1}\sup_{0\ne\Psi_h\in Y_{h,0}^k}
\frac{|\mathcal A_0(\Phi_h,\Psi_h)|}{\|\Psi_h\|_{*,h}}
+\bigl(1+C_{\mathcal E}(k)\bigr)\bigl(1+C_A\beta_A(k)^{-1}\bigr)|\gamma|_{\partial,h}.
\end{equation}
Let $S$ denote the supremum in
\eqref{eq:extended-weak-form-infsup}. On the uniform partition, $h_e$ is the same
for all edges, so the test triple $(\mathbf0,0,h_e^{-1}\gamma)$ belongs to
$Y_h^k$. It gives
$\mathcal A_\partial(\Phi_h,(\mathbf0,0,h_e^{-1}\gamma))
=|\gamma|_{\partial,h}^2$ and has norm $|\gamma|_{\partial,h}$, so
$|\gamma|_{\partial,h}\le S$. For $\Psi_h=(\mathbf z_h,w_h)\in Y_{h,0}^k$,
let
$\mu_h:=-Q_\partial(\mathbf z_h\cdot\mathbf n)$. Since
$\gamma\in U_h^k|_{\partial\Omega}$, the integral
$\int_{\partial\Omega}\gamma\,\mu_h\,d\sigma$ cancels the boundary integral in
the definition of $\mathcal A_0$, so
$\mathcal A_\partial(\Phi_h,(\Psi_h,\mu_h))=\mathcal A_0(\Phi_h,\Psi_h)$.
The projection $Q_\partial$ does not increase
$\sum_{e\in\mathcal F_\partial}h_e\|\cdot\|_{L^2(e)}^2$, so
\eqref{eq:boundary-normal-trace} gives
$\|(\Psi_h,\mu_h)\|_{*,h}\le(k+2)\|\Psi_h\|_{*,h}$, and hence
$|\mathcal A_0(\Phi_h,\Psi_h)|\le(k+2)S\|\Psi_h\|_{*,h}$. Substituting this
bound and $|\gamma|_{\partial,h}\le S$ into
\eqref{eq:trace-splitting-estimate} gives
\eqref{eq:extended-weak-form-infsup-degree}.
\end{proof}
From now on, $C_A$, $\beta_A(k)$, and $\beta_\partial(k)$ denote fixed constants for which Propositions~\ref{prop:weak-form-continuity} and~\ref{prop:weak-form-infsup} hold with these bounds.

\begin{corollary}[Bounds on the singular values]
\label{cor:singular-value-bounds}
For every $k\ge1$, the singular values in
Theorem~\ref{prop:scalar-loss-identity} satisfy
\[
\sigma_{\min}\ge\frac{\beta_\partial(k)}{\max(1,\alpha_{\max}(k+2))},
\qquad
\sigma_{\max}\le{(C_A+k+2)}\max(1,\alpha_{\min}^{-1}).
\]
In particular, $\sigma_{\min}$ has a positive lower bound independent of
$h$.
\end{corollary}
\begin{proof}
The bounds on $\alpha$ and \eqref{eq:boundary-normal-trace} give
\[
\begin{aligned}
\|(\Psi_h,\mu_h)\|_{*,h}
&\le\max(1,\alpha_{\min}^{-1})\|(\Psi_h,\mu_h)\|_{Y_h^k},\\
\|(\Psi_h,\mu_h)\|_{Y_h^k}
&\le\max\bigl(1,\alpha_{\max}(k+2)\bigr)\|(\Psi_h,\mu_h)\|_{*,h}.
\end{aligned}
\]
Part (b) of Proposition~\ref{prop:weak-form-infsup} gives the lower
bound. For $\Phi_h=(v_h,\mathbf p_h)\in X_h^k$ and
$(\Psi_h,\mu_h)\in Y_h^k$ with $\Psi_h=(\mathbf z_h,w_h)$,
\[
\mathcal A_\partial\bigl(\Phi_h,(\Psi_h,\mu_h)\bigr)
=\mathcal A_0(\Phi_h,\Psi_h)
+\int_{\partial\Omega}v_h\,\mathbf z_h\cdot\mathbf n\,d\sigma
+\int_{\partial\Omega}v_h\,\mu_h\,d\sigma,
\]
and \eqref{eq:weak-form-continuity}, the Cauchy--Schwarz inequality,
and \eqref{eq:boundary-normal-trace} bound the three terms by $C_A$, $k+1$, and
$1$ times $\|\Phi_h\|_{*,h}\|(\Psi_h,\mu_h)\|_{*,h}$, which gives the upper
bound.
\end{proof}

\begin{corollary}[Quasi-optimality of the SDG projection]
\label{cor:sdg-projection-approximation}
For $\Phi\in\mathbb V\times{\mathbb W}$ and
$\Phi_h\in X_h^k$, let $\gamma$ denote the trace on
$\partial\Omega$ of the scalar component of $\Phi-\Phi_h$. Then
\begin{equation}
\label{eq:sdg-projection-boundary-split}
\|\Phi-{\Pi_h^k}\Phi\|_{*,h}
\le K_0(k)\|\Phi-\Phi_h\|_{*,h}
+K_\partial(k)|\gamma|_{\partial,h},
\end{equation}
where
\[
K_0(k):=1+\frac{C_A}{\beta_A(k)},
\qquad
K_\partial(k):=1+C_{\mathcal E}(k)
+\beta_A(k)^{-1}\bigl(C_A(1+C_{\mathcal E}(k))+k+1\bigr).
\]
Consequently,
\begin{equation}
\label{eq:sdg-projection-best-approximation}
\|\Phi-{\Pi_h^k}\Phi\|_{*,h}
\le K(k)
\inf_{\Phi_h\in X_h^k}\|\Phi-\Phi_h\|_{*,h},
\end{equation}
with the stability factor ${K(k)}:=K_0(k)+K_\partial(k)$. Moreover,
\begin{equation}
\label{eq:degree-projection-stability}
K_0(k)\le(1+C_AC_\beta)k^4,
\qquad
K_\partial(k)\le\bigl(C_\beta(6C_A+2)+6\bigr)k^6,
\end{equation}
so $K(k)\le C_Kk^6$ with $C_K:=C_\beta(7C_A+2)+7$, which is
independent of $h$ and $k$ for fixed coefficient bounds.
\end{corollary}
\begin{proof}
Let $\Xi:=\Phi-\Phi_h=(v,\mathbf p)$, so that $\gamma=v|_{\partial\Omega}$.
Since $\Pi_h^k$ reproduces $X_h^k$, we have
$\Phi-\Pi_h^k\Phi=\Xi-\Pi_h^k\Xi$. The moment equations in
\eqref{eq:sdg-projection} show that the scalar component of $\Pi_h^k\Xi$ has
the trace $q:=Q_\partial\gamma$. On the uniform partition,
$Q_\partial$ does not increase $|\cdot|_{\partial,h}$, so
$|q|_{\partial,h}\le|\gamma|_{\partial,h}$ and
$|\gamma-q|_{\partial,h}\le|\gamma|_{\partial,h}$.
For $\Psi_h=(\mathbf z_h,w_h)\in Y_{h,0}^k$, \eqref{eq:sdg-projection} with
$\mu_h=0$ and the definition of $\mathcal A_0$ give
\[
\mathcal A_0(\Pi_h^k\Xi,\Psi_h)
=\mathcal A_0(\Xi,\Psi_h)
+\int_{\partial\Omega}(\gamma-q)\,\mathbf z_h\cdot\mathbf n\,d\sigma,
\]
and \eqref{eq:weak-form-continuity}, the Cauchy--Schwarz inequality,
and \eqref{eq:boundary-normal-trace} bound the absolute value of the right-hand
side by $\bigl(C_A\|\Xi\|_{*,h}+(k+1)|\gamma|_{\partial,h}\bigr)\|\Psi_h\|_{*,h}$.
The estimate \eqref{eq:trace-splitting-estimate} in the proof of
Proposition~\ref{prop:weak-form-infsup}, applied to $\Pi_h^k\Xi$, and
$\|\Xi-\Pi_h^k\Xi\|_{*,h}\le\|\Xi\|_{*,h}+\|\Pi_h^k\Xi\|_{*,h}$ give
\eqref{eq:sdg-projection-boundary-split}. The boundary term in
$\|\Xi\|_{*,h}$ gives $|\gamma|_{\partial,h}\le\|\Xi\|_{*,h}$, so
\eqref{eq:sdg-projection-boundary-split} implies
\eqref{eq:sdg-projection-best-approximation} after taking the infimum over
$\Phi_h$. Finally, \eqref{eq:weak-form-infsup-degree}, $C_{\mathcal E}(k)\le5k^2$, and
$k+1\le2k^2$ give \eqref{eq:degree-projection-stability}{, and
adding the two bounds gives $K(k)\le C_Kk^6$}.
\end{proof}

The following is the standard quasi-optimality estimate
(cf.~\cite[Theorem~1.35]{DE12}) applied to the exact solution.
\begin{corollary}[Quasi-optimality of the SDG finite element approximation]
\label{cor:sdg-fe-quasi-optimality}
Let $(u,\mathbf p)$ be the exact solution with boundary data $g$, and
suppose that $\mathbf p|_\tau\in H^1(\tau)^2$ for every $\tau\in\mathcal T$.
Then ${\Pi_h^k}(u,\mathbf p)=({u_{h,g}^k},{\mathbf p_{h,g}^k})$ and
\begin{equation}
\label{eq:sdg-fe-quasi-optimality}
{D_h(k)}
\le K(k)
\inf_{\Phi_h\in X_h^k}
\|(u,\mathbf p)-\Phi_h\|_{*,h}.
\end{equation}
\end{corollary}
\begin{proof}
The pair $(u,\mathbf p)$ belongs to
$\mathbb V\times{\mathbb W}$:
$u\in H^1(\Omega)\subset\mathbb V$, and
$\mathbf p\in H(\operatorname{div};\Omega)$, since \eqref{eq:weak-formulation} gives $\nabla\cdot\mathbf p=cu-f-\nabla\cdot(\mathbf bu)\in L^2(\Omega)$.
With the elementwise $H^1$-regularity,
the normal component of $\mathbf p$ is therefore continuous across every
interior edge, so $\mathbf p\in\mathbb W$.
Since $u$ and the normal components of $\mathbf p$ and $\mathbf bu$ are
single-valued across every interior edge, integration by parts on each
subrectangle and the first-order system of
Section~\ref{sec:scalar-elliptic-problem} give
$\mathcal A((u,\mathbf p),\Psi_h)=F_g(\Psi_h)$ for every $\Psi_h\in {Y_{h,0}^k}$.
Since $u|_{\partial\Omega}=g$, it follows that
$\mathcal A_\partial((u,\mathbf p),(\Psi_h,\mu_h))=F_{g,\partial}(\Psi_h,\mu_h)$
for every $(\Psi_h,\mu_h)\in Y_h^k$.
By \eqref{eq:sdg-projection}, ${\Pi_h^k}(u,\mathbf p)$ therefore solves
\eqref{eq:sdg-boundary-problem}, and uniqueness gives
${\Pi_h^k}(u,\mathbf p)=({u_{h,g}^k},{\mathbf p_{h,g}^k})$.
Now \eqref{eq:sdg-fe-quasi-optimality} follows from
\eqref{eq:sdg-projection-best-approximation} and
$\|\cdot\|_h\le\|\cdot\|_{*,h}$.
\end{proof}
The assumption on $\mathbf p$ holds, for example, when
$a\in W^{1,\infty}(\Omega)$ and $u\in H^2(\Omega)$.

\subsection{Degree convergence of the \texorpdfstring{network}{network} projection error}
\label{subsec:degree-convergence}

We show that, for a fixed neural network
trial pair, the network projection error $P_h$ in
Theorem~\ref{prop:scalar-error-estimate} tends to zero as the degree $k$ of the SDG spaces increases.
For a hidden unit
$a_{P,\ell}\tanh(\mathbf w_{P,\ell}^{\mathsf T}\boldsymbol\xi+b_{P,\ell})$ of a patch
network from Section~\ref{subsec:subrectangle-local-interior-networks}, we write
$W_{P,\ell}:=\|\mathbf w_{P,\ell}\|_1$, the sum of the absolute values of the components of $\mathbf w_{P,\ell}$. The unit is nonconstant on $\widehat P$
if and only if $a_{P,\ell}\ne0$ and $W_{P,\ell}>0$, since $\tanh$ is injective. For $W_{P,\ell}>0$, we write
$R_{P,\ell}:=\exp(\operatorname{arsinh}(\pi/(4W_{P,\ell})))$. The value
$R_{P,\ell}>1$ is the parameter of a Bernstein ellipse, the region bounded by the ellipse with foci
$\pm1$ whose semiaxes sum to $R_{P,\ell}$. Its imaginary semiaxis is $\pi/(4W_{P,\ell})$, and the profile
$t\mapsto\tanh(b_{P,\ell}+W_{P,\ell}t)$, with
$t=\mathbf w_{P,\ell}^{\mathsf T}\boldsymbol\xi/W_{P,\ell}\in[-1,1]$, extends analytically to its closure and is
bounded by one there. It differs from $\rho_\ell$ of
Appendix~\ref{subsec:weight-quadrature}, which is determined by the nearest
singularities in the coordinates of a subrectangle. The proof in Appendix~\ref{app:chebyshev} combines
Corollary~\ref{cor:sdg-projection-approximation} with a comparison pair obtained
by truncating the Chebyshev series of each hidden unit.

\begin{theorem}[Degree convergence for a fixed neural network trial pair]
\label{thm:fixed-network-degree-convergence}
Let the mesh, the coefficient functions, and all network parameters
be fixed. Suppose that at least one hidden unit is nonconstant, and let $R_*$
denote the minimum of $R_{P,\ell}$ over the nonconstant hidden units of all
networks for $u$, $p_x$, and $p_y$. Then
$R_*>1$, and there exists a constant $C>0$,
independent of $k$, such that

\begin{equation}
\label{eq:fixed-network-degree-convergence}
{P_h(k;\theta,\eta)}
\le C k^6R_*^{-k}
\qquad\text{for every }k\ge1.
\end{equation}
In particular, $P_h(k;\theta,\eta)\to0$ as $k\to\infty$.
If every hidden unit is constant, then
$P_h(k;\theta,\eta)=0$ for every $k\ge1$.
\end{theorem}

\section{Numerical experiments}
\label{sec:numerical-experiments}
\suppressfloats[t]
\par
We test the SDGNN method on four problems: a smooth Poisson problem, a convection--diffusion problem
with an internal layer, and two elliptic transmission problems with prescribed scalar jumps, across a
square interface and across a circular interface.
Section~\ref{subsec:numerical-settings} defines the problems and
the computational settings. Section~\ref{subsec:numerical-loss-error}
compares the full error with the loss and with the terms of
Theorem~\ref{prop:scalar-error-estimate}, and
Section~\ref{subsec:numerical-pinn} compares the SDGNN method with other
neural network solvers. Section~\ref{subsec:numerical-interfaces} treats the
interface problems, and Section~\ref{subsec:numerical-mesh-degree-budget}
examines how the errors depend on the mesh, on the polynomial degree, and on
the number of training iterations.
\par

\subsection{Problems and computational settings}
\label{subsec:numerical-settings}
In Examples~1--3, $\Omega=(0,1)^2$, and we use uniform Cartesian partitions with $N$ primal cells in each
coordinate direction and cell size $h=1/N$.

\textbf{Example 1.} We consider the Poisson problem with $a=1$,
$\mathbf b=\mathbf0$, and $c=0$.
The exact solution and forcing are
\[
u(x,y)=\sin(\pi x)\sin(\pi y),
\qquad
f(x,y)=2\pi^2\sin(\pi x)\sin(\pi y),
\]
with homogeneous Dirichlet boundary conditions. The diffusive variable is
\[
\mathbf p(x,y)=
\begin{pmatrix}
\pi\cos(\pi x)\sin(\pi y)\\
\pi\sin(\pi x)\cos(\pi y)
\end{pmatrix}.
\]

\textbf{Example 2.} We consider the convection--diffusion problem
\[
\begin{aligned}
-\varepsilon\Delta u-\partial_xu-\partial_yu&=0
&&\text{in }\Omega,\\
u(x,y)&=\tanh\!\left(\frac{y-x}{\varepsilon}\right)
&&\text{on }\partial\Omega,
\end{aligned}
\]
with $\varepsilon=0.1,0.03,0.01$. In the notation of the model problem,
$a=\varepsilon$, $\mathbf b=(1,1)^{\mathsf T}$, and $c=0$, so the diffusive
variable is $\mathbf p=\varepsilon\nabla u$. The boundary data generate
a transition near the diagonal $y=x$.

\textbf{Example 3.} We consider a transmission problem with
$a=1$, $\mathbf b=\mathbf 0$, and $c=0$ and a
prescribed scalar jump across $\Gamma=\partial D$, the boundary of the square inclusion
\[
D=(3/8,7/8)^2.
\]
With $w(x,y)=\sin(\pi x)\sin(\pi y)$ and
\[
v(x,y)=1+\tfrac12\cos\bigl(2\pi(x-3/8)\bigr)\cos\bigl(2\pi(y-3/8)\bigr),
\]the exact fields are
\[
u=w+\chi_Dv,
\qquad
\mathbf p=\nabla w{}+\chi_D\nabla v,
\]
where $\chi_D$ denotes the characteristic function of $D$. In $D$ and $\Omega\setminus\overline D$, these fields satisfy
\[
\mathbf p-\nabla u=\mathbf0,
\qquad
-\nabla\cdot\mathbf p=f,
\]
with $u=0$ on $\partial\Omega$, where
\[
f=2\pi^2w+4\pi^2\chi_D\cos\bigl(2\pi(x-3/8)\bigr)\cos\bigl(2\pi(y-3/8)\bigr).
\]
Since the normal derivative of $v$ vanishes on $\Gamma$, the interface conditions are
\[
[u]_\Gamma=v|_\Gamma,
\qquad
[\mathbf p\cdot\mathbf n]_\Gamma=0,
\]
where $\mathbf n$ is the outward unit normal to $D$, and each jump is
computed as the trace from inside $D$ minus the trace from outside $D$.
The prescribed jump varies along $\Gamma$ between $1/2$ and $3/2$.

We impose the prescribed scalar jump through the load functional
\[
F_\Gamma(\Psi_h)=F_0(\Psi_h)+\int_\Gamma [u]_\Gamma\,\mathbf z_h\cdot\mathbf n\,d\sigma,
\qquad\Psi_h=(\mathbf z_h,w_h)\in Y_{h,0}^k,
\]
where $F_0$ is the load functional in \eqref{eq:scalar-load-functional}
with $g=0$; the weak forms $\mathcal A$ and $\mathcal A_\partial$ are
unchanged. The interface integral enters only the constitutive residual, since
$[\mathbf p\cdot\mathbf n]_\Gamma=0$. The trial functions contain no fixed
function with the prescribed jump, so the networks must produce the jump
themselves. A constant jump would not test this, since adding it to the output
biases of the scalar networks on the patches in $D$ turns every trial pair of
Example 1 into one of this problem with the same residuals and errors. The
jump $v|_\Gamma$ varies along $\Gamma$ and cannot be produced in this way.

We use $N=4,12,20$, for which $\Gamma$ lies on $\mathcal F_p$ and the
scalar trial functions admit independent traces on the two sides of $\Gamma$.
On these meshes, the exact pair $(u,\mathbf p)$ belongs to
$\mathbb V\times{\mathbb W}$, and integration by parts on each
subrectangle gives $\mathcal A((u,\mathbf p),\Psi_h)=F_\Gamma(\Psi_h)$ for every
$\Psi_h\in Y_{h,0}^k$. Theorems~\ref{prop:scalar-loss-identity}
and~\ref{prop:scalar-error-estimate} and
Corollary~\ref{cor:sdg-fe-quasi-optimality} therefore hold for Example~3 when
$F_g$ is replaced by $F_\Gamma$ in the SDG finite element approximation and in
the residual.
On these aligned meshes, we measure the relative error in the
prescribed scalar jump by
\[
E_{[u]}=\frac{\|[u_\theta]_\Gamma-v\|_{L^2(\Gamma)}}{\|v\|_{L^2(\Gamma)}}.
\]
We evaluate the two scalar traces at the same points on $\Gamma$, one
from each adjacent subrectangle. The normal component of the diffusive
variable is continuous across $\mathcal F_p$ by construction.

We use the quadrature of Appendix~\ref{subsec:weight-quadrature} with
tolerance $10^{-10}$ and $q_{\max}=81$. Every 50 L-BFGS iterations, we
predict the orders again, which gives $\tilde q_\tau$. If $\tilde q_\tau>q_\tau$ on
some subrectangle, we replace every order $q_\tau$ by
$\max\{q_\tau,\min(\tilde q_\tau+2,q_{\max})\}$ and restart the L-BFGS memory. At the end of
every run, the loss evaluated with $q_\tau=81$ on all subrectangles differs
from the training loss by a relative amount below $6\times10^{-11}$.

We train all networks for 4,000 iterations of L-BFGS~\cite{1989Liu} with a
history size of 10, normalizing each objective by its initial value.
We compute a sparse Cholesky factorization of $M_Y$ once before
training, and each evaluation of the loss and its gradient solves one system
with $M_Y$. We train
each configuration from five random initializations. Unless stated
otherwise, the SDG spaces have degree $k=4$, the SDGNN networks have
local hidden width $r=6$, and the tables give the median of each quantity over
these initializations, with its range in parentheses. At the end of
the 4,000-iteration runs, the largest training quadrature order
$\max_\tau q_\tau$ lies between 14 and 16 in Example~1, between 14 and 62 in
Example~2, and between 15 and 22 in Example~3.

We measure the accuracy of the scalar and diffusive variables by their relative
$L^2(\Omega)$ errors
\[
E_u=\frac{\|u-u_\theta\|_{L^2(\Omega)}}{\|u\|_{L^2(\Omega)}},
\qquad
E_{\mathbf p}=
\frac{\|\mathbf p-{\mathbf p_\eta}\|_{L^2(\Omega)}}
{\|\mathbf p\|_{L^2(\Omega)}}.
\]
We also measure the full error
$E:=\|(u,\mathbf p)-(u_\theta,\mathbf p_\eta)\|_h$ in the norm
\eqref{eq:full-error-norm}, the left-hand side of
Theorem~\ref{prop:scalar-error-estimate}; unlike $E_u$ and $E_{\mathbf p}$, it
is not normalized.
In Examples 1 and 3, we evaluate these errors using the analytical
expressions for the exact scalar and diffusive fields.
In Example 2, we replace the exact fields in these definitions by
reference fields computed with the conforming Galerkin method with
continuous piecewise ${\mathbb Q}_{12}$ polynomials, on the uniform mesh of $64\times64$
squares for $\varepsilon=0.01$ and of $48\times48$ squares for
$\varepsilon=0.1,0.03$. The boundary values of this approximation interpolate
$g$ at the Gauss--Lobatto--Legendre nodes, and the reference diffusive field is
$\varepsilon$ times its gradient. Its relative $L^2$ differences from the
approximation of degree 14 on the same mesh and from that of degree 12 on the
mesh with $3/2$ times as many squares in each direction are below
$3\times10^{-10}$ for the scalar variable and $1.4\times10^{-7}$ for the
diffusive variable. In $E$, the exact solution has no jumps
across $\mathcal F_p$ and equals $g$ on $\partial\Omega$, so we use the
reference fields only for $\nabla u$ and $\mathbf p$ on the subrectangles.
We evaluate the reference fields directly at the Gauss--Legendre
quadrature points used for all field-error integrals.

We implement PINN and Mixed-PINN using DeepXDE~\cite{2021Lu} with the PyTorch
backend~\cite{2019Paszke}. SDGNN uses the Hager--Zhang line
search~\cite{2006Hager}, and the PINN methods use the strong-Wolfe line
search~\cite[Section~3.1]{2006Nocedal}. SDGNN runs on one CPU thread,
and PINN and Mixed-PINN run on a GPU, with several runs sharing each device, so
we do not compare elapsed times.

We compare the methods at two parameter levels, the larger given by SDGNN with
$N=6$, $r=6$ and the smaller by SDGNN with $N=2$, $r=12$. Each patch network has $4r+1$
parameters, so an SDGNN trial pair with $N$ primal cells in each coordinate
direction has $\bigl((N+1)^2+4N^2\bigr)(4r+1)$ parameters. PINN uses one network
with two tanh hidden layers of width 67 or 33. Mixed-PINN uses three separate
networks, each with two tanh hidden layers of width 38 or 18, for the scalar
variable and the two components of the diffusive variable. At each level, these
are the smallest widths for which the parameter count is at least that of
SDGNN (Table~\ref{tab:smooth-methods}). The PINN
interior loss is the mean squared residual of the second-order equation, and the
Mixed-PINN interior loss is the sum of the mean squared residuals of the
first-order system of Section~\ref{sec:scalar-elliptic-problem}. Both methods
evaluate these residuals at the same fixed set of 65,536 interior Sobol points
and add, with unit weight, the mean squared boundary discrepancy at 256 midpoint
samples on each side of the boundary. SDGNN imposes the boundary data through the
boundary residual $R_\partial$, whose contribution to the loss
\eqref{eq:sdgnn-loss} is fixed by the test inner product in
\eqref{eq:residual-metric-matrix} rather than by a separate weight.
We impose the Dirichlet data by boundary moments; with the data imposed strongly, the scalar approximations
of Example~2 have larger errors near the corner $(1,1)$, where the internal layer meets the boundary.

\par
\textbf{Example 4.} We consider the transmission problem of~\cite[Example~6]{2023Tseng} with a circular
interface. Let $\Omega=(-1,1)^2$, let $D=\{\mathbf x\in\Omega:|\mathbf x|<2/3\}$, and let $\Gamma=\partial D$.
The coefficients are $a=2$ in $D$ and $a=3$ in $\Omega\setminus\overline D$, $\mathbf b=\mathbf 0$, and
$c=0$, and the exact scalar field is
\[
u=\sin(4\pi x)\sin(4\pi y)+7\quad\text{in }D,
\qquad
u=5e^{-x^2-y^2}\quad\text{in }\Omega\setminus\overline D,
\]
with $\mathbf p=a\nabla u$. The forcing $f$, the boundary data $g$, and the jumps $[u]_\Gamma$ and
$[\mathbf p\cdot\mathbf n]_\Gamma$, taken from inside $D$ minus outside $D$ as in Example~3, are those of
this field.

\begin{figure}[!htbp]
\centering
\begin{tikzpicture}[scale=2.5,
  fpface/.style={draw=black!45, line width=0.35pt, dash pattern=on 1.6pt off 1.2pt},
  fuface/.style={draw=black, line width=0.7pt},
  bdface/.style={draw=black, line width=1.1pt},
  ifface/.style={draw=red!80!black, line width=1.5pt}]
\input{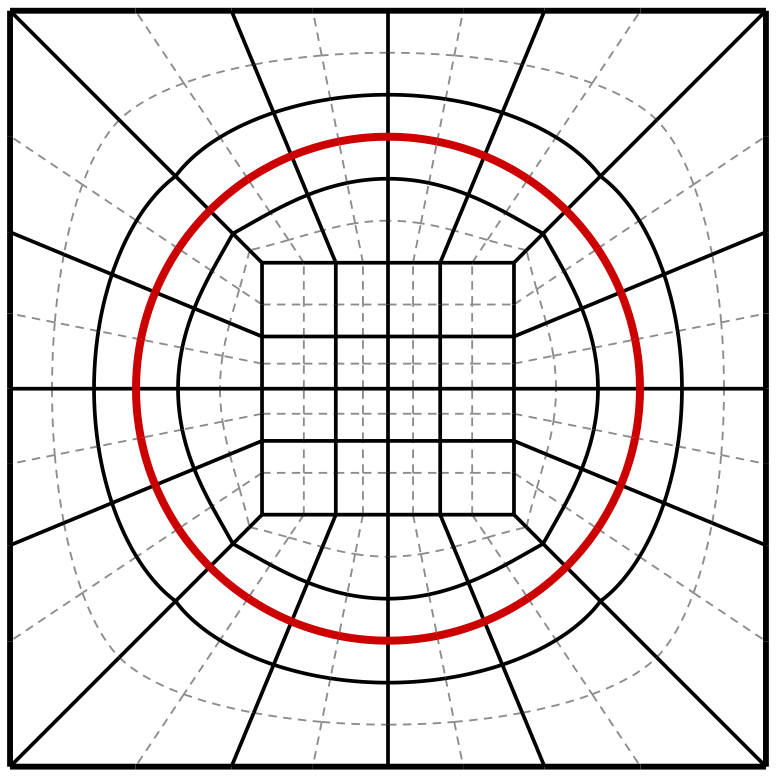}
\end{tikzpicture}
\caption{Example 4: the O-grid. Solid black lines: edges in $\mathcal F_u$; dashed gray
lines: edges in $\mathcal F_p$; red: the interface $\Gamma$, which consists of 32 edges in
$\mathcal F_p$.}
\label{fig:ogrid-mesh}
\end{figure}

Since the scalar trial functions may jump only across $\mathcal F_p$, we use an O-grid on which $\Gamma$
consists of edges in $\mathcal F_p$ (Figure~\ref{fig:ogrid-mesh}): a $4\times4$ central block and a ring of
$16\times3$ primal cells, 64 in all, each the image of the reference square under a smooth map; the refined
mesh, with an $8\times8$ central block and a ring of $32\times5$ primal cells, has 224. The SDG spaces are defined through these maps, the scalar networks are evaluated
at the physical coordinates, and the diffusive variable is the contravariant Piola transform of the
reference field given by the networks of $\mathcal P_{p_x}$ and $\mathcal P_{p_y}$, so that
$\mathbf p_\eta\cdot\mathbf n$ stays continuous across $\mathcal F_p$. The prescribed scalar jump enters the
load functional as in Example~3, and the prescribed jump of $\mathbf p\cdot\mathbf n$ enters through a fixed
vector field added to $\mathbf p_\eta$.

We use $k=4$, the widths $r=1,2,4,6$ on the O-grid and $r=1,2$ on the refined mesh, and five
initializations, with the training settings of Examples~1--3 except for the quadrature: since
Appendix~\ref{subsec:weight-quadrature} assumes affine subrectangles, we use fixed Gauss--Legendre rules with
16 points per coordinate, and doubling the points changes the final losses by at most $2.5\times10^{-6}$ in
relative terms. We train for 16,000 L-BFGS iterations, the limit of the continued runs in
Section~\ref{subsec:numerical-iteration-limit}; $E_{[u]}$ is defined as in Example~3, with $[u]_\Gamma$ in place of $v$. We compare SDGNN with
the cusp-capturing PINN of~\cite{2023Tseng}, in which a first network (400 parameters) is fitted to the
scalar jump and carries the jump of the approximation, and a second network (775 parameters) with the
additional input $|x^2+y^2-4/9|$ is then trained on a collocation loss. Both stages use the
Levenberg--Marquardt method as implemented by the authors of~\cite{2023Tseng}, with at most 10,000 updates
per stage.
\par

\FloatBarrier
\subsection{The loss as an error measure}
\label{subsec:numerical-loss-error}

\begin{figure}[!htbp]
\centering
\includegraphics[width=\linewidth]{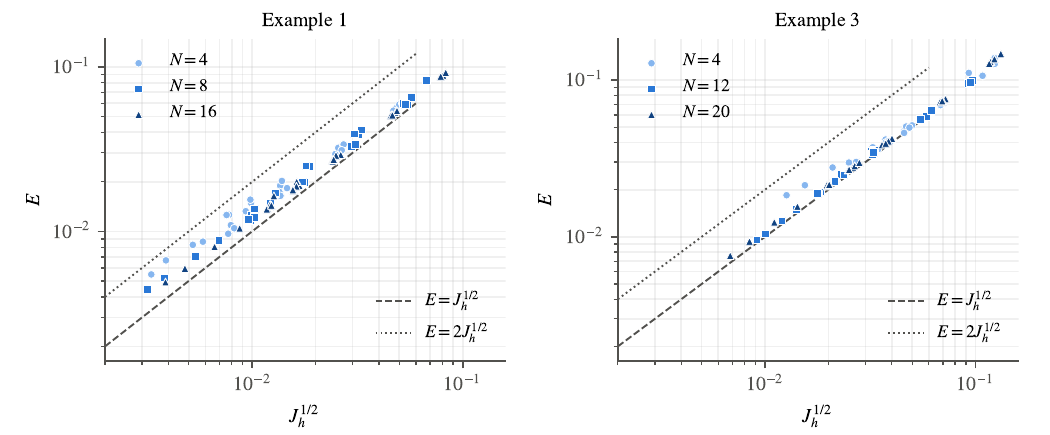}
\caption{Examples 1 and 3: full error against $J_h^{1/2}$ at iterations
500, 1,000, 2,000, 3,000, and 4,000 of all runs and at 6,000, 8,000, 12,000, and
16,000 of the runs continued in
Section~\ref{subsec:numerical-iteration-limit}, with the lines $E=J_h^{1/2}$ and
$E=2J_h^{1/2}$.}
\label{fig:loss-error-ratio}
\end{figure}

For Examples~1 and~3, the ratio of the full error to $J_h^{1/2}$ remains
between $0.97$ and $1.72$ from the 500th iteration on, independently of the mesh
(Figure~\ref{fig:loss-error-ratio}), so that the square root of the loss
serves as a measure of the error. To see which term dominates, we evaluate the terms
of the bound in Theorem~\ref{prop:scalar-error-estimate} for Example~1.
We write this bound as $E\le D_h+\sigma_{\min}^{-1}J_h^{1/2}+P_h$, where
$D_h$ and $P_h$ are the discretization error
\eqref{eq:sdg-discretization-error-degree} and the network projection error
\eqref{eq:neural-sdg-projection-difference}, $J_h$
is the loss \eqref{eq:sdgnn-loss}, and $\sigma_{\min}$ is the
constant in Theorem~\ref{prop:scalar-loss-identity}. We compute
$\sigma_{\min}^{-2}$ as the largest eigenvalue of
${A}_h^{-1}M_Y{A}_h^{-\mathsf T}M_X$ by the Lanczos method in the $M_X$
inner product, with a sparse LU factorization of ${A}_h$.
For each mesh, we
compute the SDG finite element approximation and the SDG projection of the
final SDGNN approximation.
Table~\ref{tab:smooth-error-components} gives the full error and the terms
$D_h$, $J_h^{1/2}$, and $P_h$ after 4,000 L-BFGS iterations for $N=4,8,16$.
\begin{table}[!htbp]
\centering\small
\caption{Example 1: full error $E$, the terms $D_h$, $J_h^{1/2}$, and $P_h$, and
the bound $D_h+\sigma_{\min}^{-1}J_h^{1/2}+P_h$ with
$\sigma_{\min}=5.445\times10^{-2}$ after 4,000 L-BFGS
iterations. The entries for $E$, $J_h^{1/2}$, $P_h$, and the bound are medians over the initializations, taken separately for each column.}
\label{tab:smooth-error-components}
\setlength{\tabcolsep}{4pt}
\begin{tabular}{r c c c c c}
\hline
$N$ & $E$ & $D_h$ & $J_h^{1/2}$ & $P_h$ & Bound\\
\hline
4 & $1.091\times10^{-2}$ & $1.810\times10^{-5}$ & $7.765\times10^{-3}$ & $2.078\times10^{-3}$ & $1.457\times10^{-1}$\\
8 & $1.214\times10^{-2}$ & $1.134\times10^{-6}$ & $9.957\times10^{-3}$ & $2.562\times10^{-3}$ & $1.858\times10^{-1}$\\
16 & $1.441\times10^{-2}$ & $7.091\times10^{-8}$ & $1.217\times10^{-2}$ & $2.899\times10^{-3}$ & $2.263\times10^{-1}$\\
\hline
\end{tabular}
\end{table}

On each mesh, $D_h$ is more than $500$ times smaller than $E$, and
$P_h$ is between $9\%$ and $36\%$ of $E$ for every initialization. After 4,000
L-BFGS iterations, the full error is therefore mostly the training error, which
Theorem~\ref{prop:scalar-loss-identity} bounds by
$\sigma_{\min}^{-1}J_h^{1/2}$. The loss thus measures the error that the
training can still reduce; $D_h$ does not depend on the training, and it
is negligible here. The computed value
$\sigma_{\min}=5.445\times10^{-2}$ does not change across the three meshes, consistent with the
mesh-independent lower bound of Corollary~\ref{cor:singular-value-bounds}, and
the bound $D_h+\sigma_{\min}^{-1}J_h^{1/2}+P_h$ exceeds $E$ by factors between $11$
and $16$ for every initialization.

\begin{figure}[!htbp]
\centering
\includegraphics[width=\linewidth,height=0.64\textheight,keepaspectratio]{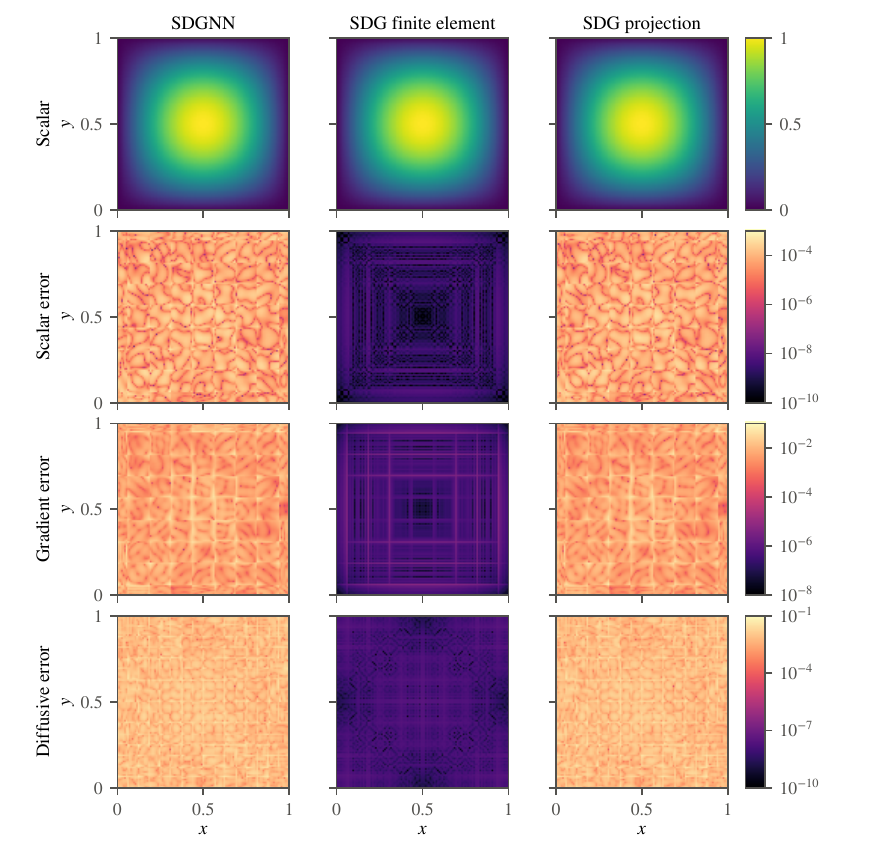}
\caption{Example 1: scalar approximations and pointwise errors for $N=8$,
from the initialization with the median $E_u$.
The columns show the final SDGNN approximation, the SDG finite element
approximation, and the SDG projection. The rows show the scalar
approximation, absolute scalar error, magnitude of the piecewise scalar
gradient error, and magnitude of the diffusive-variable error.
A common color scale is used in each row. The error scales are logarithmic,
with lower display limits of $10^{-10}$, $10^{-8}$, and $10^{-10}$, respectively.}
\label{fig:smooth-projection-fields}
\end{figure}

Figure~\ref{fig:smooth-projection-fields} shows the approximations and the
pointwise errors for $N=8$. The three scalar approximations cannot be told
apart at the scale of the plot. The error maps of the SDGNN approximation and
of its SDG projection are nearly identical, and the errors of the SDG finite
element approximation are several orders of magnitude smaller.

The ratio $E/J_h^{1/2}$ of Figure~\ref{fig:loss-error-ratio}, recorded
from iteration 500 on, is thus far below $\sigma_{\min}^{-1}\approx18$: the square root of the
loss estimates the full error more closely than the bound guarantees. In Example 2,
the computed $\sigma_{\min}$ is about $0.14\varepsilon$
to $0.18\varepsilon$, from $1.4\times10^{-2}$ for $\varepsilon=0.1$ to
between $1.5\times10^{-3}$ and $1.8\times10^{-3}$ for
$\varepsilon=0.01$. After 4,000 iterations, the difference between the SDG
projection and the SDG finite element approximation is comparable to $E$,
and $J_h^{1/2}$ underestimates this difference by factors between about $10$ and $200$, as the
small $\sigma_{\min}$ permits; the ratio $E/J_h^{1/2}$ lies between $\varepsilon^{-1}$ and
$3\varepsilon^{-1}$, and between $10$ and $290$ over the three values of
$\varepsilon$.

\subsection{\texorpdfstring{Comparison with other neural network solvers}{Comparison with other neural network solvers}}
\label{subsec:numerical-pinn}

For Example 1, we compare SDGNN, PINN, and Mixed-PINN at the two
parameter levels of Section~\ref{subsec:numerical-settings} (Table~\ref{tab:smooth-methods}).

\begin{table}[!htbp]
\centering\small
\caption{Example 1: parameter counts and relative $L^2$ errors after 4,000
L-BFGS iterations for SDGNN, PINN, and Mixed-PINN at two parameter levels. The
errors are medians with ranges in parentheses.}
\label{tab:smooth-methods}
\setlength{\tabcolsep}{4pt}
\begin{tabular}{l r c c}
\hline
Method & Parameters & $E_u$ & $E_{\mathbf p}$\\
\hline
\multicolumn{4}{l}{Larger parameter level}\\
SDGNN ($N=6$, $r=6$) & 4,825 & $2.56\,(2.51\text{--}3.14)\times10^{-4}$ & $3.19\,(3.07\text{--}3.26)\times10^{-3}$\\
PINN & 4,825 & $1.20\,(0.73\text{--}2.89)\times10^{-3}$ & $3.23\,(2.13\text{--}8.22)\times10^{-3}$\\
Mixed-PINN & 4,905 & $2.70\,(2.29\text{--}2.96)\times10^{-3}$ & $6.94\,(6.01\text{--}8.48)\times10^{-3}$\\
\multicolumn{4}{l}{Smaller parameter level}\\
SDGNN ($N=2$, $r=12$) & 1,225 & $5.19\,(4.72\text{--}5.84)\times10^{-4}$ & $2.94\,(2.61\text{--}3.06)\times10^{-3}$\\
PINN & 1,255 & $2.25\,(0.94\text{--}2.84)\times10^{-3}$ & $6.58\,(2.65\text{--}8.69)\times10^{-3}$\\
Mixed-PINN & 1,245 & $4.40\,(2.18\text{--}4.96)\times10^{-3}$ & $9.64\,(6.67\text{--}12.63)\times10^{-3}$\\
\hline
\end{tabular}
\end{table}

At both parameter levels, SDGNN attains the smallest scalar error for
every initialization in Table~\ref{tab:smooth-methods}; the median scalar
errors of PINN and Mixed-PINN are 4.3 to 10.6 times that of SDGNN. The median
diffusive-variable errors of SDGNN and PINN are close at the larger level, and
that of PINN is 2.2 times that of SDGNN at the smaller level.
Between 2,000 and 4,000 L-BFGS iterations, the median errors of PINN change by
less than 4\%, whereas those of SDGNN and Mixed-PINN decrease by factors of 1.9 to 2.4
and 1.9 to 2.7, respectively.

For Example 2 with $\varepsilon=0.01$, we compare SDGNN with PINN
and Mixed-PINN at the two parameter levels of Table~\ref{tab:smooth-methods},
with the settings of Section~\ref{subsec:numerical-settings}. Here the
first-order system reads
\[
\varepsilon^{-1}\mathbf p-\nabla u=\mathbf0,
\qquad
-\nabla\cdot\mathbf p-\partial_xu-\partial_yu=0,
\]
so the Mixed-PINN loss measures the constitutive residual on the scale
of $\varepsilon^{-1}(\mathbf p-\varepsilon\nabla u)$. The SDGNN loss, by contrast, measures it
on the scale of $\mathbf p-\varepsilon\nabla u$ (Section~\ref{subsec:matrix-representation-norm}).

\begin{table}[!htbp]
\centering\small\setlength{\tabcolsep}{3pt}
\caption{Example 2: relative $L^2$ errors and root-mean-square (RMS) errors of the scalar boundary values after
4,000 L-BFGS iterations at $\varepsilon=0.01$, at the two parameter levels of
Table~\ref{tab:smooth-methods}. The boundary errors are evaluated at the midpoints of 4,096 equal
subintervals of each side, none of which is a training point. The errors are medians with ranges in
parentheses.}
\label{tab:diagonal-methods}
\begin{tabular}{l c c c}
\hline
Method & $E_u$ & $E_{\mathbf p}$ & Boundary RMS\\
\hline
\multicolumn{4}{l}{Larger parameter level}\\
SDGNN & $1.66\,(1.52\text{--}3.05)\times10^{-2}$ & $1.96\,(1.80\text{--}2.20)\times10^{-1}$ & $3.05\,(0.80\text{--}12.2)\times10^{-3}$\\
PINN & $7.55\,(7.50\text{--}7.57)\times10^{-2}$ & $5.41\,(5.37\text{--}5.42)\times10^{-1}$ & $1.87\,(1.86\text{--}1.87)\times10^{-1}$\\
Mixed-PINN & $1.49\,(1.43\text{--}1.60)\times10^{-1}$ & $7.13\,(7.08\text{--}7.25)\times10^{-1}$ & $2.76\,(2.73\text{--}2.86)\times10^{-1}$\\
\multicolumn{4}{l}{Smaller parameter level}\\
SDGNN & $5.91\,(3.59\text{--}6.59)\times10^{-2}$ & $3.99\,(3.82\text{--}4.25)\times10^{-1}$ & $3.51\,(2.81\text{--}4.72)\times10^{-2}$\\
PINN & $7.55\,(7.53\text{--}7.57)\times10^{-2}$ & $5.35\,(5.33\text{--}5.36)\times10^{-1}$ & $1.86\,(1.86\text{--}1.86)\times10^{-1}$\\
Mixed-PINN & $1.30\,(1.19\text{--}1.52)\times10^{-1}$ & $6.98\,(6.85\text{--}7.15)\times10^{-1}$ & $2.65\,(2.53\text{--}2.79)\times10^{-1}$\\
\hline
\end{tabular}
\end{table}

At both parameter levels, SDGNN attains the smallest scalar and
diffusive-variable errors for every initialization in
Table~\ref{tab:diagonal-methods}. The PINN errors are nearly the
same at the two parameter levels and for all initializations, and they change little
after about 1,000 L-BFGS iterations. The scalar boundary errors of
PINN and Mixed-PINN are much larger than those of SDGNN at both levels, so the
field comparison includes the effects of the boundary
treatments as well as the trial functions and residual losses.

\par
Table~\ref{tab:structure-comparison} adds three trial functions trained on the SDGNN loss: the PINN
network, with $\mathbf p_\eta=a\nabla u_\theta$; the three Mixed-PINN networks; and, for each of $u$, $p_x$,
and $p_y$, a shared network that all patches use with the center of the
patch as an additional input. They are integrated with fixed Gauss--Legendre rules, with which SDGNN gives
errors within $13\%$ of those in the table.

\begin{table}[!htbp]
\centering\small\setlength{\tabcolsep}{3.5pt}
\caption{Examples 1 and 2 ($\varepsilon=0.01$): median relative $L^2$ errors after 4,000 L-BFGS
iterations at the two parameter levels of Section~\ref{subsec:numerical-settings}. The collocation rows and the
patch-network rows repeat Tables~\ref{tab:smooth-methods} and~\ref{tab:diagonal-methods}.}
\label{tab:structure-comparison}
\begin{tabular}{l l r c c c c}
\hline
 & & & \multicolumn{2}{c}{Example 1} & \multicolumn{2}{c}{Example 2}\\
Trial functions & Loss & Parameters & $E_u$ & $E_{\mathbf p}$ & $E_u$ & $E_{\mathbf p}$\\
\hline
\multicolumn{7}{l}{Larger parameter level ($N=6$)}\\
PINN network & collocation & 4,825 & $1.20\times10^{-3}$ & $3.23\times10^{-3}$ & $7.55\times10^{-2}$ & $5.41\times10^{-1}$\\
Mixed-PINN networks & collocation & 4,905 & $2.70\times10^{-3}$ & $6.94\times10^{-3}$ & $1.49\times10^{-1}$ & $7.13\times10^{-1}$\\
PINN network & SDGNN & 4,825 & $6.69\times10^{-4}$ & $2.86\times10^{-3}$ & $1.07\times10^{-1}$ & $9.71\times10^{-1}$\\
Mixed-PINN networks & SDGNN & 4,905 & $2.60\times10^{-3}$ & $1.00\times10^{-2}$ & $1.80\times10^{-1}$ & $7.71\times10^{-1}$\\
Shared networks & SDGNN & 4,887 & $2.31\times10^{-3}$ & $9.95\times10^{-3}$ & $1.59\times10^{-1}$ & $7.74\times10^{-1}$\\
Patch networks (SDGNN) & SDGNN & 4,825 & $2.56\times10^{-4}$ & $3.19\times10^{-3}$ & $1.66\times10^{-2}$ & $1.96\times10^{-1}$\\
\multicolumn{7}{l}{Smaller parameter level ($N=2$)}\\
PINN network & collocation & 1,255 & $2.25\times10^{-3}$ & $6.58\times10^{-3}$ & $7.55\times10^{-2}$ & $5.35\times10^{-1}$\\
Mixed-PINN networks & collocation & 1,245 & $4.40\times10^{-3}$ & $9.64\times10^{-3}$ & $1.30\times10^{-1}$ & $6.98\times10^{-1}$\\
PINN network & SDGNN & 1,255 & $2.78\times10^{-4}$ & $1.44\times10^{-3}$ & $6.70\times10^{-2}$ & $6.32\times10^{-1}$\\
Mixed-PINN networks & SDGNN & 1,245 & $1.41\times10^{-3}$ & $5.65\times10^{-3}$ & $1.35\times10^{-1}$ & $6.69\times10^{-1}$\\
Shared networks & SDGNN & 1,227 & $1.39\times10^{-3}$ & $5.99\times10^{-3}$ & $1.41\times10^{-1}$ & $6.64\times10^{-1}$\\
Patch networks (SDGNN) & SDGNN & 1,225 & $5.19\times10^{-4}$ & $2.94\times10^{-3}$ & $5.91\times10^{-2}$ & $3.99\times10^{-1}$\\
\hline
\end{tabular}
\end{table}

Among the trial functions trained on the SDGNN loss, the patch networks are the most accurate in
Example~2, and in Example~1 at the larger parameter level they have the smallest scalar error. In
Example~1, the PINN network trained on the SDGNN loss is more accurate at the smaller level. The
differences from the collocation rows also contain different optimizer implementations and integration
points, so the table does not isolate the effect of the loss.

In Example~3, the PINN and Mixed-PINN networks are continuous and cannot represent the prescribed jump;
the Mixed-PINN networks trained on the SDGNN loss give $E_{[u]}=1$. The shared
networks can jump across $\mathcal F_p$, but their errors are 5.6 to 8.0 times those of SDGNN in
Table~\ref{tab:jump-errors} at $N=4$, and 28 to 160 times at $N=12$ and $N=20$, where SDGNN has more
parameters.
\par

\subsection{\texorpdfstring{Interfaces}{Interfaces}}
\label{subsec:numerical-interfaces}
\label{subsec:numerical-accuracy}

In Example~3, the prescribed jump lies on the midlines of the primal
cells, and the scalar trial functions represent it through their independent
traces on the two sides of the interface, without a penalty term.
Table~\ref{tab:jump-errors} compares the relative $L^2$ errors of the
scalar and diffusive variables and of the scalar jump for
SDGNN and the SDG finite element approximations.

\begin{table}[!htbp]
\centering\small\setlength{\tabcolsep}{5pt}
\caption{Example 3: relative $L^2$ errors, including that of the scalar jump,
for the SDGNN approximations after 4,000 L-BFGS iterations (medians with ranges
in parentheses) and the SDG finite element approximations, on meshes aligned with $\Gamma$.}
\label{tab:jump-errors}
\begin{tabular}{c l c c c}
\hline
$N$ & Method & $E_u$ & $E_{\mathbf p}$ & $E_{[u]}$\\\hline
4 & SDGNN & $7.80\,(6.57\text{--}8.11)\times10^{-4}$ & $1.05\,(0.92\text{--}1.11)\times10^{-2}$ & $9.80\,(8.47\text{--}11.4)\times10^{-4}$\\
4 & SDG-FE & $8.378\times10^{-7}$ & $5.682\times10^{-6}$ & $2.637\times10^{-6}$\\
12 & SDGNN & $2.17\,(2.13\text{--}2.29)\times10^{-4}$ & $5.67\,(5.59\text{--}5.88)\times10^{-3}$ & $2.42\,(2.19\text{--}2.59)\times10^{-4}$\\
12 & SDG-FE & $3.505\times10^{-9}$ & $3.785\times10^{-8}$ & $6.627\times10^{-9}$\\
20 & SDGNN & $4.01\,(3.82\text{--}4.30)\times10^{-4}$ & $6.38\,(6.34\text{--}6.53)\times10^{-3}$ & $1.57\,(1.42\text{--}1.70)\times10^{-4}$\\
20 & SDG-FE & $2.732\times10^{-10}$ & $3.736\times10^{-9}$ & $4.745\times10^{-10}$\\\hline
\end{tabular}
\end{table}

The SDGNN approximations capture the prescribed jump on the three meshes,
with the median jump error decreasing as $N$ increases.
The SDG finite element approximations yield smaller scalar,
diffusive-variable, and jump errors on the three meshes.
At the fixed iteration limit, the scalar and diffusive-variable errors
decrease from $N=4$ to $N=12$ and increase from $N=12$ to $N=20$;
Section~\ref{subsec:numerical-iteration-limit} examines these errors at a larger iteration limit.

\begin{figure}[!htbp]
\centering
\includegraphics[width=\linewidth]{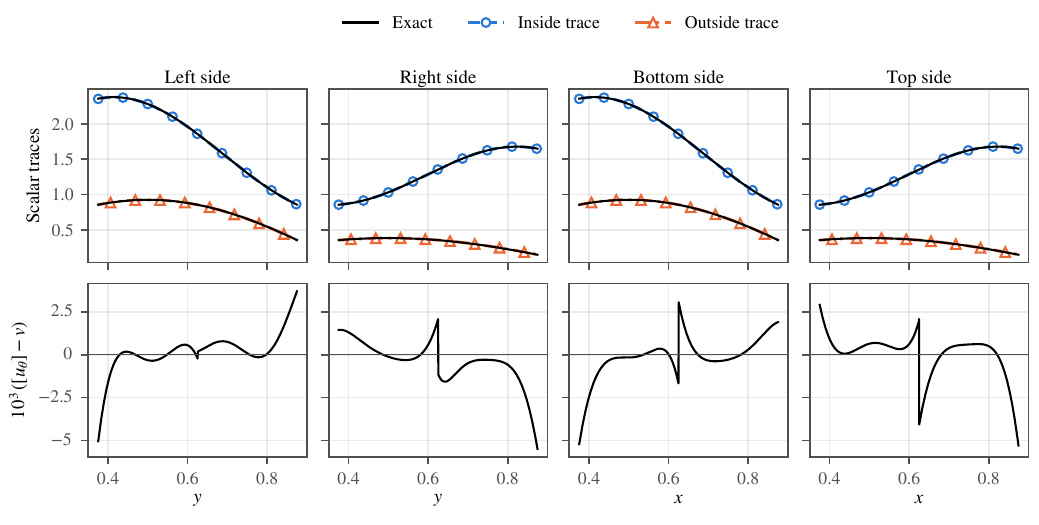}
\caption{Example 3: scalar traces and jump error for the final SDGNN
approximation with $N=4$, from the initialization with the median
$E_u$. The traces are evaluated separately on the
two sides of the interface.}
\label{fig:jump-traces}
\end{figure}

Figure~\ref{fig:jump-traces} shows the scalar traces on both sides of
the interface and the error in the prescribed jump for $N=4$.

\par
\begin{figure}[!htbp]
\centering
\includegraphics[width=\linewidth]{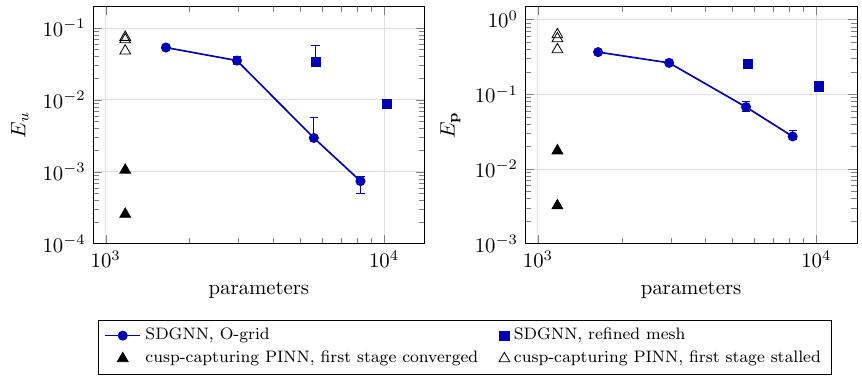}
\caption{Example 4: errors against the number of parameters. SDGNN: medians over five initializations
after 16,000 L-BFGS iterations, with bars from the smallest to the largest value. Cusp-capturing
PINN: each initialization after at most 10,000 updates per stage.}
\label{fig:ogrid-errors}
\end{figure}

In Example~4, $\Gamma$ consists of edges in $\mathcal F_p$ of the O-grid, so the scalar trial functions represent the prescribed
jump through their independent traces on the two sides of $\Gamma$, as in Example~3, and the jump is
captured consistently over the initializations: with $r=6$ on the O-grid, the median $E_{[u]}$ after 16,000
iterations is $2.1\times10^{-3}$, and in every SDGNN setting each error varies by a factor of at most 2.2 over the
five initializations (Table~\ref{tab:ogrid-sdgnn}).
The cusp-capturing PINN, which builds the interface into its network, is more accurate per parameter when
its first stage converges (Table~\ref{tab:ogrid-cusp}, Figure~\ref{fig:ogrid-errors}): SDGNN needs about seven times its 1,175 parameters to reach a comparable $E_u$.
Its first stage, however, converged for two of the five initializations; for the other three, the error in
the jump remains between $18\%$ and $24\%$. The SDGNN errors are far above those of the SDG finite element
approximation on the O-grid ($E_u=1.5\times10^{-5}$), and the analysis of
Sections~\ref{sec:weak-residual-loss} and~\ref{sec:projection-error-estimate} does not cover this mapped
mesh.
\par

\begin{table}[!htbp]
\centering\small\setlength{\tabcolsep}{4pt}
\caption{Example 4: relative $L^2$ errors, including that of the scalar jump, of the SDGNN approximations
after 16,000 L-BFGS iterations; medians over five initializations with ranges in parentheses.}
\label{tab:ogrid-sdgnn}
\begin{tabular}{l c r c c c}
\hline
Mesh & $r$ & Parameters & $E_u$ & $E_{\mathbf p}$ & $E_{[u]}$\\
\hline
O-grid & 1 & 1,645 & $5.35\,(5.27\text{--}5.87)\times10^{-2}$ & $3.66\,(3.57\text{--}3.79)\times10^{-1}$ & $1.04\,(0.89\text{--}1.10)\times10^{-1}$\\
O-grid & 2 & 2,961 & $3.52\,(3.12\text{--}4.03)\times10^{-2}$ & $2.63\,(2.51\text{--}2.80)\times10^{-1}$ & $4.46\,(3.62\text{--}5.04)\times10^{-2}$\\
O-grid & 4 & 5,593 & $2.95\,(2.64\text{--}5.63)\times10^{-3}$ & $6.73\,(5.93\text{--}8.01)\times10^{-2}$ & $7.28\,(5.59\text{--}8.78)\times10^{-3}$\\
O-grid & 6 & 8,225 & $7.42\,(5.05\text{--}8.71)\times10^{-4}$ & $2.72\,(2.50\text{--}3.31)\times10^{-2}$ & $2.08\,(1.85\text{--}3.94)\times10^{-3}$\\
refined & 1 & 5,685 & $3.41\,(3.15\text{--}5.80)\times10^{-2}$ & $2.51\,(2.46\text{--}2.74)\times10^{-1}$ & $2.88\,(2.60\text{--}3.03)\times10^{-2}$\\
refined & 2 & 10,233 & $8.85\,(8.55\text{--}9.31)\times10^{-3}$ & $1.27\,(1.19\text{--}1.29)\times10^{-1}$ & $1.35\,(1.23\text{--}1.47)\times10^{-2}$\\
\hline
\end{tabular}
\end{table}

\begin{table}[!htbp]
\centering\small\setlength{\tabcolsep}{4pt}
\caption{Example 4: errors of the cusp-capturing PINN with 1,175 parameters for each initialization,
ordered by the final loss of the first stage, after at most 10,000 Levenberg--Marquardt updates per stage.}
\label{tab:ogrid-cusp}
\begin{tabular}{c c c c c}
\hline
Initialization & Loss of first stage & $E_u$ & $E_{\mathbf p}$ & $E_{[u]}$\\
\hline
A & $7.5\times10^{-11}$ & $2.57\times10^{-4}$ & $3.22\times10^{-3}$ & $2.26\times10^{-6}$\\
B & $9.1\times10^{-11}$ & $1.06\times10^{-3}$ & $1.75\times10^{-2}$ & $2.48\times10^{-6}$\\
C & $0.45$ & $4.80\times10^{-2}$ & $3.95\times10^{-1}$ & $1.75\times10^{-1}$\\
D & $0.77$ & $6.88\times10^{-2}$ & $6.28\times10^{-1}$ & $2.29\times10^{-1}$\\
E & $0.87$ & $7.43\times10^{-2}$ & $5.55\times10^{-1}$ & $2.42\times10^{-1}$\\
\hline
\end{tabular}
\end{table}

\subsection{\texorpdfstring{Mesh, degree, and iteration limit}{Mesh, degree, and iteration limit}}
\label{subsec:numerical-mesh-degree-budget}

\subsubsection*{Mesh refinement}

For Example~2, Table~\ref{tab:diagonal-errors} gives the relative scalar
and diffusive-variable errors of SDGNN and of the SDG finite element
approximations for $N=4,8,16$ against the reference fields in
Section~\ref{subsec:numerical-settings}. The finite element columns give the discretization error on these meshes, measured by the
relative $L^2$ errors.

\begin{table}[!htbp]
\centering\small\setlength{\tabcolsep}{4pt}
\caption{Example 2: relative $L^2$ errors of SDGNN after 4,000 L-BFGS
iterations and of the degree-4 SDG finite element approximations
defined by \eqref{eq:sdg-boundary-problem} on the same meshes.
The SDGNN errors are medians with ranges in parentheses. Errors are evaluated against
the reference fields in Section~\ref{subsec:numerical-settings}.}
\label{tab:diagonal-errors}
\begin{tabular}{c r c c c c}
\hline
$\varepsilon$ & $N$ & $E_u$ (SDGNN) & $E_{\mathbf p}$ (SDGNN) & $E_u$ (SDG-FE) & $E_{\mathbf p}$ (SDG-FE)\\
\hline
0.1 & 4 & $1.95\,(1.49\text{--}2.06)\times10^{-3}$ & $1.57\,(1.32\text{--}1.67)\times10^{-2}$ & $3.370\times10^{-5}$ & $4.026\times10^{-4}$\\
0.1 & 8 & $7.61\,(7.32\text{--}8.73)\times10^{-4}$ & $9.57\,(9.00\text{--}9.77)\times10^{-3}$ & $2.022\times10^{-6}$ & $2.839\times10^{-5}$\\
0.1 & 16 & $5.37\,(5.02\text{--}5.83)\times10^{-4}$ & $8.93\,(8.41\text{--}9.39)\times10^{-3}$ & $8.008\times10^{-8}$ & $1.569\times10^{-6}$\\
0.03 & 4 & $8.50\,(6.28\text{--}9.37)\times10^{-3}$ & $5.49\,(5.36\text{--}6.25)\times10^{-2}$ & $1.123\times10^{-3}$ & $1.282\times10^{-2}$\\
0.03 & 8 & $2.94\,(2.72\text{--}3.67)\times10^{-3}$ & $3.38\,(3.24\text{--}3.54)\times10^{-2}$ & $9.476\times10^{-5}$ & $2.045\times10^{-3}$\\
0.03 & 16 & $1.82\,(1.66\text{--}2.01)\times10^{-3}$ & $2.69\,(2.56\text{--}2.80)\times10^{-2}$ & $4.253\times10^{-6}$ & $1.473\times10^{-4}$\\
0.01 & 4 & $3.97\,(2.77\text{--}5.77)\times10^{-2}$ & $2.70\,(2.21\text{--}2.97)\times10^{-1}$ & $9.949\times10^{-3}$ & $1.053\times10^{-1}$\\
0.01 & 8 & $1.98\,(1.46\text{--}2.16)\times10^{-2}$ & $1.60\,(1.46\text{--}1.74)\times10^{-1}$ & $1.578\times10^{-3}$ & $2.534\times10^{-2}$\\
0.01 & 16 & $6.63\,(5.26\text{--}7.22)\times10^{-3}$ & $9.93\,(8.81\text{--}10.2)\times10^{-2}$ & $1.444\times10^{-4}$ & $5.301\times10^{-3}$\\
\hline
\end{tabular}
\end{table}

For each diffusion coefficient, the median SDGNN errors of both
fields decrease across
$N=4,8,16$ after 4,000 L-BFGS iterations. At each $N$, the errors are
larger for smaller diffusion coefficients.
The discretization error decreases faster with $N$ than the SDGNN error:
for $\varepsilon=0.01$, the median SDGNN
scalar error is 4.0 times the finite element error at $N=4$ and 46 times at
$N=16$, so that at this iteration limit the SDGNN error is limited by the training rather
than by the discretization.

\subsubsection*{Degree}
\label{subsec:numerical-degree}

\begin{figure}[!htbp]
\centering
\includegraphics[width=\linewidth]{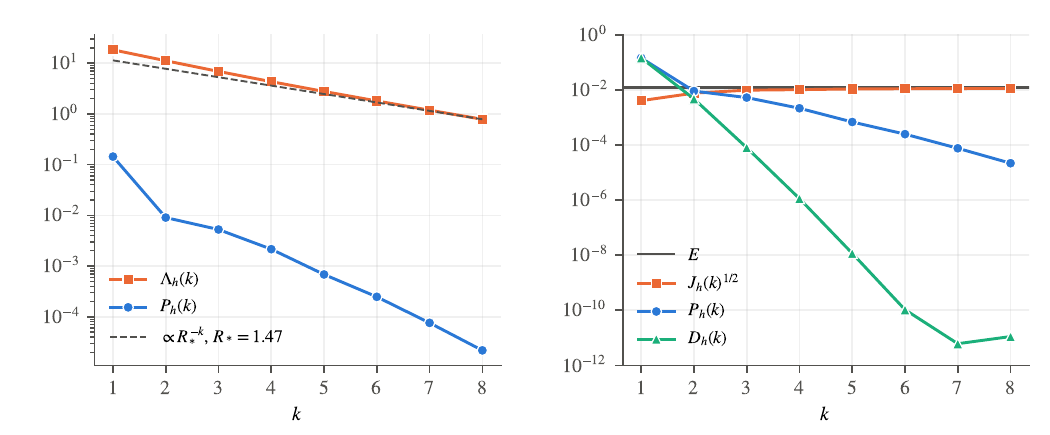}
\caption{Example 1 with $N=8$, one trained network. Left: the
network projection error $P_h(k)$ and the quantity $\Lambda_h(k)$ of
\eqref{eq:chebyshev-global-error-bound}, which bounds the distance between the
network and its Chebyshev comparison pair, against the degree $k$ of the SDG spaces, with a
reference line proportional to $R_*^{-k}$. Right: the full error $E$, the
discretization error $D_h(k)$, $J_h(k)^{1/2}$, and $P_h(k)$.}
\label{fig:degree-convergence}
\end{figure}

To examine Theorem~\ref{thm:fixed-network-degree-convergence}, we fix
the network of Example 1 with $N=8$, trained with $k=4$, for which $E_u$ is the median over
the initializations, and we compute the terms of
Theorem~\ref{prop:scalar-error-estimate} for the SDG spaces of degree
$k=1,\dots,8$. For each $k$, we use
the quadrature orders predicted for degree $k$, and raising every order by 6
leaves $P_h$ unchanged in the first six digits.
Figure~\ref{fig:degree-convergence} shows the results. The network projection error
decreases from $1.4\times10^{-1}$ at $k=1$ to $2.2\times10^{-5}$ at $k=8$, by
factors between $2.4$ and $3.5$ per degree for $k\ge3$. The quantity
$\Lambda_h(k)$ of Appendix~\ref{app:chebyshev} bounds the distance between the
network and its Chebyshev comparison pair; it decreases by factors between
$1.5$ and $1.7$ per degree, approaches the rate $R_*^{-k}$ with
$R_*=1.47$, and exceeds $P_h$ by factors between $1.3\times10^{2}$ and
$3.6\times10^{4}$. The bound $Ck^6R_*^{-k}$ of
Theorem~\ref{thm:fixed-network-degree-convergence} increases with $k$ for
$k<6/\log R_*\approx15.7$, so Figure~\ref{fig:degree-convergence} shows how
$P_h$ and $\Lambda_h$ decrease, not how the proven bound decreases. The
discretization error $D_h(k)$ reaches $10^{-11}$ at $k=7$. The full error $E$
does not depend on $k$, and $J_h(k)^{1/2}$ increases with $k$, from $34\%$ of
$E$ at $k=1$ to $93\%$ at $k=8$. Since $Y_h^k\subset Y_h^{k+1}$ and the test inner product does not
depend on $k$, $J_h(k)$ cannot decrease as $k$ increases, so the loss measures a
larger part of the error as the test space grows.

Training the networks of Example 1 with $N=8$ for $k=2,\dots,6$ instead of
$k=4$, with the other settings unchanged, gives a median full error $9.4$ times
larger for $k=2$ and $1.33$ times larger for $k=3$ than for $k=4$, and $12\%$
and $13\%$ smaller for $k=5$ and $k=6$. For $k=2$, $P_h$ is about as large as
$E$, and the median of $E/J_h^{1/2}$ is $15.7$: with a small test space, most of
the error lies in the part that the loss does not see. For $k\ge4$, the median
of $E/J_h^{1/2}$ is between $1.07$ and $1.25$. Over the same degrees,
$\sigma_{\min}$ decreases from $5.4\times10^{-2}$ to $2.9\times10^{-2}$.

\subsubsection*{Iteration limit}
\label{subsec:numerical-iteration-limit}

At the common limit of 4,000 L-BFGS iterations, the full errors in
Example 1 increase with $N$ (Table~\ref{tab:smooth-error-components}). The SDG finite element errors are much
smaller: $D_h$ is at most $1.9\times10^{-3}$ times $E$ for each initialization in
Table~\ref{tab:smooth-error-components}, and the SDG finite element rows of
Table~\ref{tab:jump-errors} are at least two orders of magnitude below the
SDGNN rows. Since the ratio $E/J_h^{1/2}$ stays between $0.97$ and $1.72$
(Section~\ref{subsec:numerical-loss-error}), the errors at the iteration limit
follow the loss reached within 4,000 iterations.
In Example 1, the median final loss increases from $6.0\times10^{-5}$ to
$1.5\times10^{-4}$ as $N$ increases from 4 to 16 and the number of
network parameters from 2,225 to 32,825. In Example 3, it is $1.0\times10^{-3}$, $3.2\times10^{-4}$,
and $3.9\times10^{-4}$ for $N=4,12,20$, with 2,225, 18,625, and 51,025
parameters. The number of parameters grows like $N^2$, and the patch
networks exchange information only through the loss, so that a fixed number of
iterations can reduce the loss less on a finer mesh.

\begin{figure}[!htbp]
\centering
\includegraphics[width=\linewidth]{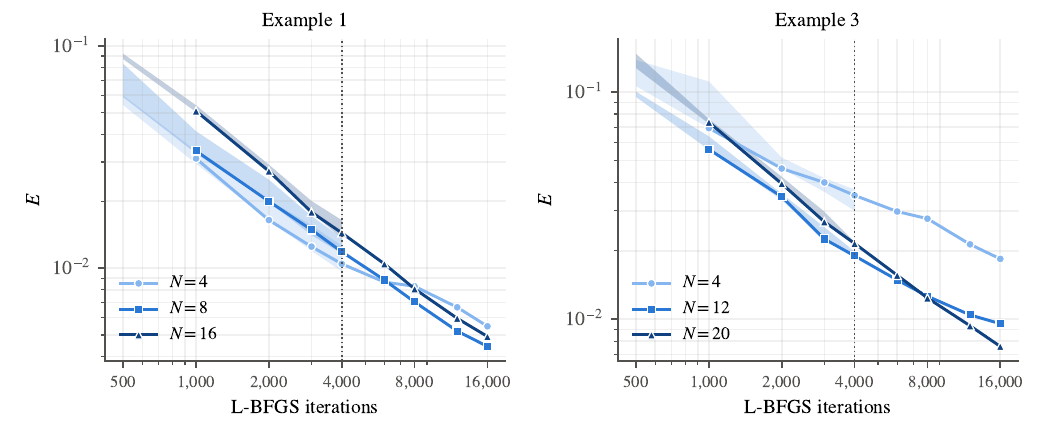}
\caption{Examples 1 and 3: full error against L-BFGS iterations.
The lines show the initialization continued to 16,000 iterations, and the
bands show the range over the initializations up to 4,000 iterations.}
\label{fig:iteration-limit}
\end{figure}

To examine the effect of the iteration limit, we continue one initialization
on each mesh to 16,000 L-BFGS iterations with the other settings unchanged
(Table~\ref{tab:iteration-limit} and Figure~\ref{fig:iteration-limit}). A larger
number of iterations restores the dependence on the mesh. At 16,000 iterations, the scalar error decreases as $N$ increases in
both examples. In Example 1, the diffusive-variable errors and the full errors take similar
values on all meshes; in Example 3, they also decrease with $N$. The full
errors are still decreasing at 16,000 iterations (Figure~\ref{fig:iteration-limit}). On the coarsest
mesh, $P_h$ increases between 4,000 and 16,000 iterations, to $82\%$ of $E$ in
Example 3; on the finer meshes, $P_h$ is between $16\%$ and $38\%$ of $E$ at 16,000
iterations. These results describe errors at
prescribed iteration limits and do not establish a mesh convergence rate for
SDGNN.

\begin{table}[!htbp]
\centering\small\setlength{\tabcolsep}{5pt}
\caption{Examples 1 and 3: errors of one SDGNN initialization after 4,000
and 16,000 L-BFGS iterations, with the other settings of
Tables~\ref{tab:smooth-error-components} and~\ref{tab:jump-errors}.}
\label{tab:iteration-limit}
\begin{tabular}{c r r c c c c}
\hline
Example & $N$ & Iterations & $E_u$ & $E_{\mathbf p}$ & $E$ & $P_h$\\
\hline
1 & 4 & 4,000 & $3.91\times10^{-4}$ & $3.19\times10^{-3}$ & $1.05\times10^{-2}$ & $9.38\times10^{-4}$\\
1 & 4 & 16,000 & $2.01\times10^{-4}$ & $1.32\times10^{-3}$ & $5.47\times10^{-3}$ & $1.77\times10^{-3}$\\
1 & 8 & 4,000 & $1.68\times10^{-4}$ & $3.61\times10^{-3}$ & $1.19\times10^{-2}$ & $2.56\times10^{-3}$\\
1 & 8 & 16,000 & $7.17\times10^{-5}$ & $1.18\times10^{-3}$ & $4.44\times10^{-3}$ & $1.59\times10^{-3}$\\
1 & 16 & 4,000 & $1.97\times10^{-4}$ & $4.65\times10^{-3}$ & $1.44\times10^{-2}$ & $2.51\times10^{-3}$\\
1 & 16 & 16,000 & $3.65\times10^{-5}$ & $1.47\times10^{-3}$ & $4.93\times10^{-3}$ & $1.89\times10^{-3}$\\
3 & 4 & 4,000 & $8.06\times10^{-4}$ & $1.04\times10^{-2}$ & $3.51\times10^{-2}$ & $5.39\times10^{-3}$\\
3 & 4 & 16,000 & $3.29\times10^{-4}$ & $5.43\times10^{-3}$ & $1.84\times10^{-2}$ & $1.50\times10^{-2}$\\
3 & 12 & 4,000 & $2.14\times10^{-4}$ & $5.64\times10^{-3}$ & $1.90\times10^{-2}$ & $3.01\times10^{-3}$\\
3 & 12 & 16,000 & $1.34\times10^{-4}$ & $2.82\times10^{-3}$ & $9.56\times10^{-3}$ & $1.53\times10^{-3}$\\
3 & 20 & 4,000 & $4.04\times10^{-4}$ & $6.50\times10^{-3}$ & $2.15\times10^{-2}$ & $2.72\times10^{-3}$\\
3 & 20 & 16,000 & $6.83\times10^{-5}$ & $2.15\times10^{-3}$ & $7.59\times10^{-3}$ & $2.03\times10^{-3}$\\
\hline
\end{tabular}
\end{table}

\FloatBarrier
\section{\texorpdfstring{Conclusion}{Conclusion}}\label{sec:conclusions}

We have presented the SDGNN method for scalar elliptic problems, in which
neural network trial functions defined on the patches of a staggered mesh
satisfy the continuity conditions of the SDG finite element method for all
values of the network parameters. The SDGNN loss is the squared dual norm of the weak
residual on $Y_h^k$ defined by the test Gram matrix $M_Y$, the
Dirichlet data enter through boundary moments, and the weak residual requires
no derivatives of the networks. We evaluate the weak residual with
Gauss--Legendre quadrature
whose orders are predicted from the network weights during training. Since the
SDG system is square and inf-sup stable, the loss is equivalent to the square
of the training error, the distance between the SDG projection of the network
and the SDG finite element approximation
(Theorem~\ref{prop:scalar-loss-identity}), and the full error is bounded by the
sum of the discretization error, a multiple of the square root of the loss,
and the network projection error (Theorem~\ref{prop:scalar-error-estimate}). On uniform square meshes,
this multiple, $\sigma_{\min}^{-1}$, is bounded independently of the mesh size
for each degree (Corollary~\ref{cor:singular-value-bounds}). On these meshes and for a
fixed network, the network projection error is bounded by $Ck^6R_*^{-k}$, where
$R_*>1$ depends on the network weights
(Theorem~\ref{thm:fixed-network-degree-convergence}).

\par
The numerical experiments are consistent with the analysis. In Examples~1 and~3, the ratio
$E/J_h^{1/2}$ stayed between $0.97$ and $1.72$ from the 500th iteration on, so that
the loss served as a measure of the error. The network projection error of a
trained network decreased from $1.4\times10^{-1}$ at $k=1$ to
$2.2\times10^{-5}$ at $k=8$ in Example~1 with $N=8$. The prescribed jump of
Example~3 was represented on the three meshes without a penalty term, with
errors that were consistent over the initializations,
and the same held for the circular interface of Example~4 on a mapped
mesh. In Example~1 and in Example~2 with $\varepsilon=0.01$, SDGNN attained
smaller scalar errors than PINN and Mixed-PINN at both parameter levels.
\par

\par
The method has the following limitations. The SDG finite element
approximation on the same mesh is more accurate than the SDGNN approximation
and, for a single linear problem,
cheaper. The number of network parameters grows like $N^2$, and at a fixed
number of iterations the error need not decrease with $N$, as in Example~1; a larger number of
iterations is needed on finer meshes. On smooth problems, a single global
network trained on the SDGNN loss attains comparable errors at
the same parameter count.
In Example~2, the ratio $E/J_h^{1/2}$ reached $290$
for $\varepsilon=0.01$, where $\sigma_{\min}$ was between $0.15\varepsilon$
and $0.18\varepsilon$, so that the loss underestimates the error in the
convection-dominated regime. A prescribed jump requires a mesh fitted to the
interface, as in discontinuous Galerkin methods. The analysis is restricted to
linear scalar elliptic problems in two dimensions on rectangular meshes, which excludes the mapped
mesh of Example~4, and the degree-explicit
estimates of Section~\ref{sec:projection-error-estimate} to uniform square
meshes. At a fixed degree, the loss does not control $P_h$
(Remark~\ref{rem:loss-kernel-example}), and the constant in
Theorem~\ref{thm:fixed-network-degree-convergence} depends on the network.
\par

\FloatBarrier
\section*{Acknowledgments}

G.-M.G. was partially supported by the Fulbright U.S. Scholar Program,
which is sponsored by the U.S. Department of State and the
Korean-American Educational Commission.
C.-Y.J. was supported by a grant from the National Research Foundation
of Korea (NRF) funded by the Korea government (MSIT)
(No. 2023R1A2C1003120) and by the 2026 Research Fund (1.260002.01)
of UNIST (Ulsan National Institute of Science \& Technology).
B.K. was supported by the National Research Foundation of Korea (NRF)
grant funded by the Korea government (MSIT) (00560003).
K.K. was supported by the National Research Foundation of Korea (NRF)
grants funded by the Korea government (MSIT)
(Nos. 2023R1A2C1003120 and 00560003).

\appendix
\section[Proof of Lemma~\ref{lem:sdg-b-infsup}]{Proof of Lemma~\ref{lem:sdg-b-infsup}}
\label{app:degree-b-proof}

\begin{proof}
We make the constants in the construction of
\cite[proof of Theorem~1]{KJN21} explicit.
Let $q\in\mathbb P_{k-1}(0,1)$ and $t_k=1/(2k^2)$.
Lemma~\ref{lem:polynomial-trace} with $n=k-1$ and $d=1$ gives
\[
\int_0^{t_k} |q(t)|^2\,dt
\le t_k\|q\|_{L^\infty(0,1)}^2
\le\frac12\|q\|_{L^2(0,1)}^2.
\]
Hence $\int_{t_k}^1|q(t)|^2\,dt\ge\frac12\|q\|_{L^2(0,1)}^2$, and since
$t\ge t_k$ on $(t_k,1)$,
\begin{equation}
\label{eq:weighted-polynomial-lower-bound}
\int_0^1 t|q(t)|^2\,dt
\ge\frac1{4k^2}\int_0^1|q(t)|^2\,dt.
\end{equation}
Reflection gives the same estimate with weight $1-t$.

For $v_h\in U_{h,0}^k$, the gradient and jump contributions to
$\|v_h\|_{1,h}^2$ are
\[
\|v_h\|_\nabla^2=\sum_{\tau\in\mathcal T}\|\nabla v_h\|_{L^2(\tau)}^2,
\qquad
{\|v_h\|_{[\,]}}^2=\sum_{e\in\mathcal F_p}h_e^{-1}\|[v_h]_e\|_{L^2(e)}^2.
\]
The field $\mathbf z_h^{\nabla}\in W_h^k$ from the proof of
Proposition~\ref{prop:scalar-discrete-uniqueness} has zero normal trace
on $\mathcal F_p$.
Its affine weights satisfy $0\le{\omega_{\tau,i}}\le1$.
These bounds and \eqref{eq:weighted-polynomial-lower-bound}, applied in
the normal coordinate of each subrectangle, give
\[
\|\mathbf z_h^{\nabla}\|_{0,h}\le \|v_h\|_\nabla,
\qquad B(\mathbf z_h^{\nabla},v_h)\ge\frac{\|v_h\|_\nabla^2}{4k^2}.
\]
For each $e\in\mathcal F_p$, $\pi_e$ denotes the orthogonal projection
onto the line containing $e$. Since the primal cells are squares of side length $h$, the two
subrectangles adjacent to $e$ are squares of side length $h_e=h/2$.
On these subrectangles, the field $\mathbf s_e$ is defined by
\[
\mathbf s_e(\mathbf x)
=-h_e^{-1}[v_h]_e(\pi_e\mathbf x)
 L_k\!\left(1-\frac{2\operatorname{dist}(\mathbf x,e)}{h_e}\right)
 \mathbf n_e,
\]
and it is zero elsewhere; here $L_k$ is the Legendre polynomial of degree $k$
normalized by $L_k(1)=1$.
The tangential degree of $[v_h]_e$ and the normal degree of $L_k$
are at most $k$, so $\mathbf s_e|_\tau\in {\mathbb Q}_k(\tau)^2$ on each
adjacent subrectangle.
Its normal trace equals $-h_e^{-1}[v_h]_e$ on $e$ and vanishes on every
other edge in $\mathcal F_p$.
Thus $\mathbf s_e\in {W_h^k}$.
Legendre orthogonality gives
$\int_\tau\mathbf s_e\cdot\nabla v_h\,d\mathbf x=0$
on both adjacent subrectangles, since the normal derivative of $v_h$ has
normal degree at most $k-1$.
For $\mathbf z_h^{[\,]}=\sum_{e\in\mathcal F_p}\mathbf s_e$, the vanishing
volume integrals and the prescribed normal traces give
\[
B(\mathbf z_h^{[\,]},v_h)=\|v_h\|_{[\,]}^2.
\]
On each subrectangle, the contributions from its two edges in $\mathcal F_p$
occupy perpendicular vector components.
Since $|L_k|\le1$ on $[-1,1]$, integration over the two
subrectangles adjacent to each edge and summation over edges bound the volume
contribution by $2\|v_h\|_{[\,]}^2$.
The normal trace contribution is $\|v_h\|_{[\,]}^2$, so
$\|\mathbf z_h^{[\,]}\|_{0,h}^2\le3\|v_h\|_{[\,]}^2$.
The field $\mathbf z_{v_h}=(\mathbf z_h^{\nabla}+\mathbf z_h^{[\,]})/2$ therefore satisfies
\[
\begin{aligned}
\|\mathbf z_{v_h}\|_{0,h}
&\le{\frac12\bigl(\|v_h\|_\nabla+\sqrt3\,\|v_h\|_{[\,]}\bigr)}
\le\sqrt{\|v_h\|_\nabla^2+\|v_h\|_{[\,]}^2}=\|v_h\|_{1,h},\\
B(\mathbf z_{v_h},v_h)
&\ge\frac{\|v_h\|_\nabla^2}{8k^2}+\frac{\|v_h\|_{[\,]}^2}{2}
\ge\frac1{8k^2}\|v_h\|_{1,h}^2.
\end{aligned}
\]
These inequalities prove \eqref{eq:sdg-b-test-function} with the value
of $\beta_B(k)$ in \eqref{eq:sdg-b-infsup}. %
The power of $k$ in $\beta_B(k)$ appears not to be sharp: for $N\le4$
and $k\le8$, the computed inf-sup constant in \eqref{eq:sdg-b-infsup} does not depend on
$N$ and agrees with $((2k+1)/(k^2+2k+3))^{1/2}$, which decreases only like $k^{-1/2}$. %
We expect the gap to come from the weights in $\mathbf z_h^{\nabla}$.
\end{proof}

\section{A discrete extension of boundary traces}
\label{app:boundary-trace-extension}

We use the uniform partition and the seminorm
$|\cdot|_{\partial,h}$ of
Section~\ref{subsec:sdg-projection-approximation}.

\begin{lemma}[Discrete extension of boundary traces]
\label{lem:boundary-trace-extension}
For every $k\ge1$, there exists a linear map
$\mathcal E_h:U_h^k|_{\partial\Omega}\to U_h^k$ such that
$(\mathcal E_h\gamma)|_{\partial\Omega}=\gamma$ and
\[
\|\mathcal E_h\gamma\|_{1,h}\le C_{\mathcal E}(k)|\gamma|_{\partial,h},
\qquad
{C_{\mathcal E}}(k)=\tfrac54(k+1)^2,
\]
for every $\gamma\in U_h^k|_{\partial\Omega}$.
\end{lemma}
\begin{proof}
Let $s:=h/2$ denote the side length of the subrectangles, so that $h_e=s$ for every $e\in\mathcal F$. For
$d>0$ and $q\in\mathbb P_k(0,d)$,
\begin{equation}
\label{eq:legendre-derivative}
\|q'\|_{L^2(0,d)}\le\frac{(k+1)^2}{d}\|q\|_{L^2(0,d)}.
\end{equation}
Since $L_n''$ has degree less than $n$, integration by parts gives
$\int_{-1}^1(L_n')^2\,dt=[L_nL_n']_{-1}^{1}=n(n+1)$, so the derivative of the
orthonormal polynomial $\sqrt{(2n+1)/2}\,L_n$ has squared $L^2(-1,1)$ norm
$(2n+1)n(n+1)/2$, and
\[
\sum_{n=0}^k\frac{(2n+1)n(n+1)}2=\frac{k(k+1)^2(k+2)}4\le\frac{(k+1)^4}4.
\]
The Cauchy--Schwarz inequality applied to the expansion of $q$, as
in the proof of Lemma~\ref{lem:polynomial-trace}, and scaling to $(0,d)$ give
\eqref{eq:legendre-derivative}.

Let $\omega(t):=1-t/{s}$ for $0\le t\le{s}$. On a subrectangle
$\tau=(x_0,x_0+{s})\times(0,{s})$ with one side $e$ on the bottom edge of
$\Omega$ and no other side on $\partial\Omega$, we set
$\mathcal E_h\gamma(x,y):=\gamma(x)\omega(y)$. On the corner subrectangle
$(0,{s})^2$, we set
\[
\mathcal E_h\gamma(x,y)
:=\gamma_B(x)\omega(y)+\gamma_L(y)\omega(x)-\gamma(0,0)\omega(x)\omega(y),
\]
where $\gamma_B$ and $\gamma_L$ are the restrictions of $\gamma$ to the
bottom side $e_B$ and the left side $e_L$ of this subrectangle. The other
sides and corners are treated in the same way, and $\mathcal E_h\gamma:=0$ on the
subrectangles with no side on $\partial\Omega$. On each subrectangle $\tau$,
$\mathcal E_h\gamma|_\tau\in {\mathbb Q}_k(\tau)$ for $k\ge1$, and
$(\mathcal E_h\gamma)|_{\partial\Omega}=\gamma$. Across an edge in $\mathcal F_u^0$,
the two traces agree because $\gamma$ is continuous at the primal vertices.
Hence $\mathcal E_h\gamma\in U_h^k$.

On the subrectangle $\tau$ above, \eqref{eq:legendre-derivative}
gives
\[
\|\partial_y\mathcal E_h\gamma\|_{L^2(\tau)}^2={s}^{-1}\|\gamma\|_{L^2(e)}^2,
\qquad
\|\partial_x\mathcal E_h\gamma\|_{L^2(\tau)}^2=\frac{s}{3}\|\gamma'\|_{L^2(e)}^2
\le\frac{(k+1)^4}{3}{s}^{-1}\|\gamma\|_{L^2(e)}^2.
\]
On the corner subrectangle, the same bounds, Lemma~\ref{lem:polynomial-trace} for
$\gamma(0,0)$ on $e_B$ and on $e_L$, and
$(\xi_1+\xi_2+\xi_3)^2\le2\xi_1^2+4\xi_2^2+4\xi_3^2$, applied to each partial derivative with
$\xi_1$ the term that contains $\gamma_B'$ or $\gamma_L'$, give
\[
\|\nabla \mathcal E_h\gamma\|_{L^2((0,{s})^2)}^2
\le\Bigl(\tfrac23(k+1)^4+\tfrac43(k+1)^2+4\Bigr){s}^{-1}
\bigl(\|\gamma_B\|_{L^2(e_B)}^2+\|\gamma_L\|_{L^2(e_L)}^2\bigr).
\]
Across an edge $e\in\mathcal F_p$ with an endpoint on $\partial\Omega$, the
jump of $\mathcal E_h\gamma$ is $[\gamma]_m\,\omega$, where $[\gamma]_m$ is the jump
of $\gamma$ at the midpoint $m$ of a primal boundary edge. Lemma~\ref{lem:polynomial-trace}
gives
\[
h_e^{-1}\bigl\|[\gamma]_m\,\omega\bigr\|_{L^2(e)}^2
=\frac{|[\gamma]_m|^2}{3}
\le\frac23(k+1)^2{s}^{-1}
\bigl(\|\gamma\|_{L^2(e_-)}^2+\|\gamma\|_{L^2(e_+)}^2\bigr),
\]
where $e_\pm$ are the two edges in $\mathcal F$ on $\partial\Omega$ that meet
at $m$. The other edges in $\mathcal F_p$ carry no jump, and each edge in
$\mathcal F$ on $\partial\Omega$ enters one jump term. With
$\kappa:=(k+1)^2\ge4$, the coefficient
$\tfrac23\kappa^2+\tfrac43\kappa+4$ of the corner subrectangle exceeds the
coefficient $1+\tfrac13\kappa^2$ of the other subrectangles, and adding the
coefficient $\tfrac23\kappa$ of the jump term for each boundary edge
gives
\[
\|\mathcal E_h\gamma\|_{1,h}^2
\le\Bigl(\tfrac23\kappa^2+2\kappa+4\Bigr)|\gamma|_{\partial,h}^2
\le\tfrac{17}{12}\kappa^2|\gamma|_{\partial,h}^2,
\]
where the last step uses $2\kappa+4\le\tfrac34\kappa^2$.
Since $\sqrt{17/12}\le5/4$, the lemma follows.
\end{proof}

\section{Chebyshev estimates and proof of
Theorem~\ref{thm:fixed-network-degree-convergence}}\label{app:chebyshev}

To prove Theorem~\ref{thm:fixed-network-degree-convergence}, we bound
the right-hand side of \eqref{eq:sdg-projection-boundary-split}
for a comparison pair built from the network. We first
record a consequence of \cite[Theorems~8.1 and~8.2]{T13} for a single tanh
unit.

\begin{lemma}[Chebyshev approximation of tanh]
\label{lem:tanh-chebyshev-approximation}
For $W>0$ and $b\in\mathbb R$, let $f(t)=\tanh(b+Wt)$ and
$R=\exp(\operatorname{arsinh}(\pi/(4W)))>1$.
For each integer $m\ge0$, let $S_mf=\sum_{n=0}^{m}c_nT_n$ denote the
partial sum of degree $m$ of the Chebyshev series of $f$, where $T_n$ are
the Chebyshev polynomials and $c_n$ the Chebyshev coefficients of $f$. Then
\begin{equation}
\label{eq:tanh-chebyshev-tail-bounds}
\begin{aligned}
\|f-S_mf\|_{L^\infty(-1,1)}&\le \lambda_m(R),\\
\|f'-(S_mf)'\|_{L^\infty(-1,1)}&\le \mu_m(R),
\end{aligned}
\end{equation}
where, for $\rho>1$,
\[
\lambda_m(\rho)=\frac{2\rho^{-m}}{\rho-1},
\qquad
\mu_m(\rho)=2\sum_{n=m+1}^{\infty}n^2\rho^{-n}.
\]
\end{lemma}
\begin{proof}
The function $\tanh z$ is analytic for $|\operatorname{Im}z|<\pi/2$, and
for $z=x+iy$ with $|y|<\pi/2$,
\[
|\tanh z|^2
=\frac{\cosh2x-\cos2y}{\cosh2x+\cos2y}.
\]
If $|y|\le\pi/4$, then $\cos2y\ge0$, so $|\tanh z|\le1$.
The Bernstein ellipse with parameter $R$ has imaginary semiaxis
$(R-R^{-1})/2=\sinh(\log R)=\pi/(4W)$, so
$|\operatorname{Im}(b+Wt)|\le\pi/4$ for every $t$ in the closed ellipse.
Thus, $f$ is analytic in a neighborhood of the closed ellipse, and $|f|\le1$
on the closed ellipse.
The first estimate follows from \cite[Theorem~8.2]{T13}.
By \cite[Theorem~8.1]{T13}, the Chebyshev coefficients $c_n$ of $f$
satisfy $|c_n|\le2R^{-n}$ for $n\ge1$.
Markov's inequality gives
$\|T_n'\|_{L^\infty(-1,1)}\le n^2$
(cf.~\cite[Exercise~3.9]{T13}).
The differentiated series therefore converges uniformly, and summing
its tail gives the second estimate.
\end{proof}

\subsubsection*{Truncation of a patch network}
For a patch network from Section~\ref{subsec:subrectangle-local-interior-networks},
we suppress the component and parameter labels and write
\[
N_P(\boldsymbol\xi)
=c_P+\sum_{\ell=1}^{r}a_{P,\ell}
\tanh(\mathbf w_{P,\ell}^{\mathsf T}\boldsymbol\xi+b_{P,\ell}),
\]
with $W_{P,\ell}$ and $R_{P,\ell}$ as in
Section~\ref{subsec:degree-convergence}.
For each hidden unit with $W_{P,\ell}>0$, the normalized argument
$\mathbf w_{P,\ell}^{\mathsf T}\boldsymbol\xi/W_{P,\ell}$ lies in
$[-1,1]$ for $\boldsymbol\xi\in\widehat P$, and we write
$f_{P,\ell}(t)=\tanh(b_{P,\ell}+W_{P,\ell}t)$. Truncating these hidden units at
degree $k$ and keeping
the constant terms exact gives
\begin{equation}
\label{eq:patch-chebyshev-polynomial}
\begin{aligned}
N_{P,k}(\boldsymbol\xi)
&=c_P+\sum_{\ell:W_{P,\ell}=0}a_{P,\ell}\tanh b_{P,\ell}\\
&\quad+\sum_{\ell:W_{P,\ell}>0}a_{P,\ell}
(S_{k}f_{P,\ell})\!\left(
\frac{\mathbf w_{P,\ell}^{\mathsf T}\boldsymbol\xi}{W_{P,\ell}}
\right).
\end{aligned}
\end{equation}
The polynomial $N_{P,k}$ has total degree at most $k$, so it belongs to
${\mathbb Q}_{k}(\widehat P)$. Applying Lemma~\ref{lem:tanh-chebyshev-approximation}
unit by unit and using the chain rule gives
\begin{equation}
\label{eq:patch-chebyshev-error-bounds}
\begin{aligned}
\|N_P-N_{P,k}\|_{L^\infty(\widehat P)}
&\le\sum_{\ell:W_{P,\ell}>0}|a_{P,\ell}|\,\lambda_{k}(R_{P,\ell}),\\
\|\partial_{\xi_i}(N_P-N_{P,k})\|_{L^\infty(\widehat P)}
&\le\sum_{\ell:W_{P,\ell}>0}|a_{P,\ell}|
\frac{|w_{P,\ell,i}|}{W_{P,\ell}}\,\mu_{k}(R_{P,\ell}),
\qquad i=1,2.
\end{aligned}
\end{equation}

\subsubsection*{A comparison pair in $X_h^k$}
We apply \eqref{eq:patch-chebyshev-polynomial} to each component network,
using superscripts $u$, $p_x$, and $p_y$.
The pair ${\widetilde{\Phi}_h^k}=({v_h^k},{\mathbf z_h^k})$ is defined by
\begin{equation}
\label{eq:chebyshev-sdg-pair}
\begin{aligned}
{v_h^k}|_P
&=N^u_{P,k}\circ\boldsymbol\xi_P,
&&P\in\mathcal P_u,\\
({\mathbf z_h^k})_j|_P
&=N^{p_j}_{P,k}\circ\boldsymbol\xi_P,
&&P\in\mathcal P_{p_j},\quad j\in\{x,y\}.
\end{aligned}
\end{equation}
The scalar and vector components have degree at most $k$ in each
coordinate on every subrectangle.
The shared patch polynomials give scalar continuity across
$\mathcal F_u^0$ and normal-component continuity across $\mathcal F_p$.
Consequently, ${\widetilde{\Phi}_h^k}\in X_h^k$.

\subsubsection*{Patchwise error bounds}
For $s\in\{u,p_x,p_y\}$ and $P\in\mathcal P_s$, the value error bound
$V_P^s(k)$ is defined by the first line of
\eqref{eq:scalar-chebyshev-error-constants}. For $P\in\mathcal P_u$, the patch
widths $H_{P,i}=x_{P,i}^{\max}-x_{P,i}^{\min}$ give
$\partial_i\xi_{P,i}=2/H_{P,i}$, and the gradient error bound
$G_P(k)$ is defined by the second line:
\begin{equation}
\label{eq:scalar-chebyshev-error-constants}
\begin{aligned}
V_P^s(k)
&=\sum_{\ell:W^{s}_{P,\ell}>0}|a^{s}_{P,\ell}|\,
\lambda_{k}(R^{s}_{P,\ell}),\\
G_P(k)^2
&=\sum_{i=1}^2
\left[
\frac{2}{H_{P,i}}\sum_{\ell:W^u_{P,\ell}>0}|a^u_{P,\ell}|
\frac{|w^u_{P,\ell,i}|}{W^u_{P,\ell}}\,
\mu_{k}(R^u_{P,\ell})
\right]^2.
\end{aligned}
\end{equation}
Then \eqref{eq:patch-chebyshev-error-bounds} gives
\begin{equation}
\label{eq:scalar-chebyshev-patch-error}
\begin{aligned}
\|u_\theta-{v_h^k}\|_{L^\infty(P)}&\le V_P^u(k),\\
\|\nabla(u_\theta-{v_h^k})\|_{L^\infty(P;\mathbb R^2)}
&\le G_P(k),
\end{aligned}
\end{equation}
and, for $P\in\mathcal P_{p_j}$ with $j\in\{x,y\}$,
\begin{equation}
\label{eq:vector-chebyshev-patch-error}
\|(\mathbf p_\eta-{\mathbf z_h^k})_j\|_{L^\infty(P)}
\le V_P^{p_j}(k),
\end{equation}
where the vector-valued $L^\infty$ norm in
\eqref{eq:scalar-chebyshev-patch-error} uses the Euclidean norm on
$\mathbb R^2$.

\subsubsection*{The bound $\Lambda_h$}
For $\Phi=(u_\theta,\mathbf p_\eta)$, we combine the patch estimates
\eqref{eq:scalar-chebyshev-patch-error} and
\eqref{eq:vector-chebyshev-patch-error} in the norm $\|\cdot\|_{*,h}$.
For $e\in\mathcal F_p$, $P_e^\pm$ are the scalar patches containing
$\tau_e^\pm$.
The component label $j(e)$ is $x$ for a vertical edge and $y$ for a
horizontal edge, and $P_e^p\in\mathcal P_{p_{j(e)}}$ denotes the component
patch containing both $\tau_e^+$ and $\tau_e^-$.
For $e\in\mathcal F_\partial$, $P_e$ denotes the scalar patch
containing the subrectangle with side $e$.
On each $e\in\mathcal F_p$, the scalar jump error is bounded by
$V_{P_e^+}^u(k)+V_{P_e^-}^u(k)$, and the single-valued normal trace error
is bounded by $V_{P_e^p}^{p_{j(e)}}(k)$.
On each $e\in\mathcal F_\partial$, the trace $\gamma$ of
$u_\theta-{v_h^k}$ is bounded by $V_{P_e}^u(k)$. The quantities
$\Lambda_h^\partial(k;\theta,\eta)\ge0$ and $\Lambda_h(k;\theta,\eta)\ge0$ are
defined by
\[
\Lambda_h^\partial(k;\theta,\eta)^2
:=\sum_{e\in\mathcal F_\partial}
\bigl[V_{P_e}^u(k)\bigr]^2
\]
and
\begin{equation}
\label{eq:chebyshev-global-error-bound}
\begin{aligned}
\Lambda_h(k;\theta,\eta)^2
&=\sum_{P\in\mathcal P_u}|P|G_P(k)^2
+\sum_{e\in\mathcal F_p}
\bigl(V_{P_e^+}^u(k)+V_{P_e^-}^u(k)\bigr)^2\\
&\quad+\sum_{j\in\{x,y\}}\sum_{P\in\mathcal P_{p_j}}
|P|\bigl[V_P^{p_j}(k)\bigr]^2
+\sum_{e\in\mathcal F_p}h_e^2
\bigl[V_{P_e^p}^{p_{j(e)}}(k)\bigr]^2
+\Lambda_h^\partial(k;\theta,\eta)^2.
\end{aligned}
\end{equation}
Since $|e|=h_e$, the five terms in
\eqref{eq:chebyshev-global-error-bound} bound, in this order, the scalar
gradient, scalar jump, vector $L^2$, normal trace, and boundary trace
contributions to $\|\Phi-{\widetilde{\Phi}_h^k}\|_{*,h}^2$ in
\eqref{eq:projection-input-norm}, and the last contribution is
$|\gamma|_{\partial,h}^2$. Hence
\[
\|\Phi-{\widetilde{\Phi}_h^k}\|_{*,h}
\le{\Lambda_h(k;\theta,\eta)},
\qquad
{|\gamma|_{\partial,h}\le\Lambda_h^\partial(k;\theta,\eta)}.
\]
Since ${\widetilde{\Phi}_h^k}\in{X_h^k}$,
\eqref{eq:sdg-projection-boundary-split} with
$\Phi_h={\widetilde{\Phi}_h^k}$ and $\|\cdot\|_h\le\|\cdot\|_{*,h}$ yield
\begin{equation}
\label{eq:neural-sdg-chebyshev-estimate}
{P_h(k;\theta,\eta)}
\le K_0(k)\Lambda_h(k;\theta,\eta)+K_\partial(k)\Lambda_h^\partial(k;\theta,\eta).
\end{equation}

\begin{proof}[Proof of Theorem~\ref{thm:fixed-network-degree-convergence}]
Since $\operatorname{arsinh}$ is positive on
$(0,\infty)$, $R_*>1$. For $\rho>1$ and $m\ge0$, the definitions in
Lemma~\ref{lem:tanh-chebyshev-approximation} and the inequality
$n\le(m+1)(n-m)$ for $n\ge m+1$ give
\[
\lambda_m(\rho)=\lambda_0(\rho)\rho^{-m},
\qquad
\mu_m(\rho)\le\mu_0(\rho)(m+1)^2\rho^{-m}.
\]%
In the sums defining $V_P^u(k)$, $V_P^{p_j}(k)$, and $G_P(k)$, a hidden
unit with $a_{P,\ell}=0$ contributes no term, and every other hidden unit
satisfies $R_{P,\ell}\ge R_*$. Since the mesh and all network parameters are
fixed, these sums and the sums in \eqref{eq:chebyshev-global-error-bound} have
finitely many terms, and their number does not depend on $k$. The bounds above
and $(k+1)^2\le4k^2$ therefore give positive constants $C_1$ and $C_2$, independent of
$k$, such that
\[
\Lambda_h(k;\theta,\eta)\le C_1k^2R_*^{-k},
\qquad
{\Lambda_h^\partial(k;\theta,\eta)}\le C_2R_*^{-k},
\qquad k\ge1.
\]
By \eqref{eq:degree-projection-stability}, both terms in
\eqref{eq:neural-sdg-chebyshev-estimate} are at most constant
multiples of $k^6R_*^{-k}$, and \eqref{eq:fixed-network-degree-convergence}
holds with
\[
C=(1+C_AC_\beta)C_1+\bigl(C_\beta(6C_A+2)+6\bigr)C_2.
\]
Since $R_*>1$, $k^6R_*^{-k}\to0$.

If every hidden unit is constant, each patch network $N_P$ is constant
on $\widehat P$. Each component of $(u_\theta,\mathbf p_\eta)$ is then constant on
each of its patches, so $(u_\theta,\mathbf p_\eta)\in X_h^k$. Since
$\Pi_h^k$ reproduces $X_h^k$, it follows that
$P_h(k;\theta,\eta)=0$.
\end{proof}

\FloatBarrier

\section{\texorpdfstring{Quadrature determined by the network weights}{Quadrature determined by the network weights}}
\label{subsec:weight-quadrature}

Theorems~\ref{prop:scalar-loss-identity} and~\ref{prop:scalar-error-estimate}
and the estimates of Section~\ref{sec:projection-error-estimate} assume that the
integrals in \eqref{eq:scalar-residual-components} are evaluated exactly and that
the linear systems are solved exactly. In computations, we evaluate the
integrals that contain $u_\theta$ or $\mathbf p_\eta$ on each subrectangle $\tau$ by the tensor-product Gauss--Legendre rule
with $q_\tau$ points in each coordinate. An edge integral of a network
trace is evaluated by the Gauss--Legendre rule with $q_\tau$ points along the
edge, where $\tau$ is the subrectangle from which the trace is taken. The terms
without the networks, such as those with $f$ and $g$, are evaluated once with
fixed Gauss--Legendre rules. We choose
$q_\tau$ from the network weights, without evaluating residuals.

For $\rho>1$, the Bernstein ellipse with parameter $\rho$ is the open
region bounded by the ellipse with foci $\pm1$ whose semiaxes sum to $\rho$. A
point $z\in\mathbb C\setminus[-1,1]$ lies on the boundary of the Bernstein
ellipse with parameter $|z+\sqrt{z^2-1}|$, where the branch of the square root is
chosen so that this value exceeds one~\cite[Chapter~8]{T13}. Let
$\hat{\boldsymbol\xi}$ denote the coordinates of $\tau$ mapped to $[-1,1]^2$ by the
affine map of Section~\ref{subsec:subrectangle-local-interior-networks} with $\tau$
in place of $P$. In these coordinates, a hidden unit of a network whose patch
contains $\tau$ reads
$a_\ell\tanh(\hat{\mathbf w}_\ell^{\mathsf T}\hat{\boldsymbol\xi}+\hat b_\ell)$, where
$\hat{\mathbf w}_\ell$ and $\hat b_\ell$ are obtained from its weights by the affine
map between the coordinates of $\tau$ and those of its patch. As a function of
one coordinate with the other fixed, it is analytic near $[-1,1]$, and its
singularities nearest to $[-1,1]$ are the points where
$\hat{\mathbf w}_\ell^{\mathsf T}\hat{\boldsymbol\xi}+\hat b_\ell=\pm i\pi/2$. If $a_\ell=0$, or if the unit does not depend on that coordinate, there are no such points, and the unit or coordinate is skipped. Let
$\rho_\ell>1$ be the smallest Bernstein ellipse parameter of these points over
both coordinates and the values of the other coordinate in $[-1,1]$.

For an integrand analytic in the Bernstein ellipse with
parameter $\rho$ and bounded there by $M$, the $q$-point Gauss rule has error
at most $\tfrac{64}{15}M\rho^{-2(q-1)}/(\rho^2-1)$
\cite[Theorem~19.3]{T13}. We apply this bound to the unit times a test polynomial
of degree $k$ with maximum one on $[-1,1]$, with $\rho=\rho_\ell^{0.9}$,
which keeps the ellipse away from the singularities, and with
$M=|a_\ell|\rho^k$, which, by Bernstein's lemma for the test polynomial, would bound this product on the ellipse if
$|\tanh|$ did not exceed one there, and we take the
smallest $q$ for which it does not exceed a prescribed tolerance. The order $q_\tau$ is the
largest of these values over the units of the three networks whose
patches contain $\tau$, raised to $k+1$ if it is smaller and lowered to a
prescribed maximum order $q_{\max}$ if it is larger. Since $|\tanh|$ can exceed one on this
ellipse, $q_\tau$ is a predicted order rather than a guaranteed accuracy.
Section~\ref{subsec:numerical-settings} compares the computed losses with a
higher-order rule.

\bibliographystyle{amsplain-doi}
\bibliography{references}

\end{document}